\documentclass[11pt]{article}
\usepackage{amssymb,amsthm,amsmath,float}
 \usepackage{mathtools}   
\usepackage[total={6.5in,8.75in}, top=1.2in, left=0.9in, includefoot]{geometry}
\usepackage{graphicx}
\usepackage{algorithm}
\usepackage{algpseudocode}
\usepackage{color} 

\usepackage[pdftex,bookmarks, colorlinks, citecolor =blue, breaklinks]{hyperref}
\usepackage{url}

\newcommand{\Eq}[1]{(\ref{eq:#1})}
\newcommand{\Th}[1]{Th.~\ref{thm:#1}}

\newcommand{\Sec}[1]{\S \ref{sec:#1}}
\newcommand{\Fig}[1]{Fig.~\ref{fig:#1}}
\newcommand{\Tbl}[1]{Table~\ref{tbl:#1}}
\newcommand{\App}[1]{App.~\ref{app:#1}}
\newcommand{\Alg}[1]{Algorithm~\ref{alg:#1}}

\newcommand{\InsertFig}[4]
{\begin{figure}[h!t]
       \centerline{
         \includegraphics[width=#4]{./figures/#1}
       }
       \caption{{\footnotesize  #2}
       \label{fig:#3}}
\end{figure}}

\newcommand{\InsertFigTwo}[5] {
\begin{figure}[h!t]
       \centerline{
         \includegraphics[width=#5]{./figures/#1}
         \hskip 0.5in
         \includegraphics[width=#5]{./figures/#2}
       }
       \caption{{\footnotesize  #3}
       \label{fig:#4}}
\end{figure}}

 \newcommand{\InsertFigThree}[6] {
\begin{figure}[h!t]
       \centerline{
\renewcommand{\arraystretch}{0.01}
         \begin{tabular}{ccc}
         \includegraphics[width=#6]{./figures/#1}& \includegraphics[width=#6]{./figures/#2}
        \includegraphics[width=#6]{./figures/#3}  & 
        \end{tabular}
       }
       \caption{{\footnotesize  #4}
       \label{fig:#5}}
\end{figure}}

 \newcommand{\InsertFigFour}[7] {
\begin{figure}[h!t]
       \centerline{
\renewcommand{\arraystretch}{0.01}
         \begin{tabular}{cc}
         \includegraphics[width=#7]{./figures/#1}&  \includegraphics[width=#7]{./figures/#2} \\
        \includegraphics[width=#7]{./figures/#3}  &  \includegraphics[width=#7]{./figures/#4}
        \end{tabular}
       }
       \caption{{\footnotesize  #5}
       \label{fig:#6}}
\end{figure}}

\newcommand{\bN}{{\mathbb{ N}}}
\newcommand{\bQ}{{\mathbb{ Q}}}
\newcommand{\bR}{{\mathbb{ R}}}

\newcommand{\bT}{{\mathbb{ T}}}
\newcommand{\bZ}{{\mathbb{ Z}}}

\newcommand{\cD}{{\cal D}}

\newcommand{\cL}{{\cal L}}
\newcommand{\cO}{{\cal O}}
\newcommand{\cP}{{\cal P}}
\newcommand{\cR}{{\cal R}}

\newcommand{\cT}{{\cal T}}

\newcommand{\eps}{\varepsilon}

\newcommand{\WB} {\mathit{WB}}
\newcommand{\Znorm}[1] {{\| {#1} \|_\bZ}}
\newcommand{\digT}{{dig^{(T)}}}
\newcommand{\tv}[2] {{\begin{pmatrix} {#1}\\{#2} \end{pmatrix}}}

\newtheorem{thm}{Theorem}

\newcommand{\beq}[1]{\begin{equation}\label{eq:#1}}
\newcommand{\eeq}{\end{equation}}

\newenvironment{se}[1]{\equation\label{eq:#1}\aligned}{\endaligned\endequation}
\newcommand{\bsplit}[1]{\begin{se}{#1}}
\newcommand{\esplit}{\end{se}}

\newenvironment{example}[1][]
  {
	\setlength \leftmargini {1.0em}		
	\setlength \topsep {0.5em}			
	\begin{quote}
	{\it Example#1} }
	{\end{quote}
  }
\newcommand{\bexam}[1][:]{\begin{example}[#1]}
\newcommand{\eexam}{\end{example}}

\title{Birkhoff Averages and the Breakdown of Invariant Tori in Volume-Preserving Maps}
\author{J.D.~Meiss and E. Sander \thanks 
      {
        JDM was supported in part by NSF grant DMS-181248. ES was supported in part 
        by the Simons Foundation under Award 636383. ES and JDM acknowledge
        support from NSF grant DMS-140140 while they were at residence at the 
        Mathematical Sciences Research Institute in Berkeley, CA, during the Fall 2018 semester.
        Useful conversations with Keith Briggs--who pointed out \cite{Clarkson97a}---and with Robert MacKay are gratefully acknowledged. 
      }
    \\
 \begin{tabular}{ll}
	Department of Applied Mathematics 		&Department of Mathematical Sciences	\\
    University of Colorado					&George Mason University				\\
	Boulder, CO 80309-0526, USA 					&Fairfax, VA 22030, USA			\\
	James.Meiss@colorado.edu				&esander@gmu.edu	\\ 
\end{tabular}
}
\date{\today}
\begin{document}
\maketitle
\tableofcontents      
\newpage			  

\begin{abstract}
\vspace*{1ex}
\noindent

In this paper, we develop numerical methods 
based on the weighted Birkhoff average for
studying two-dimensional invariant  tori for volume-preserving maps.
The methods do not rely on symmetries, such as time-reversal
symmetry, nor on approximating tori by periodic orbits. 
The rate of convergence of the average gives a sharp distinction 
between chaotic and regular dynamics and allows accurate computation of rotation vectors
for regular orbits. Resonant and rotational tori are distinguished by computing
the resonance order of the rotation vector to a given precision.
Critical parameter values, where tori are destroyed, are computed by a
sharp decrease in convergence rate of the Birkhoff average. 
We apply these methods for a three-dimensional generalization of Chirikov's standard map: 
an angle-action map with two angle variables. Computations on grids
in frequency and perturbation amplitude allow estimates of the
critical set. We also use continuation to follow tori with fixed rotation vectors.
We test three conjectures for cubic fields that have been proposed to give locally robust invariant tori.

\end{abstract}

\section{Introduction}

The dynamics of an integrable Hamiltonian or volume-preserving system
consists of quasi-periodic motion on invariant tori. As such a system is
smoothly perturbed, KAM theory implies that some of these tori persist, but some are
replaced by isolated periodic orbits, resonances, and chaotic regions.
Typically, as the perturbation grows, more of the tori are destroyed.
For two-dimensional maps, the robust tori are circles on which the
dynamics is conjugate to rigid rotation with a Diophantine rotation number.
It is conjectured from careful numerical study that the most robust of these invariant circles
have rotation numbers that are ``noble"---they are in the quadratic
field of the golden mean, or equivalently they have continued fractions
with an infinite tail of ones \cite{MacKay83,MacKay92c}. A similar result
for higher dimensional tori has not been found even though,
as we discuss below, there has been considerable research and conjecture
on a suitable generalization.

Here we investigate the existence of tori for a map $f: \bT^d \times \bR^k \to \bT^d \times \bR^k$ of the angle-action form
\bsplit{AngleActionMap}
	x' = x + \Omega(y') \\
	y' = y + \eps F(x) . 
\esplit
We view $x \in \bT^d = \bR^d/\bZ^d$ as angle variables, taken
modulo one, and $y \in \bR^k$ as action variables. The function $\Omega
: \bR^k \to \bT^d$ is the \textit{frequency map} and $F : \bT^d \to
\bR^k$ is the \textit{force}. This family of maps is the composition of two
volume-preserving shears, e.g., $(x,y) \mapsto (x+\Omega(y),y)$ and
$(x,y) \mapsto (x,y+ \eps F(x))$, and hence is always volume-preserving. If $k = d$,
\Eq{AngleActionMap} is symplectic if the
force and frequency maps are gradients: $F(x) = -\nabla V(x)$,
$\Omega(y) = \nabla S(y)$. Two prominent and well-studied examples of such maps are
Chirikov's standard area-preserving map \cite{Chirikov79a}, and Froeschl\'e's
four-dimensional symplectic map \cite{Froeschle73}.

When $\eps = 0$ the dynamics of \Eq{AngleActionMap} is simple: the actions are constant, and every
orbit lies on a ``horizontal" $d$-torus
\beq{HorizontalTorus}
	H(y) = \{(x,y) : x \in \bT^d\} .
\eeq
When $\eps =0$, the dynamics of $f|_H$ is simply horizontal translation by $\omega = \Omega(y)$, i.e.,
every orbit on $H$ has rotation vector $\omega$.
More generally an orbit $\{(x_t,y_t): t \in \bZ \}$ has rotation vector $\omega$ if the limit
\beq{rotNum}
	\omega = \lim_{T \to \infty} \frac{1}{T} \sum_{t=0}^{T-1} \Omega(y_t)
\eeq
exists. Of course, if $\eps = 0$, then this is simply the value of the frequency map on the conserved action. 

We say that a $d$-dimensional torus $\cT$ is \textit{rotational} if it is homotopic to $H(0)$. If, in addition,
the torus is invariant under \Eq{AngleActionMap} and $f|_{\cT}$ is conjugate to rigid translation with a rotation vector $\omega$ we will denote the torus by $\cT_\omega$. 

When the force $F$ and frequency map $\Omega$ are analytic and a twist condition is satisfied, KAM theory shows that there are tori with ``Diophantine'' rotation vectors that persist when $\eps$ is nonzero but small \cite{Cheng90b,Xia92}.
A vector $\omega $ is defined to be Diophantine, denoted $\omega \in \cD$, where 
\beq{Diophantine}
	\cD = \bigcup_{c>0} \left\{\omega : \quad |m \cdot \omega -n | > \frac{c}{\|m\|_\infty^d} ,
	\quad
	\forall (m,n) \in \bZ^d\setminus \{0\} \times \bZ \right\}.
\eeq
By contrast, we say that $\omega$ is \textit{resonant} if there exists a nonzero $m \in \bZ^d$ such that
\beq{Resonance}
	m \cdot \omega  = n \in \bZ.
\eeq
For area-preserving maps, i.e., \Eq{AngleActionMap} with $d=k=1$, resonances correspond to periodic orbits 
where $\omega = p/q$ is rational. 
Elliptic periodic orbits are typically surrounded by island chains, 
and hyperbolic orbits have stable and unstable manifolds that typically
intersect transversely, giving rise to chaotic motion. For $d>1$ island chains are replaced by resonant tubes, 
and these are also surrounded by chaotic zones. Typically as the parameter $\eps$ grows, 
so do the regions of chaos and resonance, destroying more of the rotational tori. Our goal in this paper is to use computations of \Eq{rotNum} to investigate this destruction.

Our major tool is the \textit{weighted Birkhoff average} that was introduced in \cite{Das16a,Das16b,Das17}.
We will use this to compute the rotation vector 
\Eq{rotNum} efficiently and accurately. The rigorous convergence
results given in \cite{Das18} imply that the same method can be used to distinguish and remove
chaotic orbits, leaving only regular behavior.  We previously used this method to compute rotational 
circles for the Chirikov standard map and several other 2D  maps \cite{Sander20}. 

For two-dimensional maps, a computed rotation number has been used to efficiently find transport barriers 
\cite{Szezech13} and to find the breakup of circles in Chirikov's standard map \cite{Abud15} 
and in nontwist maps \cite{Santos18}. 
The gradient of $\omega$ was also used as an indicator of stickiness \cite{Santos19}.
A number of methods have been proposed for computing a scalar rotation number accurately, 
based on recurrence \cite{Efstathiou01}, conjugacy to rigid rotation on a circle \cite{Seara06, Luque09},
or recurrence times using Slater's method \cite{Mayer88, Altmann06, Zou07}.

There have been a number of studies of the existence and breakup of
tori for angle-action maps. A major focus of these studies is to attempt
to identify the subset of the Diophantine frequency vectors for which
the invariant tori are \textit{locally robust}; that is, more resistant
to destruction than nearby vectors. The noble numbers that are robust
for the area-preserving case are quadratic irrationals --- in the field
$\bQ[\phi]$ of the golden mean --- and it is known that more generally the
class of algebraic numbers contains Diophantine vectors
\cite{Cassels57}. As a result,  it has long been thought that a robust $d$-torus would
typically have a frequency vector that is formed from a basis for a
degree-$d$ algebraic field.

We will study the case $d=2$ and $k=1$ where the expected robust
rotation vectors are cubic irrationals. One can classify the cubic
fields by the discriminant of the minimal generating polynomial. There
are---at least---three natural conjectures about which of these fields
should replace $\bQ[\phi]$. For example, Hu and Mao \cite{Hu88} studied
a map on $\bT^2$, the case $d=2$ and $k=0$, looking at tori in the cubic
field with discriminant $D = -44$. This field is ``natural'' from the
point of view of a generalization of the continued fraction, the
Jacobi-Perron algorithm (JPA). Just like the golden-mean has a continued
fraction with elements all equal to one, there is a basis for the
$D=-44$ field with a period-one JPA expansion consisting of ``all
ones''.  The robustness of tori with frequency vectors in this same
field were also studied in \cite{Tompaidis96b}. Both of these studies
used periodic orbits to approximate the tori.

A higher dimensional case, $d=2$, $k=1$, was studied by Artuso et al \cite{Artuso92}. They fixed the frequency map to be $\Omega = (y,\delta)$ so that only the first component depends upon the action variable. Given an irrational value for $\delta$, this map is a quasiperiodically forced area-preserving map. In this paper frequencies from another cubic field, that with $D = -23$, were studied. This field corresponds to the so-called spiral mean proposed by \cite{Kim86} when they developed a generalization of the Farey, or Stern-Brocot, tree expansion. The spiral mean is distinguished by its period-one expansion (that spirals) on this tree. Artuso again used periodic orbits to approximate the incommensurate frequency vector, and generalized \textit{Greene's residue criterion} \cite{Greene79} for this case. The residue criterion essentially conjectures that a torus exists only if sequences of periodic orbits that converge to it have linearizations with bounded eigenvalues, as measured by the trace, or a scaled version of this that Greene called the residue.

Greene's residue was also used in \cite{Fox13} to study a fully three-dimensional case, again looking for breakup of tori by studying sequences of periodic orbits that approach given incommensurate vectors on the generalized Farey tree. This map will be the focus of the paper below, see \Sec{vpm}. Later, Fox and Meiss \cite{Fox16} computed tori directly from their conjugacy to a rigid rotation, using the efficient, parameterization method \cite{Haro16} to compute Fourier series.

Studies of the breakup of tori for four-dimensional, symplectic maps of the form \Eq{AngleActionMap} include \cite{Kook89b} who computed periodic orbits and the frequency map on the Kim-Oslund tree for the Froeschl\'e map. Later \cite{Bollt93} studied a complex extension of this map and tori in the spiral mean field as well as several quartic irrational vectors. Attempts to extend Greene's residue method to the 4D case include \cite{Tompaidis96b, Vrahatis96, Zhou01b}, though to our knowledge, no one has found a generalization of the renormalization, or self-similarity property that is observed in the 2D case.  

A third cubic field that has been proposed to replace $\bQ[\phi]$ corresponds to the cubic $D=49$ field \cite{Lochak92}. This field has the smallest discriminant of all the totally real fields, and is conjectured to have bases with the largest value of the (linear) Diophantine constant (see \Sec{CriticalSets}) among all vectors for $d=2$ \cite{Cusick74}. Lochak argued that the linear approximation constant is more appropriate from the point of view of KAM theory, than the simultaneous constant, as these numbers appear in the small denominators in the Fourier series expansions for the conjugacy functions of tori. The maximal Diophantine property, of course, would generalize the similar, proven property of the golden mean for $d=1$. 

The rest of this paper proceeds as follows. 
In \Sec{vpm} we introduce the standard three-dimensional, volume-preserving model that
we study in this paper. 
Section~\ref{sec:wba} describes the weighted Birkhoff average. 
In \Sec{Tori} we describe methods for distinguishing regular behavior from 
chaotic dynamics, and for distinguishing resonant from rotational tori. 
In \Sec{CriticalSets} we consider locally robust tori and the critical
surface and in \Sec{Continuation} study the continuation of tori with rotation vectors in cubic algebraic fields.
We conclude in \Sec{Conclusions} and describe some of the many problems that remain open.

\section{Standard Volume-Preserving Map}\label{sec:vpm}
A three-dimensional analog to Chirikov's area-preserving map and Froeshl\'e's four-dimensional symplectic map was obtained in \cite{Dullin12a}. This normal form corresponds to  \Eq{AngleActionMap} with $(x,y) \in \bR^2 \times \bR^1$ and the frequency map and force
\bsplit{Parabola}
	\Omega(y,\delta) &= (y+\gamma, -\delta + \beta y^2)\;, \\
	F(x) &= -a\sin(2\pi x_1) - b \sin(2 \pi x_2) -c \sin(2\pi(x_1-x_2)) \;.
\esplit
We will think of five of the parameters as fixed, choosing
\beq{StdParams}
	\gamma = \tfrac12(\sqrt{5}-1) \approx 0.61803 \;, \quad
	\beta = 2 \;, \quad
	a = b = c = 1 \;.
\eeq
This leaves two essential parameters, $\delta$ and $\eps$, that will vary for our computations. 
Note that for each $\delta$, the image $\Omega(y,\delta)$ is a parabola in $\bR^2$: only invariant tori with rotation vectors that lie on this curve exist in the integrable case $\eps = 0$. However, we take $\delta$
to be an essential parameter. Allowing it to vary makes the frequency map $\Omega: \bR^2 \to \bR^2$ a diffeomorphism. 

More generally, suppose that the initial point $(x,y) = (0,y_0)$ lies on a rotational, invariant two-torus with rotation vector $\omega$. We call such a torus $\cT_\omega(\eps,y_0,\delta)$, labeling it with parameters 
$\eps$ and $\delta$ as well as the initial action. Note that $\cT_\omega(0,y_0,\delta) = H(y_0$) 
and that a Cantor set of Diophantine tori are preserved when $\eps \ll 1$, 
according to the volume-preserving version of KAM theory \cite{Cheng90b, Xia92}.

Previous computational studies of invariant
tori for this map include studies of ``crossing orbits" giving parameter thresholds for the ``last torus'' that
divides vertically separated points \cite{Meiss12a}, a version of Greene's residue criterion to find critical tori with given rotation vectors---tori at the threshold of destruction \cite{Fox13}, and the parameterization method to
numerically compute tori and their breakup thresholds \cite{Fox16}.

\InsertFigFour{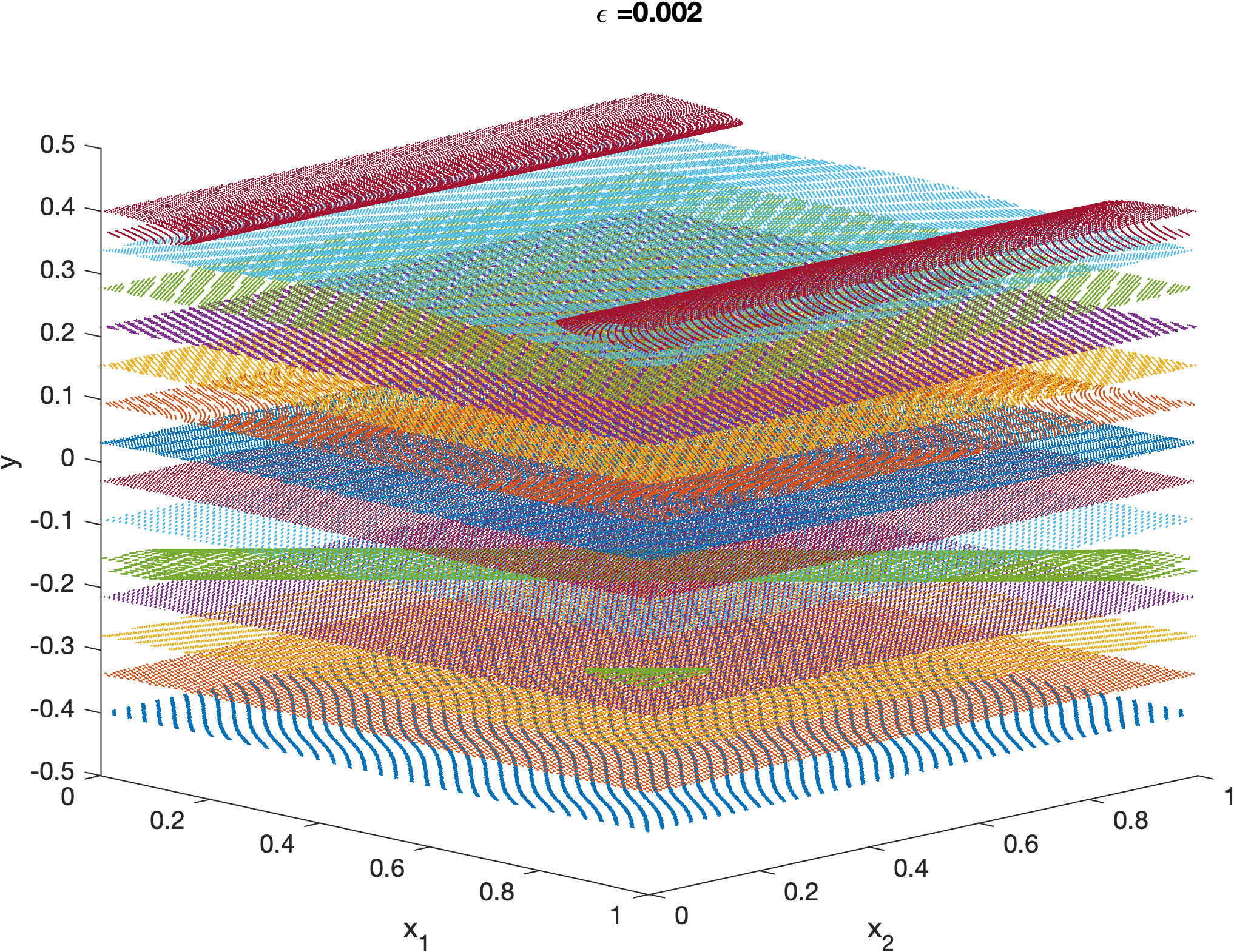}{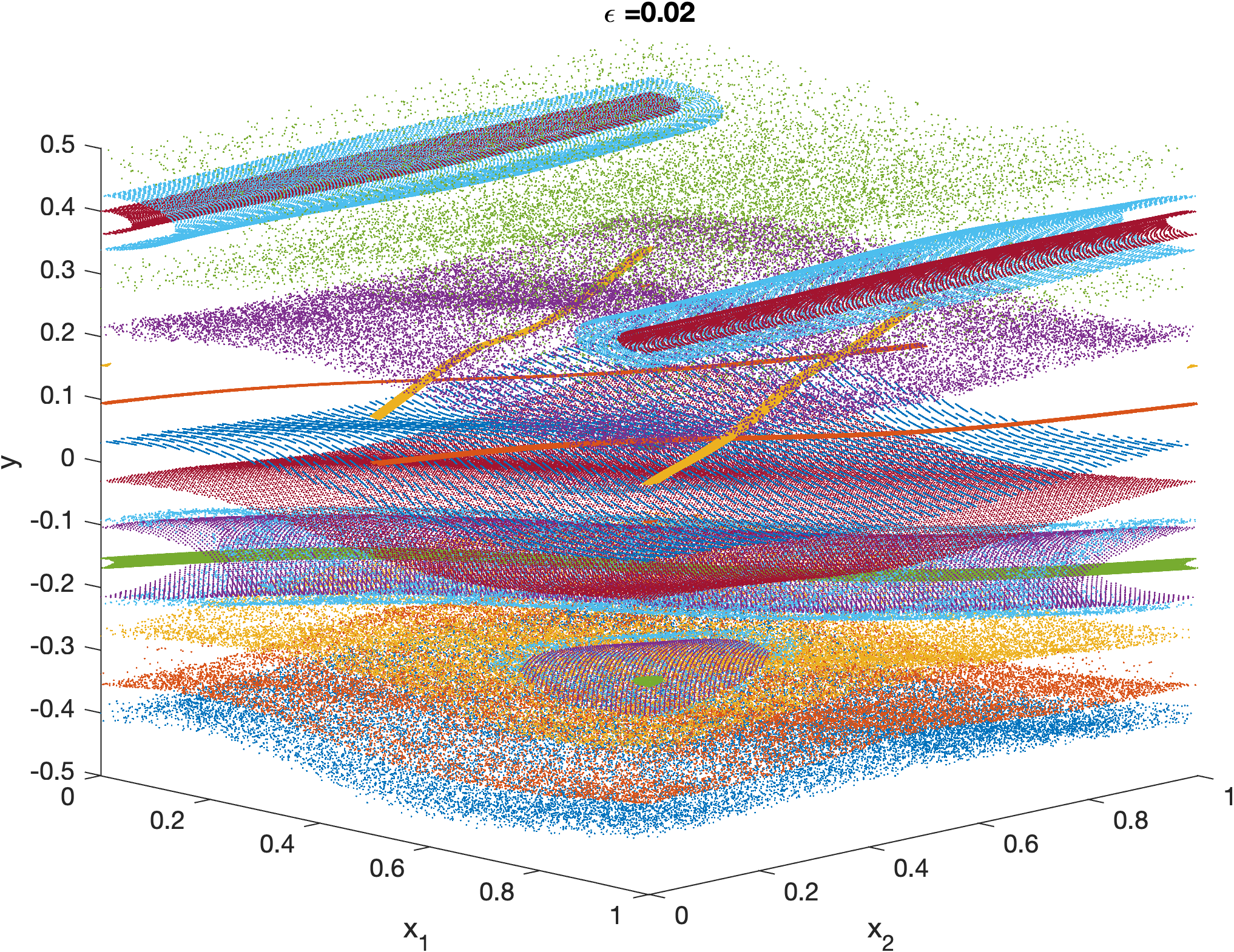}{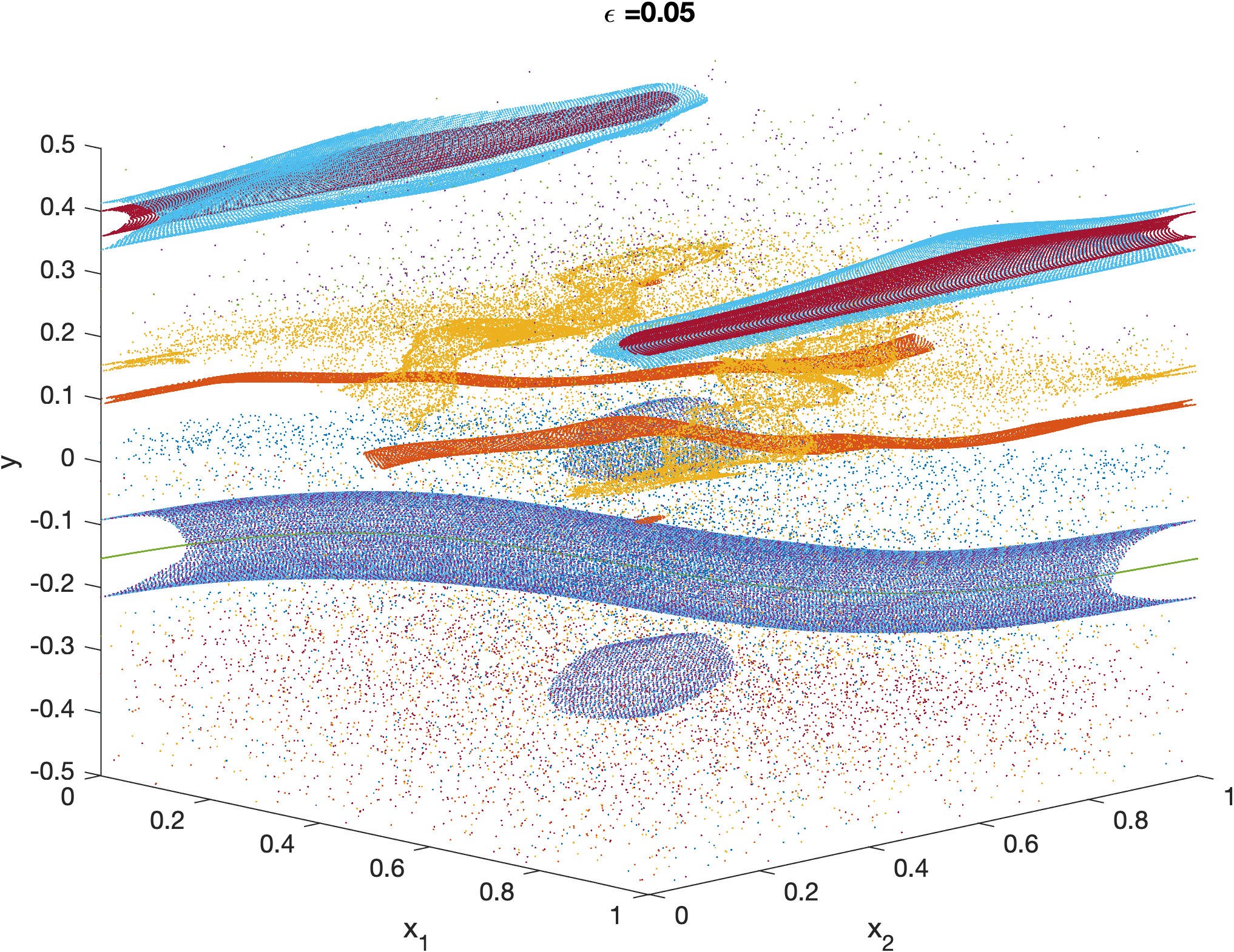}{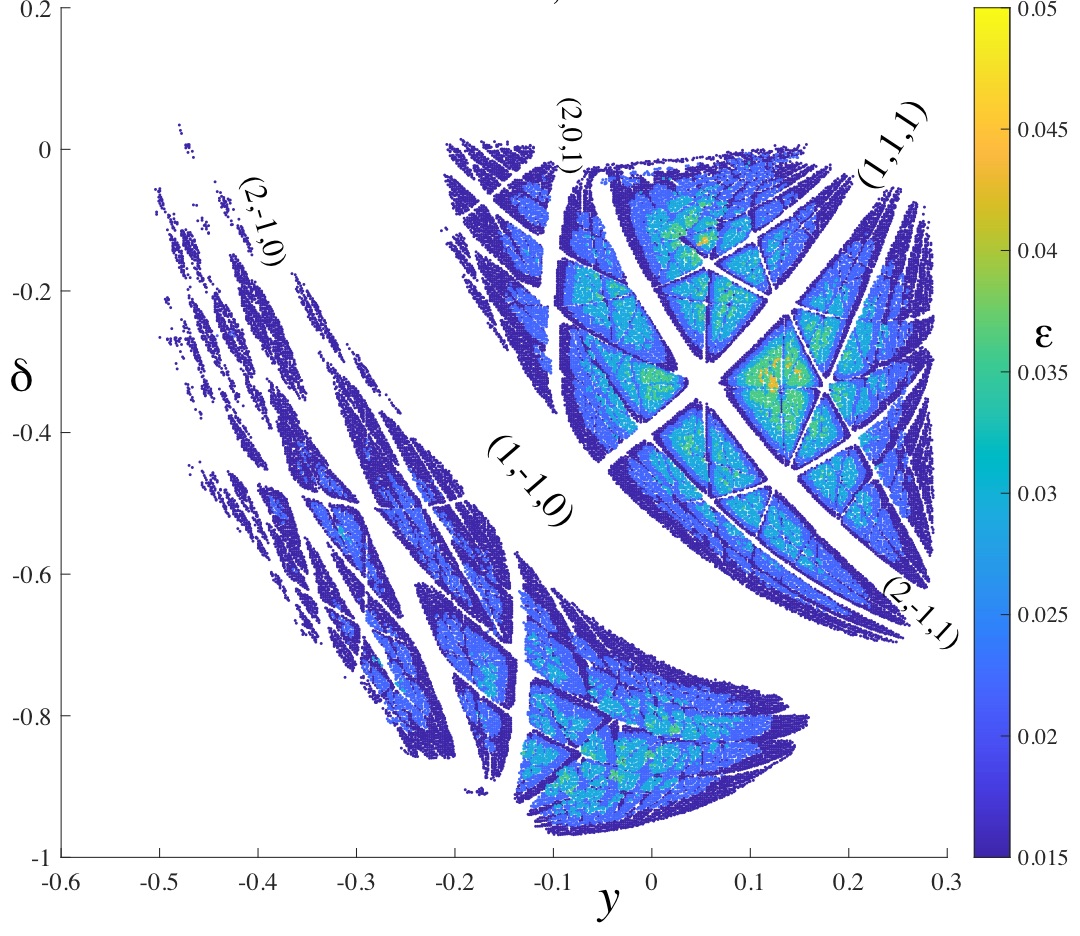}{
(a-c) Orbits for~\Eq{AngleActionMap} and \Eq{Parabola} 
for $\delta = -0.4$ and three values of $\eps$. 
All initial conditions have $x_1=0, x_2=0$. Each image shows 
fourteen different orbits with $y_0 \in [-0.4,0.5]$. As $\eps$ increases, the number of
rotational tori decreases, and at $\eps = 0.05$, only resonant
tori and chaotic orbits are visible in the figure. 
(d) Values of $(y_0,\delta)$ for which there are rotational tori with $\omega \in [0,1]^2$. 
The color represents the largest value of $\eps$ with a corresponding rotational torus. 
}{zdelta}{3in}

The first three panels of \Fig{zdelta} show examples of orbits for
\Eq{AngleActionMap} with \Eq{Parabola} for three values of $\eps$, and
$\delta=-0.4$ (other values of $\delta$ exhibit similar behavior). As
predicted by KAM theory, when $\eps \ll 1$ the typical orbits appear to
be dense on (rotational) two-tori $\cT_\omega$ that are graphs over the
angles $x$ with $y$ nearly constant. 
Even in \Fig{zdelta}{(a)}, however one can see two resonant
tubes. These are driven by the primary resonances, \Eq{Resonance}, of
the force $F$, and correspond to the phases $x_1$ and $x_1-x_2$
remaining constant mod 1 (the third driven resonance, where $x_2$ is
constant is out of range of the figure). 
These resonant tubes correspond to the resonances in \Tbl{LowOrderResonance}  such that 
$(m,n) = (1,0,1)$ and $(1,-1,0)$ respectively.
As $\eps$ grows, more of the orbits become chaotic and other resonances become visible as tube-like structures, in
particular the resonant tubes with  $(m,n) = (2,-1,1)$ and $(2,1,2)$  seen in \Fig{zdelta}(b,c).

\begin{table}[htp]
\begin{center}
\begin{tabular}{c|c||c|c}
$(m,n)$   & $y$  &$(m,n)$   & $y$ \\
\hline   
(1,1,1)   & -0.481 &(2,-1,0)  & -0.317  \\
(1,-1,0)  & -0.164 &(0,2,1)   & -0.224  \\
(1,1,1)   & -0.019 &(2,0,1)   & -0.118  \\
(1,0,1)   & 0.382  &(2,-1,1)  & 0.090   \\
		&		   &(2,1,2)   & 0.157   \\
		&		   &(0,2,1)   & 0.224   \\
		&		   &(1,2,2)   & 0.276   \\
		&		   & (1,1,2)   & 0.494  
\end{tabular}
\end{center}
\caption{Action values for resonances up to order two in the range $|y| < 0.5$ when $\eps = 0$ and $\delta = -0.4$.}
\label{tbl:LowOrderResonance}
\end{table}%

We restrict our interest to tori with rotation vectors \Eq{rotNum},
in a fixed range, $\omega \in [0,1]\times [0,1]$. When $\eps = 0$, 
\Eq{Parabola} implies that $\omega$ depends linearly on $\delta$ and
quadratically on $y_0$. Indeed, each resonance in \Eq{Module} defines
a parabola (or a vertical line if $m_2 = 0$), $m \cdot
\Omega(y_0,\delta) = n$, in the $(y_0,\delta)$ plane. When $\eps$ is
relatively small, we expect that persisting rotational tori will have
rotation vectors that at least approximate this quadratic relationship.
To illustrate this we compute the rotation vector $\omega$  using the
methods described in \Sec{wba} and \Sec{Tori} below. \Fig{zdelta}{(d)}
shows the values of $(y_0,\delta)$ for which there are rotational tori
with $\omega \in [0,1]^2$. The color represents the largest $\eps \in
[0.015, 0.045]$ for which a rotational torus exists for a given
$(y_0,\delta)$. The parabolic relationship between $y_0$ and $\delta$ is
still clear in this image---the gaps represent initial conditions for
which the corresponding orbits are resonant or chaotic. Several of the
low order resonances are labeled in the figure. As $\eps$ grows there are fewer
initial conditions that lie on rotational tori, indicated by the dearth
of yellow in \Fig{zdelta}{(d)}.  

Our primary goal is to study the persistence of rotational tori of \Eq{AngleActionMap} with conditions
\Eq{Parabola} as $\eps$ grows from $0$.
We observe in \Sec{Tori} that, for $\omega \in [0,1]^2$,
there are rotational tori only when $\eps < 0.051$.
In \Sec{CriticalSets} we compute the most
robust torus in subsets of this $\omega$ region;
that is, the rotational torus with the largest maximum $\eps$ value in the subset. 

In our calculations, we use the $\eps = 0$ approximation to determine appropriate ranges for $y_0$ and $\delta$, setting $(y_0,\delta) = \Omega^{-1}(p_1,p_2)$, i.e., inverting the frequency map \Eq{Parabola}, to obtain
\beq{ydeltaGrid}
	(y_0,\delta) \in \cP = \{ (p_1-\gamma, \beta(p_1-\gamma)^2-p_2) : -0.05 \le p_1,p_2 \le 1.05 \}.
\eeq
The added $0.05$ in $\cP$ is a buffer to cover all $\omega \in [0,1]^2$ as $\eps$ grows.
Our calculations indicate this  buffer is sufficient; indeed, we have checked that values of $(y_0,\delta)$ outside this range do not give such $\omega$ values. 

\section{Weighted Birkhoff Averages}\label{sec:wba}

A finite-time Birkhoff average on an orbit of a map $f:M \to M$ beginning at a point 
$z \in M$  for any function  $h: M \to \bR$ is the sum
\beq{Birkhoff}
	B_T(h)(z) = \frac{1}{T} \sum_{t=0}^{T-1} h \circ f^t(z) .
\eeq
This average need not converge rapidly. Even if the orbit lies on a smooth
invariant torus with irrational rotation vector, the convergence rate of \Eq{Birkhoff}
is $\cO(T^{-1})$, caused by edge effects for the finite orbit segment.
By contrast, for the chaotic case, the convergence rate of \Eq{Birkhoff} is observed
to be $\cO(T^{-1/2})$, in essence as implied by the central limit theorem~\cite{Levnajic10}. 

The convergence of \Eq{Birkhoff} on a quasiperiodic set can be significantly improved by 
using the method of weighted Birkhoff averages developed in~\cite{Das16a, Das16b, Das17}.
Since the source of error in the calculation of a time average for a quasiperiodic set 
is due to the lack of smoothness at the  ends of the orbit, we use a windowing method 
similar to the methods used in signal processing. 
Let 
\[
	g(t) \equiv \left\{ \begin{array}{ll}  e^{-[t(1-t)]^{-1}}  & t \in (0,1) \\
	         								0	& t \le 0 \mbox{ or } t \ge 1 
						\end{array} \right.\;, 
\]
be an exponential bump function that converges to zero with infinite smoothness 
at $0$ and $1$, i.e.,
 $g^{(k)}(0) = g^{(k)}(1) = 0$ for all $k \in \bN$. 
To estimate the Birkhoff average of a function $h: M \to \bR$  efficiently and accurately 
for a length $T$ segment of an orbit, we modify \Eq{Birkhoff} to compute 
\beq{WB}
	\WB_{T}(h)(z) = \sum_{t = 0}^{T-1} w_{t,T} \,\, h \circ f^t(z) \;, 
\eeq
where 
\bsplit{SmoothedAve}
    w_{t,T} &= \frac{1}{S} g\left(\tfrac{t}{T}\right) ,
	&S &= \sum_{t=0}^{T-1} g\left(\tfrac{t}{T}\right)  \;. 
\esplit
That is, the weights $w$ are chosen to be normalized and evenly spaced 
values along the curve $g(t)$. For a quasiperiodic orbit, the 
infinitely smooth convergence of $g$ to zero at the edges of the definition interval
preserves the smoothness of the original orbit. 
Indeed it was shown in \cite{Das18} that for a $C^\infty$ map $f$, a quasiperiodic orbit $\{f^t(z)\}$ 
with Diophantine rotation number, and a $C^\infty$ function $h$, it follows that \Eq{WB} is super-convergent: there are constants $c_n$, such that for all $n \in \bN$
\beq{WBerror}
	\left |\WB_{T}(h)(z) - \lim_{N \to \infty}B_N(h)(z) \right| <  c_n T^{-n} .
\eeq
Several papers~\cite{Gomez10a, Laskar92, Luque14} include a similar method to compute frequencies 
with a $\sin^2(\pi s)$ function instead of  a bump function, but this function is fourth order 
smooth rather than infinitely smooth at the two ends, implying that the method converges as 
$\cO(T^{-4})$, see e.g., \cite[Fig. 7]{Das17}. In addition to converging more rapidly, the weighted Birkhoff 
average \Eq{WB} is relatively straightforward to implement. 

\section{Computing Tori}\label{sec:Tori}
Using the weighted Birkhoff average \Eq{WB}, we can compute an approximation to a rotation vector 
$\omega$ as $\WB_T(\Omega)$ for the frequency map \Eq{Parabola}. To discern whether an orbit lies on a
rotational torus $T_\omega$  we must first distinguish chaotic from 
regular orbits, and then distinguish resonances from nonresonant tori.

By contrast to the case of regular orbits, when an orbit is chaotic (i.e., has positive Lyapunov exponents),
\Eq{WB} typically converges much more slowly; in general it converges no more rapidly than the unweighted 
average of a random signal, i.e., with an error $\cO(T^{-1/2})$ \cite{Levnajic10, Das17}. We see
in \Sec{Chaos} that this distinction is valid for the map \Eq{AngleActionMap} as well.

Given a regular orbit, in \Sec{Resonances}, we use resulting the high-precision computation 
$\omega \simeq \WB_T(\Omega)$ to define an approximate resonance order.
This allows a distinction, up to some precision, between those orbits that have a commensurate frequency vector and those that appear to be nonresonant.

\subsection{Distinguishing Chaos}\label{sec:Chaos}

To establish the distinction between chaotic and regular orbits, we estimate
the number of digits of accuracy in the weighted Birkhoff average. Following \cite{Sander20},
we compute \Eq{WB} for two segments of 
an orbit, using iterates $\{ 1,\dots,T\}$ and $\{ T+1,\dots, 2T \}$; 
these values should be approximately equal when $T$ is large since the Birkhoff average depends only on the choice of orbit. A comparison  of these gives the error estimate
\beq{digits}
	dig^{(T)}_h = -\log_{10} \left|\WB_{T}(h)(z)-\WB_{T}(h)(f^T(z)) \right| ,
\eeq
i.e., the number of consistent digits beyond the decimal point for the  two approximations of $\WB(h)$.
If, for a modest value of $T$, $dig^{(T)}_h$ is relatively large, then the convergence is relatively rapid, meaning the  orbit is regular.  On the other hand, if $dig^{(T)}_h$ is small, then convergence is slow, with the implication that the orbit is chaotic. 

Using $h = \Omega$, \Eq{Parabola}, the accuracy of the calculation of both components of $\omega$
is 
\beq{digi}
	\digT = \min \{ dig^{(T)}_{\Omega_1}, dig^{(T)}_{\Omega_2} \} ,
\eeq
In addition to distinguishing chaotic and regular orbits, $\digT$ can be  
used to estimate the precision of $\omega$. We will use these precision estimates in \Sec{CriticalSets} when 
we consider number theoretic properties of the frequencies of the  robust tori. 

\Fig{digitsviterates} shows the behavior of orbits 
when $\delta=-0.4$ and $\eps = 0.02$ for a set of initial conditions along the line $x_1 = x_2 = 0$.
Panel (a) is the slice with $|x_2| \le 0.005$ through \Fig{zdelta}(b); it clearly shows a strongly chaotic region for 
$y \lesssim -0.41$, then narrower chaotic bands for $-0.41 \lesssim y \lesssim -0.07$, followed by a region of tori and resonances up to $y \approx 0.27$, and finally a mixed regular/chaotic region up to $y = 0.5$. The largest resonant regions shown correspond to $m = (1,-1)$ near $y = -0.16$, and $m = (1,0)$ near $y = 0.38$; we also saw these in \Fig{zdelta} and \Tbl{LowOrderResonance}. 

\Fig{digitsviterates}(b) shows corresponding values of $\digT$ for three values of $T$. 
Even when $T = 10^3$ (blue points), the calculations can distinguish
between strongly chaotic---where $\digT \sim 2-4$---and regular
orbits---where $\digT \sim 8$. However there is a population of orbits,
especially those near the edges of resonant regions, that have
intermediate values of $\digT$, and for these it is harder to obtain a
definitive classification. When $T=3(10)^4$ (orange points), the 
values of $\digT$ become better separated---with regular orbits having $\digT \sim 13-14$ and chaotic 
orbits still having $\digT \sim 2-4$ ---but 
there are still 
a fair number of intermediate values, again especially at the region edges. 
However, when $T = 10^6$ (red points),
regular orbits predominantly achieve full floating point accuracy,
$\digT \sim 14-16$, while chaotic orbits still have $\digT \lesssim 5$. 

\Fig{digitsviterates}(c) shows the components of $\omega$ computed for $T=10^6$. Each
of the flat intervals corresponds to a range of $y$ values in a resonant
tube, and these are bounded by chaotic regions where the computed values
of $\omega$ vary rapidly with $y$ and the corresponding accuracy is low.

\InsertFigThree{orbitslice}{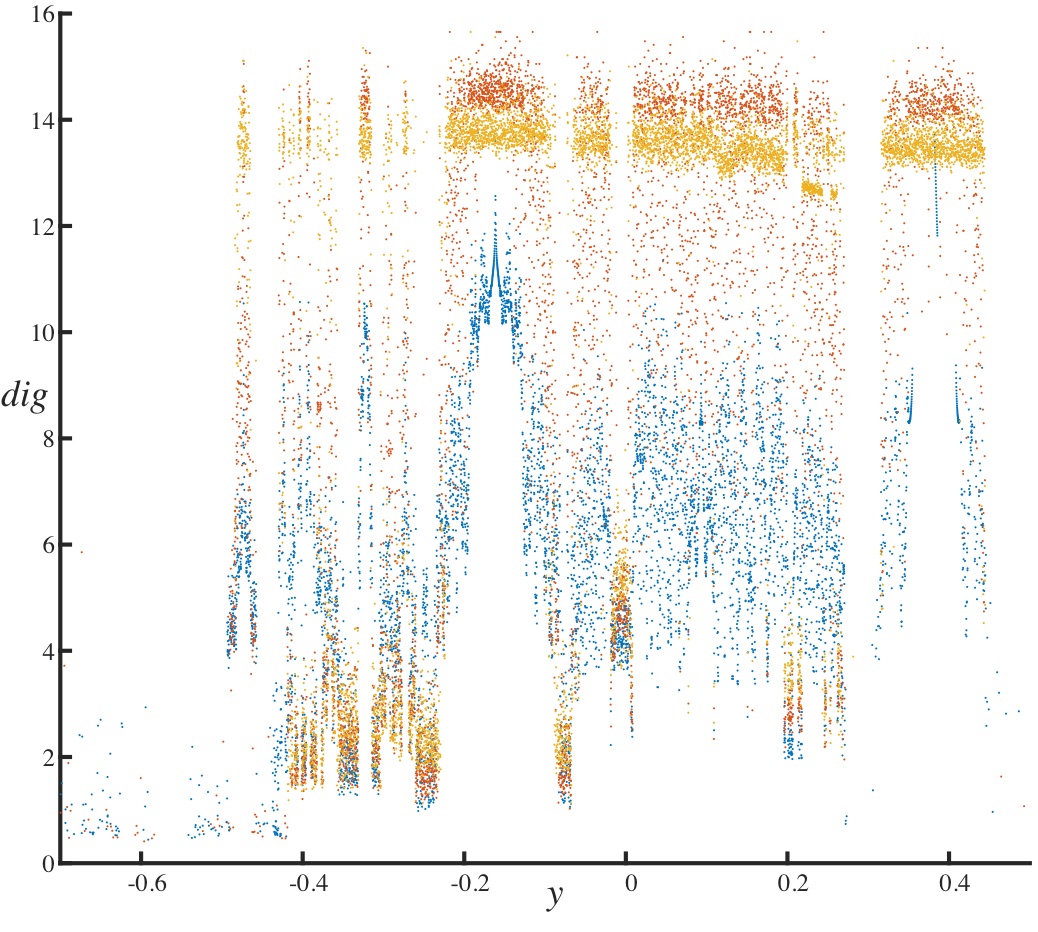}{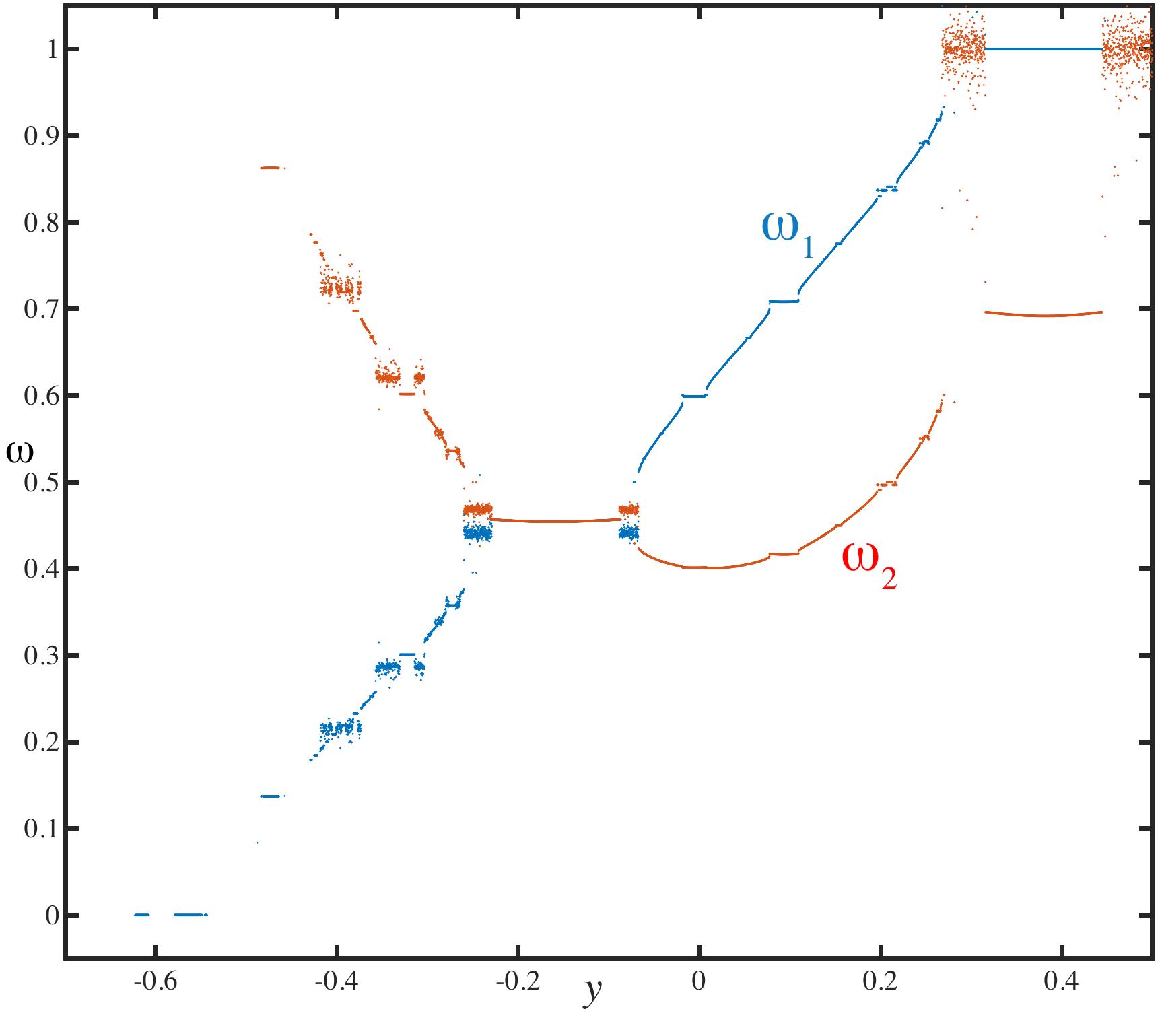}{
Orbits of \Eq{AngleActionMap} for $\delta=-0.4$, and $\eps = 0.02$ with initial conditions $(0,0,y)$
and a grid of $10^4$ initial $y$ values between $-0.7$ and $0.5$ 
(a) Slice for $|x_2| \le 0.005$  showing the $(y,x_1)$ phase space. 
(b) $\digT$, \Eq{digi}, as a function of initial $y$ for 
$T= 10^3$ (blue), $3(10)^4$ (orange), $10^6$ (red). The value
of $\digT$ changes significantly more at the edge of resonant tubes 
than it does in the middle of the tube, and for small $T$ the accuracy of the computation is low.
(c) The two components of the rotation vector $\omega$ using $T=10^6$.  
}{digitsviterates}{2.3in}

The dichotomy between the values of $\digT$ for chaotic and regular
orbits is also reflected in histograms of $\digT$, shown in
\Fig{histc8}. These show the fraction of those orbits with $\omega \in
[0,1]^2$ that have a given value of $\digT$ for a range of $\delta$,
$\eps$, and initial conditions $(0,0,y_0)$,  with $y_0$ and $\delta$
chosen such that $(y_0,\delta) \in \cP$, \Eq{ydeltaGrid}. 
For the smaller $T$, panel (a), there
are clear peaks near $\digT = 2$ and $14$ corresponding to chaotic and
regular orbits, respectively, but there is also a broad shoulder with $8
< \digT < 13$ that corresponds to orbits for which the distinction is
less clear. Note, however, that in panel (b), where $T = 10^6$, this
middle peak has moved to larger values of $\digT$, leaving only a
smaller tail just below the peak at $\digT = 14$. 

\InsertFigTwo{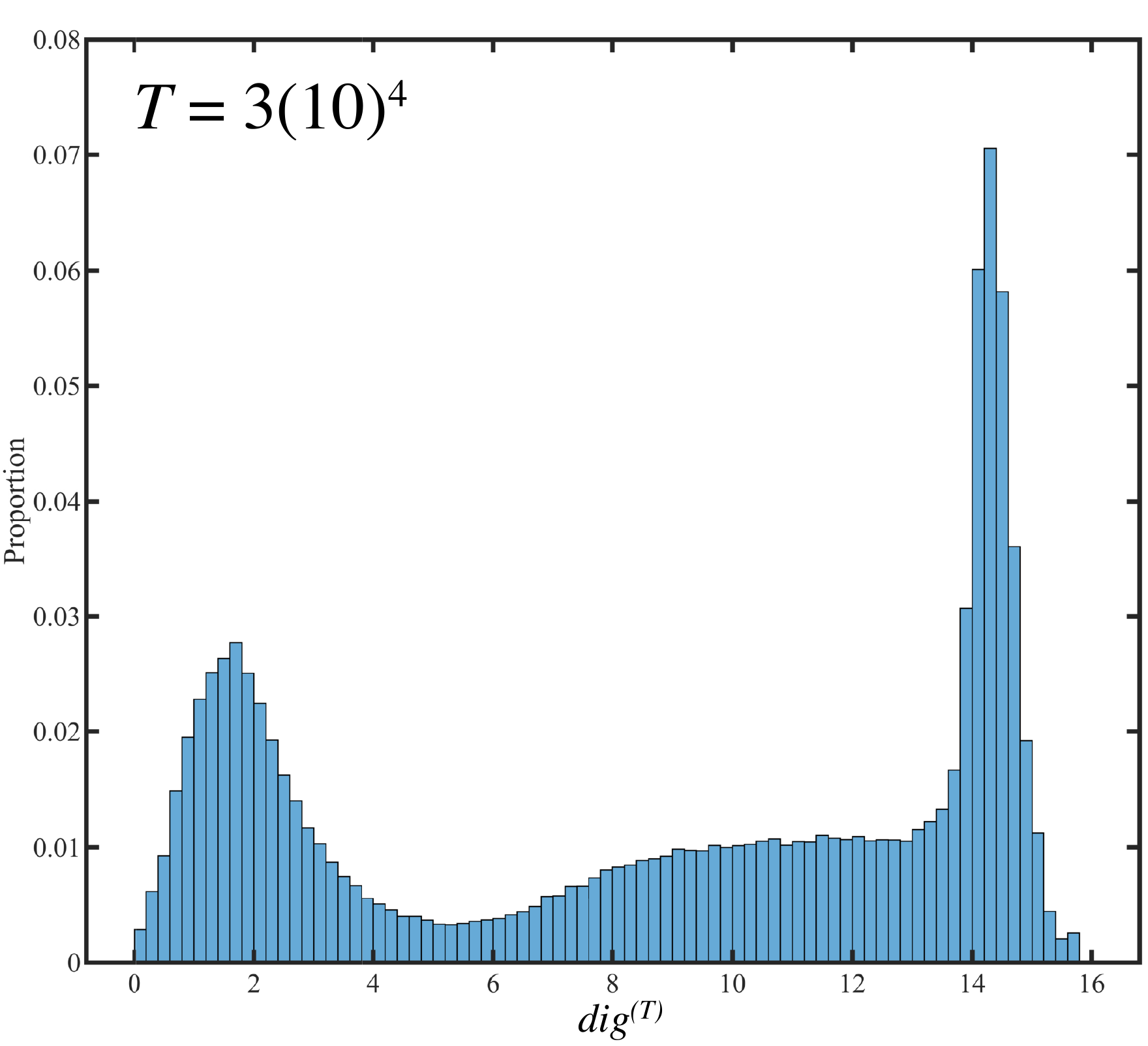}{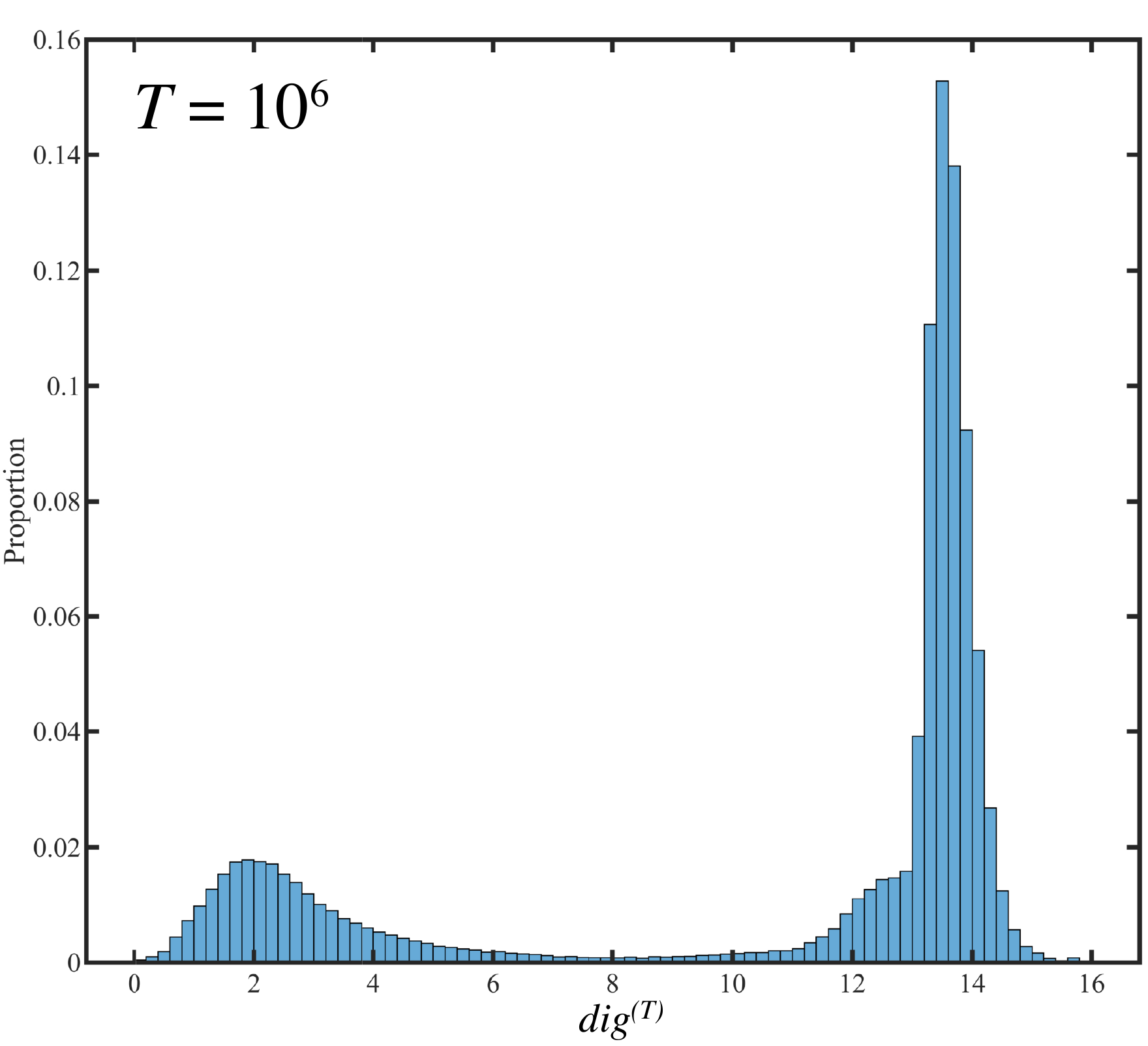}{
Frequency histograms of $\digT$ for (a) $T=3(10)^4$ and (b) $T = 10^6$ 
for a $100 \times 100$ grid on the domain $\cP$, \Eq{ydeltaGrid}, for initial conditions ($0,0,y_0)$
and 50 values of $\eps \in [0.005,0.055]$. 
}{histc8}{3in}

As we saw in \Fig{zdelta}, as $\eps$ increases the chaotic region
expands. Many of the chaotic orbits leave the interval $-0.7< y< 0.5$;
as a consequence, the computed frequencies for these orbits will be
outside $[0,1]^2$. We think of the orbits that leave
the $\omega$ range as essentially \textit{unbounded}, though we cannot
guarantee that there are no rotational tori acting as barriers at larger
or smaller action values.
 \Fig{inRange}(a) shows the proportion of the $10^4$ orbits in \Fig{histc8} that 
are \textit{bounded}  as $\eps$ grows. Since $\cP$ includes a buffer, at $\eps =  0$, the 
proportion is $1/1.1 \approx 91\%$, but  
 by the time $\eps = 0.055$, that proportion has dropped to about 10\%.  
Only the bounded orbits were used in
the histograms in \Fig{histc8}. 

Figure \ref{fig:inRange}(b) shows for bounded orbits   
how the proportions of values of $\digT$ that are small, intermediate, and large
vary with respect to $T$.  As $T$ grows, the fraction of orbits in the
intermediate range, $\digT \in (5,11)$, decreases, and the corresponding
fractions in the lower and upper ranges saturate.
In order to choose a measure to provide a good separation between
order and chaos, we set $T = 10^6$ for subsequent computations.

\InsertFigTwo{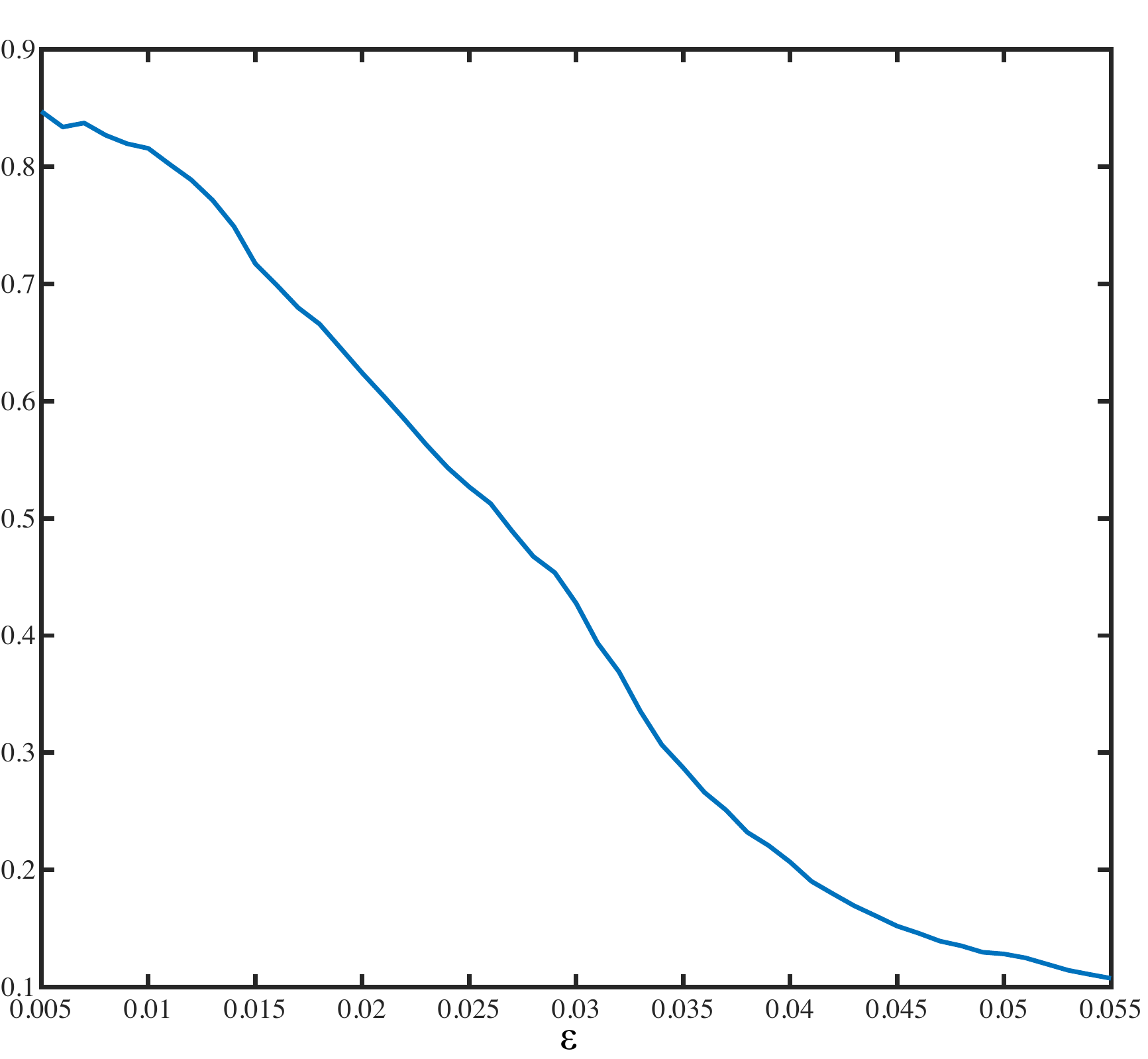}{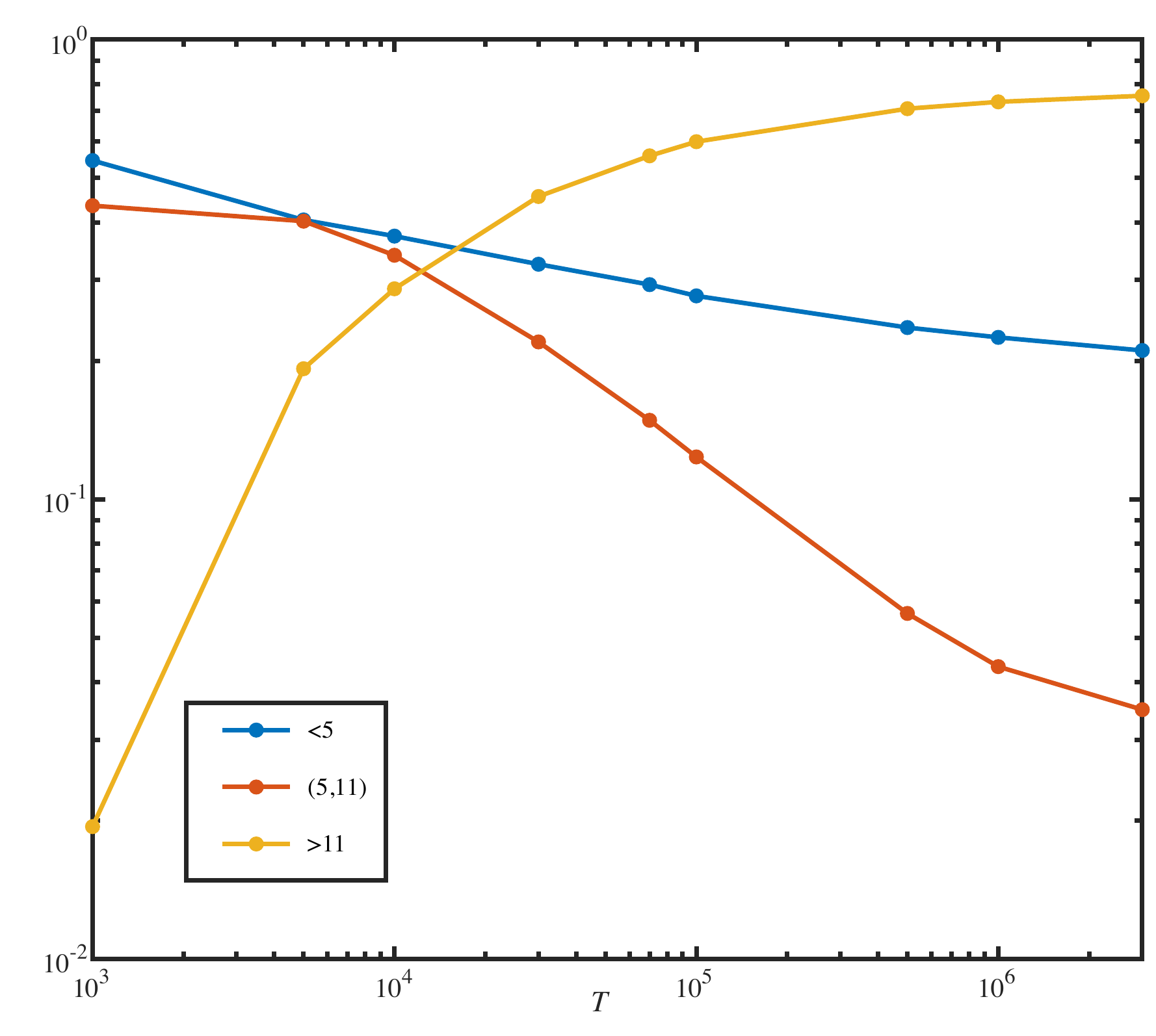}{
(a)  The proportion of orbits in $\cP$ with $\omega \in [0,1]^2$  
as a function of $\eps$. 
As the tori are destroyed, only a small fraction of 
orbits have computed rotation vectors in this range. 
(b) Log-log plot of the proportion of bounded orbits with $\digT$ $<5$ (blue), between $5$ and $11$ (red), 
and $>11$ (yellow) as a function of $T$. 
}{inRange}{3in}

Fixing $T = 10^6$, we next need to choose 
an appropriate cutoff value for $dig^{(10^6)}$ to distinguish order from chaos. 
\Fig{propchaos} shows the proportion of orbits with $dig^{(10^6)} < 7,9,11$. 
Even  though varying the cutoff does give quantitative differences, the proportions
have the same qualitative form as $\eps$ varies.
In particular, the proportions grow when $\eps < 0.03$, reflecting the increasing fraction
of bounded orbits that are chaotic. Beyond the peak at $\eps = 0.03$, the fraction of
unbounded orbits increases rapidly as the tori, which act as transport barriers,
are destroyed, allowing the escape of previously trapped, chaotic orbits.
Since we want to be conservative in classifying an
orbit as a regular---as well as to guarantee 
that $\omega$ has high accuracy---we  use the cutoff 
\beq{chaosCutoff}
	dig^{(10^6)} > 11
\eeq
to declare that an orbit is ``nonchaotic.''

\InsertFig{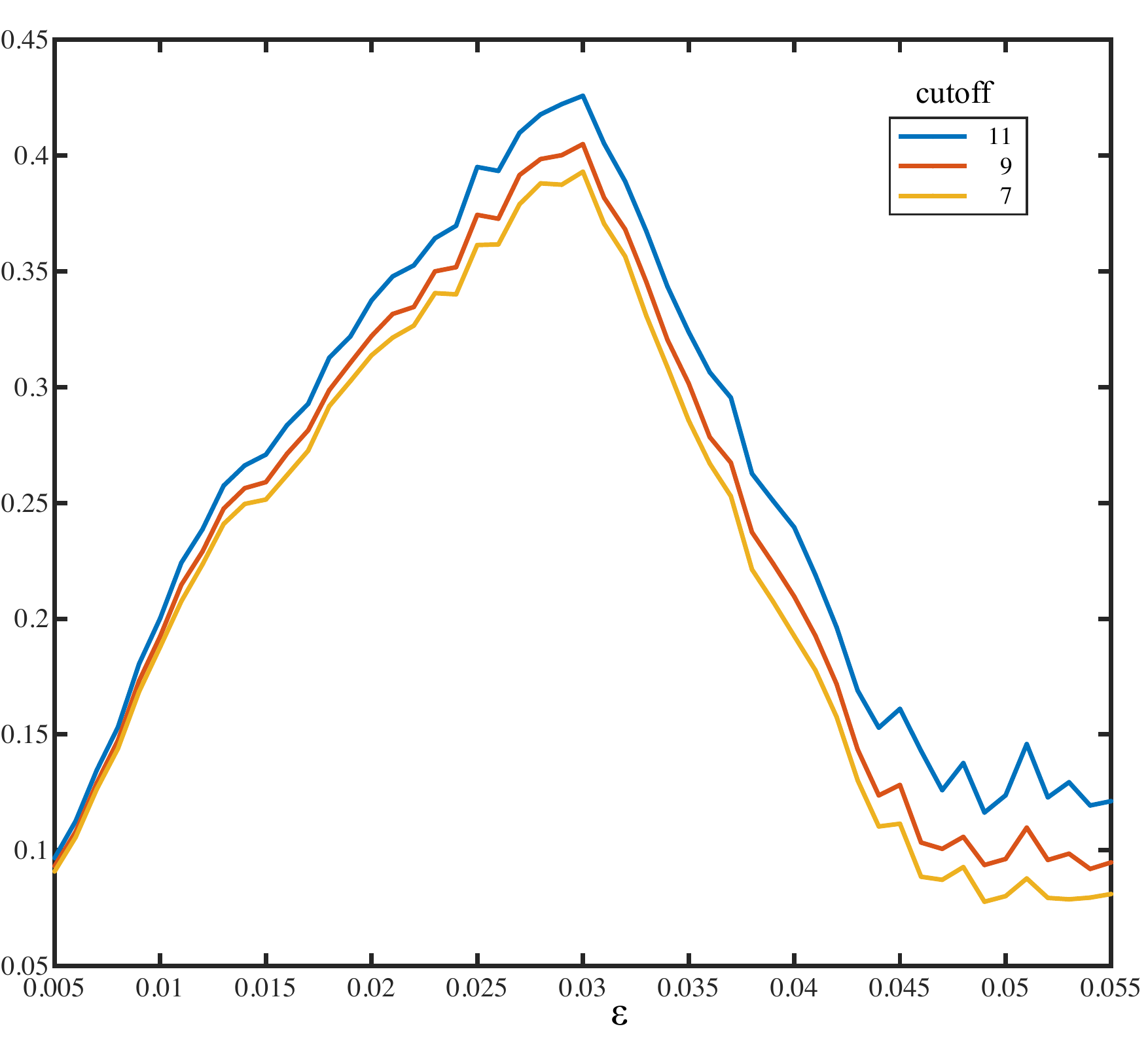}
{The proportion of  orbits that would labeled as ``chaotic'' in $\cP$ with $\omega \in [0,1]^2$ 
using the criteria  $\digT< 7,9,11$ respectively with $T = 10^6$. 
In each case the proportion peaks near $\eps = 0.03$. By criterion \Eq{chaosCutoff},  the 
blue (11) curve shows the proportion of chaotic orbits.}
{propchaos}{3in}

\Fig{crit3d} shows the set of frequencies for the nonchaotic, bounded orbits as a
function of $\eps$ as identified using the criterion \Eq{chaosCutoff}. 
Note that this number drops significantly for large values of $\eps$. 

\InsertFig{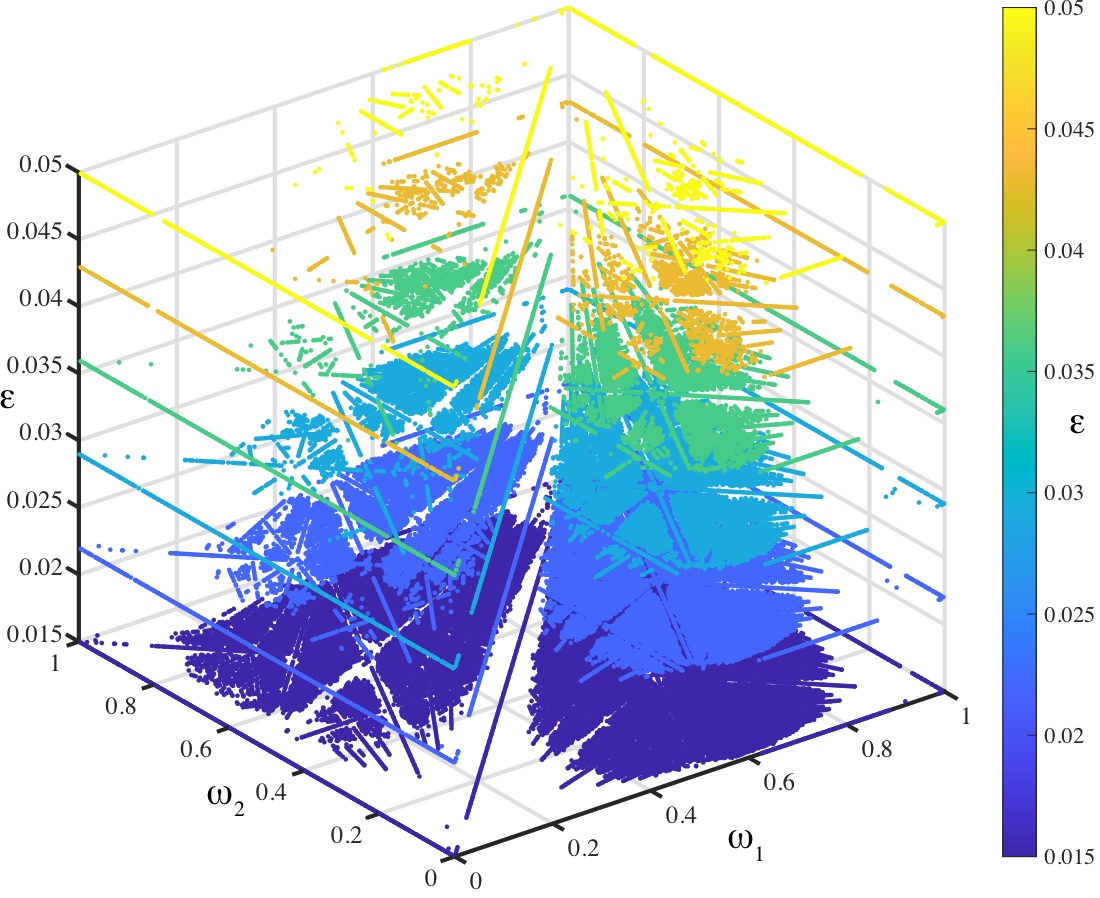}
{Rotation vectors $\omega \in [0,1]^2$ for the $1,613,136$
nonchaotic, bounded orbits, computed using a $1000^2$ grid in the domain $\cP$, \Eq{ydeltaGrid}, and 
 $\eps \in [0.015,0.022,0.029,0.036,0.43,0.05]$, determined using the criteria in \Eq{chaosCutoff}. As
 we see in \Sec{Resonances}, the straight lines correspond to orbits trapped in low-order resonant tubes.}
{crit3d}{3in}

\subsection{Distinguishing Resonances}\label{sec:Resonances}

In this section, we seek a numerical method to distinguish between resonant and incommensurate  vectors.
For a given  $\omega \in \bR^d$, define the resonant module
\beq{Module}
	\cL(\omega) \equiv \{ m \in \bZ^d: m\cdot\omega \in \bZ \} .
\eeq
We say that $\omega$ is \textit{ incommensurate} when $\cL(\omega) =
\{0\}$. By contrast, $\omega$, is \textit{resonant} if \Eq{Module} is
nontrivial, i.e., if  there is a nonzero vector $m \in \bZ^d$ that
satisfies \Eq{Resonance}. Of course, if $m, m' \in \cL(\omega)$, then so
are $m + m'$ and $km$ for any $k \in \bZ$: the set \Eq{Module} of such
vectors is a module. The length, $M = \|m\|_1$, of the smallest
(nonzero) integer vector $m$ in $\cL(\omega)$ is the \textit{order} of
the resonance.

The \textit{rank} of a resonant frequency $\omega$ is the dimension of $\cL(\omega)$.
We say that a frequency is \textit{rational} if $\dim(\cL(\omega)) = d$. 
In this case there is a $(p,q) \in \bZ^{d+1}$ so that $\omega =
\frac{p}{q}$; i.e., $\omega \in \bQ^d$. When $d > 1$, every rational
frequency is resonant, but the converse need not be not true. For
example, the vector $(\sqrt{2}, 2+3\sqrt{2})$ is resonant with $(m,n) = (-3,1,2)$,
but it is not rational. For this example the resonance order is $M = 4$.

Since we can compute the frequency vector for an orbit only to finite precision, we can only evaluate
resonance up to some precision. If $\omega \in \bR^d$ is $(m,n)$ resonant, then it lies in the codimension-one plane
\beq{ResonantPlanes}
	\cR_{m,n} = \{\alpha \in \bR^d :  m \cdot \alpha -n = 0\}.
\eeq
The collection of resonant vectors is
\beq{ResonantDomain}
	\cR = \bigcup_{m,n \in\bZ^{d+1} \setminus \{0\}} \cR_{m,n}.
\eeq
For the case $d=2$ of interest here, the lines up to order $M = 8$ are shown for a portion of the $\omega$-plane in \Fig{SpiralResLines}. 
Of course $\cR$ is dense in $\bR^d$, as are the Diophantine vectors \Eq{Diophantine}.

\InsertFig{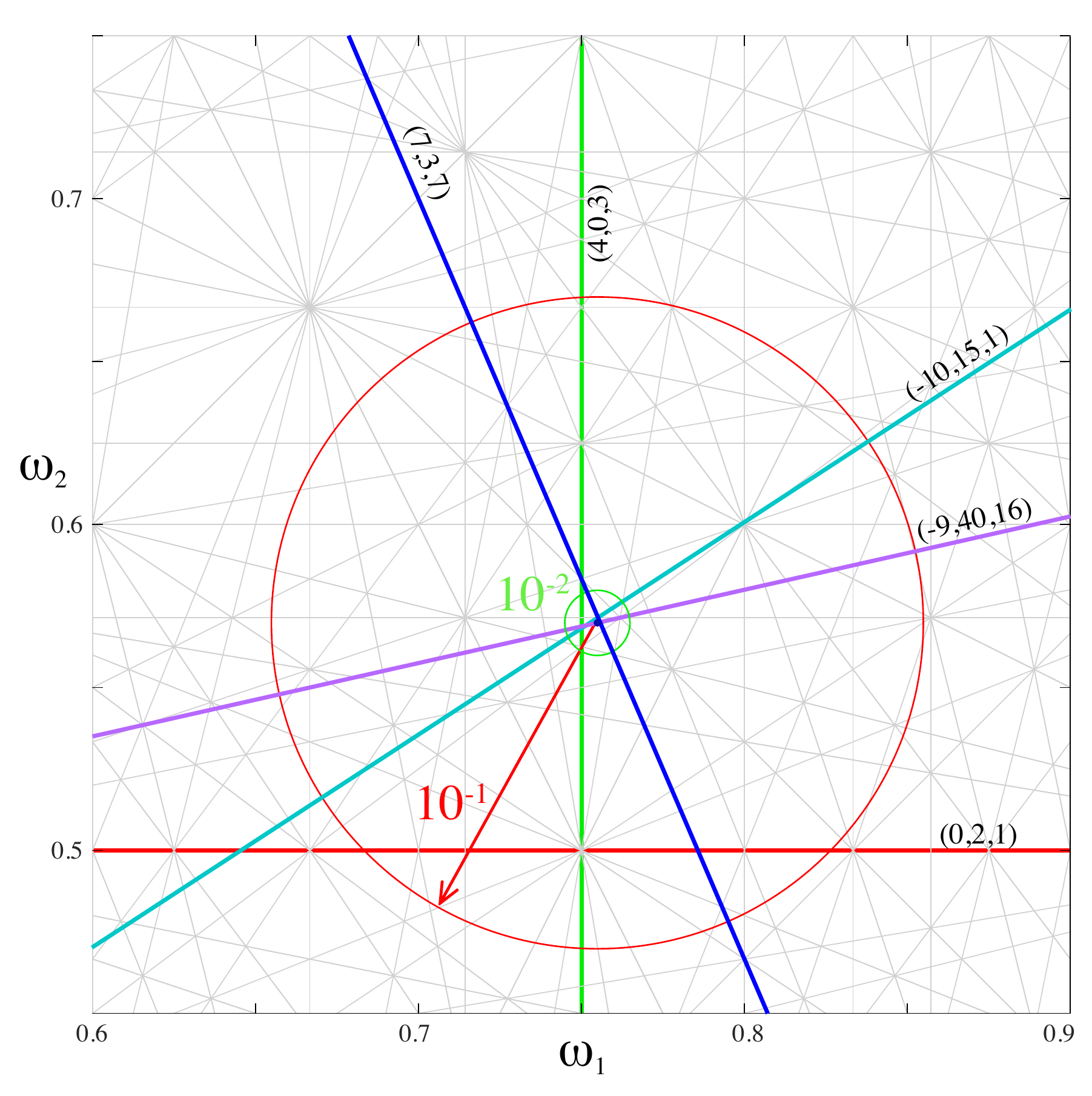}{Resonant lines (grey), $\cR_{m,n}$, up to order $M = 8$. Also shown are those of minimal order for the spiral mean frequency $(\sigma^{-1},\sigma^{-2})$, see \Eq{SpiralField}, for five values of
$\rho$, $10^{-1}$ to $10^{-5}$ from \Tbl{SpiralResonances}.}{SpiralResLines}{4in}

We will say a vector $\omega$ is $(m,n)$ resonant \textit{to precision $\rho$} if
the resonant plane intersects a ball of radius $\rho$ about $\omega$, i.e., if 
\beq{ResonancePlane}
	\cR_{m,n} \cap B_\rho(\omega) \neq \emptyset .
\eeq
Using the Euclidean norm, the minimum distance between the resonant plane and the point $\omega$ is
\beq{minDist}
	\Delta_{m,n}(\omega) = \min_{\alpha \in \cR_{m.n}} \| \alpha - \omega\|_2 
	                     = \frac{|m \cdot \omega-n|}{\|m\|_2} .
\eeq
Thus we say that $\omega$ is $(m,n)$ resonant to precision $\rho$, whenever $\Delta_{m,n}(\omega) \le \rho$.

Given a vector $\omega$ and a precision $\rho$, what is the smallest order resonance with $\Delta_{m,n}(\omega) \le \rho$? For the one-dimensional case ($d=1$), the answer to this question can be efficiently computed using the Stern-Brocot (or Farey) tree. Indeed, as we previously noted \cite{Sander20}, for any $\rho>0$, the rational $\tfrac{p}{q}$ with the smallest denominator in the interval $[\omega-\rho, \omega +\rho]$ is the first such rational on the tree that falls in that interval. The Stern-Brocot tree is essentially the generalization of the continued fraction to include ``intermediates" as well as convergents of $\omega$.

As far as we know, there is no generalization of this result for $d > 1$. Since there are finitely many $m \in \bZ$ with $\|m\|_1 \le M$, a brute force computation is of course possible for modest values of $M$. For example, given any $\rho >0$, and ignoring issues of floating point arithmetic, \Alg{ResonanceOrder} will return
\beq{Mres}
	M(\omega,\rho) = \min \{\|m\|_1 : \Delta_{m,n}(\omega) \le \rho\} ,
\eeq
which we could call $\rho$-order of $\omega$.

\begin{algorithm}  
\caption{Minimal resonance order \Eq{Mres} to precision $\rho$ for $\omega \in \bR^d$.}\label{alg:ResonanceOrder}
	\begin{algorithmic}
	\Procedure{ResonanceOrder}{$\omega$,$\rho$}
	\State $M = 0$, $\Delta = 1$	
	\While{$\Delta > \rho$}		
		\State $M \leftarrow M+1$
	\For{$m_2 = -M$ to $M$}		
		\State $m_1 = M- |m_2|$
		\State $n = \mbox{round}(m \cdot \omega)$
		\State $\Delta = \min \left(\Delta, \frac{|m \cdot \omega -n|}{\|m\|_2}\right)$
	\EndFor

	\EndWhile\\
	\Return $M$
	\EndProcedure
	\end{algorithmic}
\end{algorithm}

As an example, consider the so-called \textit{spiral} \cite{Kim86} or \textit{plastic} \cite{Stewart96} mean, 
the real solution to
\beq{SpiralField}
	\sigma^3 = \sigma +1 \Rightarrow \sigma \approx 1.324717957244746 .
\eeq
This generates an algebraic field $\bQ[\tau]$ with integral basis $(1,\tau,\tau^2)$.
The sequence of minimal order resonances for the frequency 
$(\sigma^{-1},\sigma^{-2}) = (\sigma^2-1,\sigma-\sigma^2+1)$
with tolerances  $\rho = 10^{-j}$ for $j$ up to $14$ is shown in \Tbl{SpiralResonances}. For example, $M(\omega, 10^{-9}) = 1119$. Note that this vector is Diophantine \Eq{Diophantine} \cite{Cusick72}. The first five optimal resonant lines are shown in \Fig{SpiralResLines}.

\begin{table}[htp]
\begin{center}
\begin{tabular}{r|r| r r r}
$\log_{10}(\rho)$ & $\|m\|_1$ & $m_1$	&$m_2$	&$n$\\
\hline
-1	& 2		&0		&2		&1\\ 
-2	& 4		&4		&0		&3\\
-3	& 10	&7		&3		&7\\ 
-4	& 25	&-10	&15		&1\\ 
-5	& 49	&-9		&40		&16\\ 
-6	& 96	&7		&89		&56\\ 
-7	& 208	&171	&-37	&108\\ 
-8	& 387	&316	&71		&279\\
-9	& 1119	&-350	&769	&174\\ 
-10	& 2064	&-176	&1888	&943\\ 
-11	& 4306	&3952	&354	&3185\\ 
-12	& 10322	&6783	&3539	&7137\\ 
-13	& 24301	&10676	&-13625	&295\\ 
-14	& 48897	&-10971	&37926	&13330\\
\end{tabular}
\end{center}
\label{tbl:SpiralResonances}
\caption{Optimal resonances and resonance orders for the frequency $\omega = (\sigma^{-1},\sigma^{-2})$,
seen in \Fig{SpiralResLines}, as the precision $\rho$ decreases. 
Note that without loss of generality, we can assume that $n$ is nonnegative.
}
\end{table}

The resonance orders \Eq{Mres} in \Tbl{SpiralResonances} grow as a power of the inverse of the precision $\rho$ 
with the best-fit 
\[
	M(\omega,\rho) \simeq  0.944\,\rho^{-0.336}. 
\]
Here is the intuition as to why this occurs:
by \Th{Minkowski} in \App{Diophantine}, for each $K>0$,  there is an $m \in \bZ^2$ with 
$\|m\|_\infty \le K$ such that for $p=2$
\beq{MinkowskiBound}
	\Delta_{m,n} \le \frac{1}{\|m\|_2 K^p}.
\eeq
Furthermore, a result of Laurent, see e.g., \cite[p.~693]{Waldschmidt12}, implies that 
for almost all $\omega \in \bR^2$, it is 
not possible to satisfy this equation for any value of $p>2$. Therefore, we expect that the typical value of 
$\|m\|_\infty$ is close to the maximal value, i.e., that the 
satisfaction of this bound requires that $\|m\|_\infty \sim K$.
Then since the norms of $m$ are equivalent, choosing $K = \rho^{-1/3}$, we get 
$\Delta_{m,n} \lesssim \rho$, for $\|m\|_1 \sim \rho^{-1/3}$.

More generally, for a given $\rho$, we computed the minimal resonance order \Eq{Mres} for a set of equi-distributed, random $\omega \in [0,1]^2$, see \Fig{ResonanceRandom}(a). The log of these values have mean $\langle \log_{10} M(\omega,10^{-9})\rangle = 2.92$ and standard deviation $0.171$, though the distribution differs significantly from a normal with the same mean and deviation (the red curve in the figure). For the $10^4$ randomly chosen $\omega$ we found that
\beq{Mmax}
	M(\omega,10^{-9}) \le m_{max} = 3841,
\eeq
and only six cases had $M > 2500$. A similar distribution holds for other values of $\rho$. As shown in \Fig{ResonanceRandom}(b), the mean of $\log(M)$ depends linearly on $\log(\rho)$, with the best fit
\beq{MeanRes}
	\langle \log_{10} M(\omega,\rho) \rangle = -0.334\log_{10}(\rho) -0.091,
\eeq
which is again consistent with \Eq{MinkowskiBound}.
\InsertFigTwo{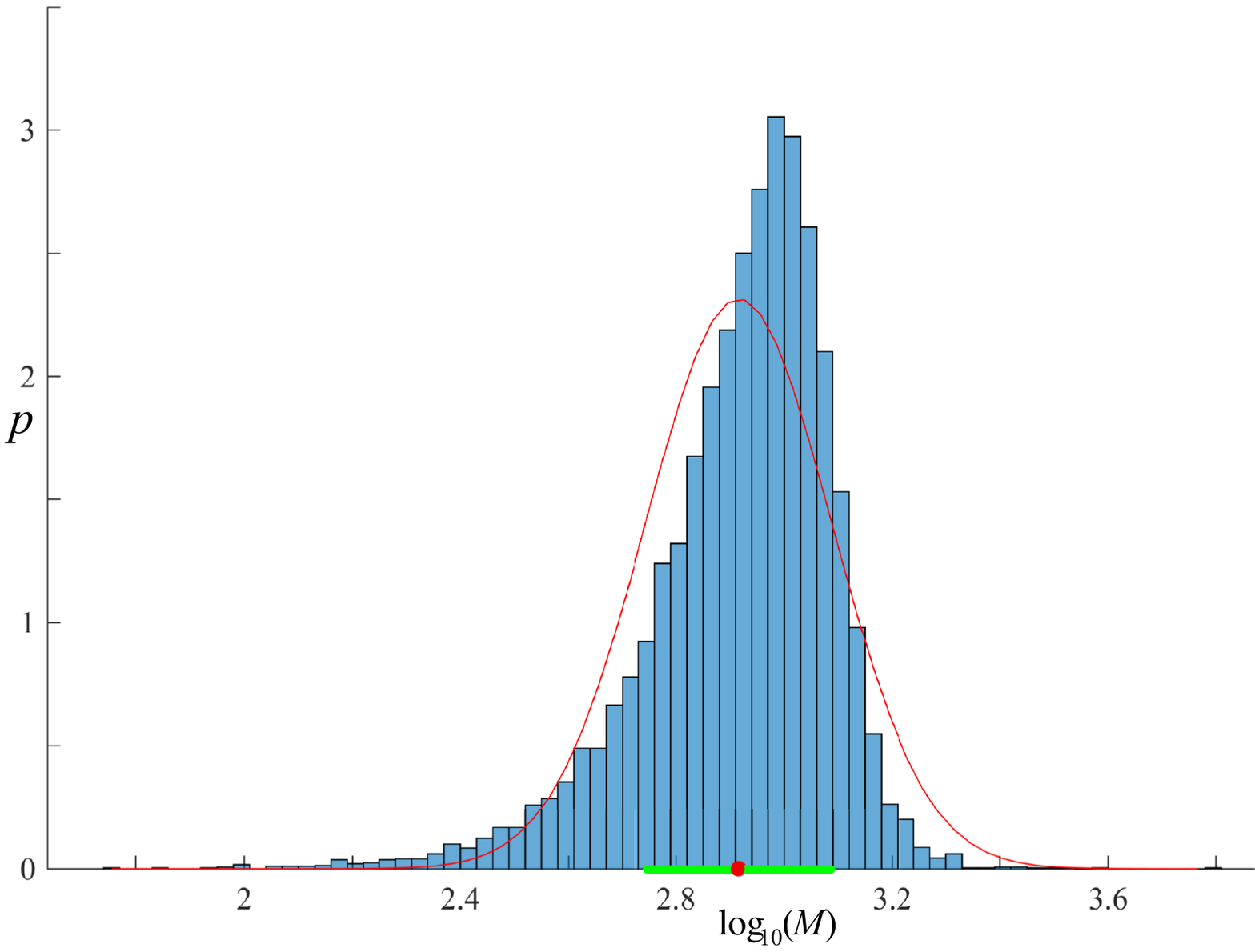}{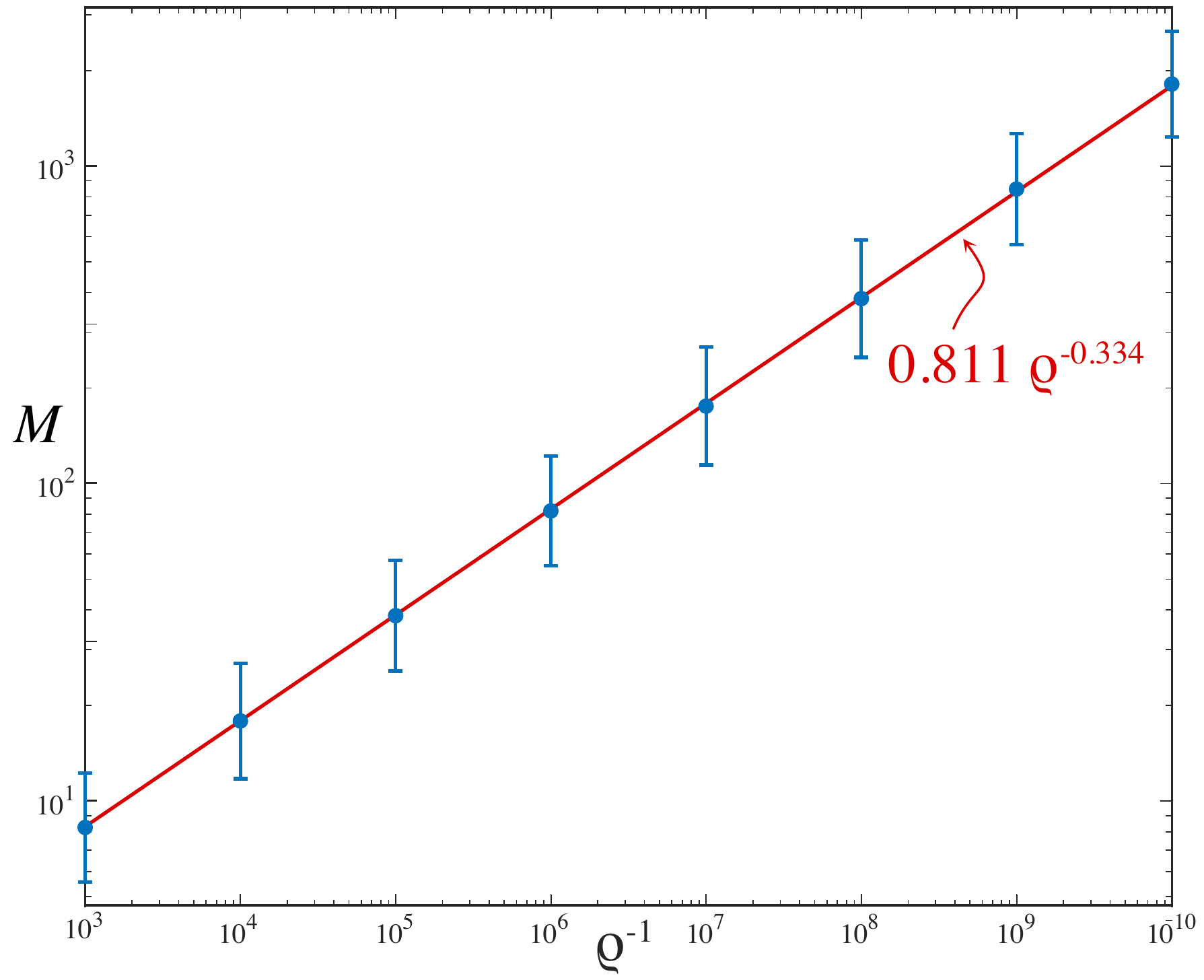}{(a) Probability density of the logarithms of 
minimal resonance orders \Eq{Mres} for $10^4$ random vectors with precision $\rho = 10^{-9}$. 
The dot (red) on the horizontal axis shows the mean, $2.92$, 
the thick line (green) shows one standard deviation $0.171$,
and the curve (red) shows the normal distribution with these parameters.
(b) A log-log plot of the mean and standard deviation of resonance order as 
a function of precision for a sample of $2000$ random vectors. 
The line (red) is the least squares fit \Eq{MeanRes}. }{ResonanceRandom}{3in}

The computation of the minimal resonance order is applied to the dynamical frequency vectors in \Fig{nonChaoticFreqs} using $\rho = 10^{-9}$. The  data corresponds to the nonchaotic orbits on a grid of $(y,\delta)$ for $\eps = 0.043$.
The orbits with $M < 8$ (dark blue), clearly lie on the low-order resonant lines (grey) shown in the figure. Of the nonchaotic orbits at this value of $\eps$, only $76$ have $M > 250$, and only $79$ have $M > 200$. It is clear that for this value of $\eps$, there are very few rotational tori. 
A similar picture is obtained when more values of $\eps$ are included in \Fig{rot1K}. The left panel shows the frequency vectors for the nonchaotic orbits of \Fig{crit3d}, now projected onto the $\omega$ plane. Note that many of the resonant lines lie in gaps in the figure, with clear clusters of points along the resonances. Indeed, if we change the color scale to be $\log_{10} M(\omega,\rho)$ for $\rho = 10^{-9}$, the resonances again show up as dark blue lines, see \Fig{rot1K}(b).

\InsertFig{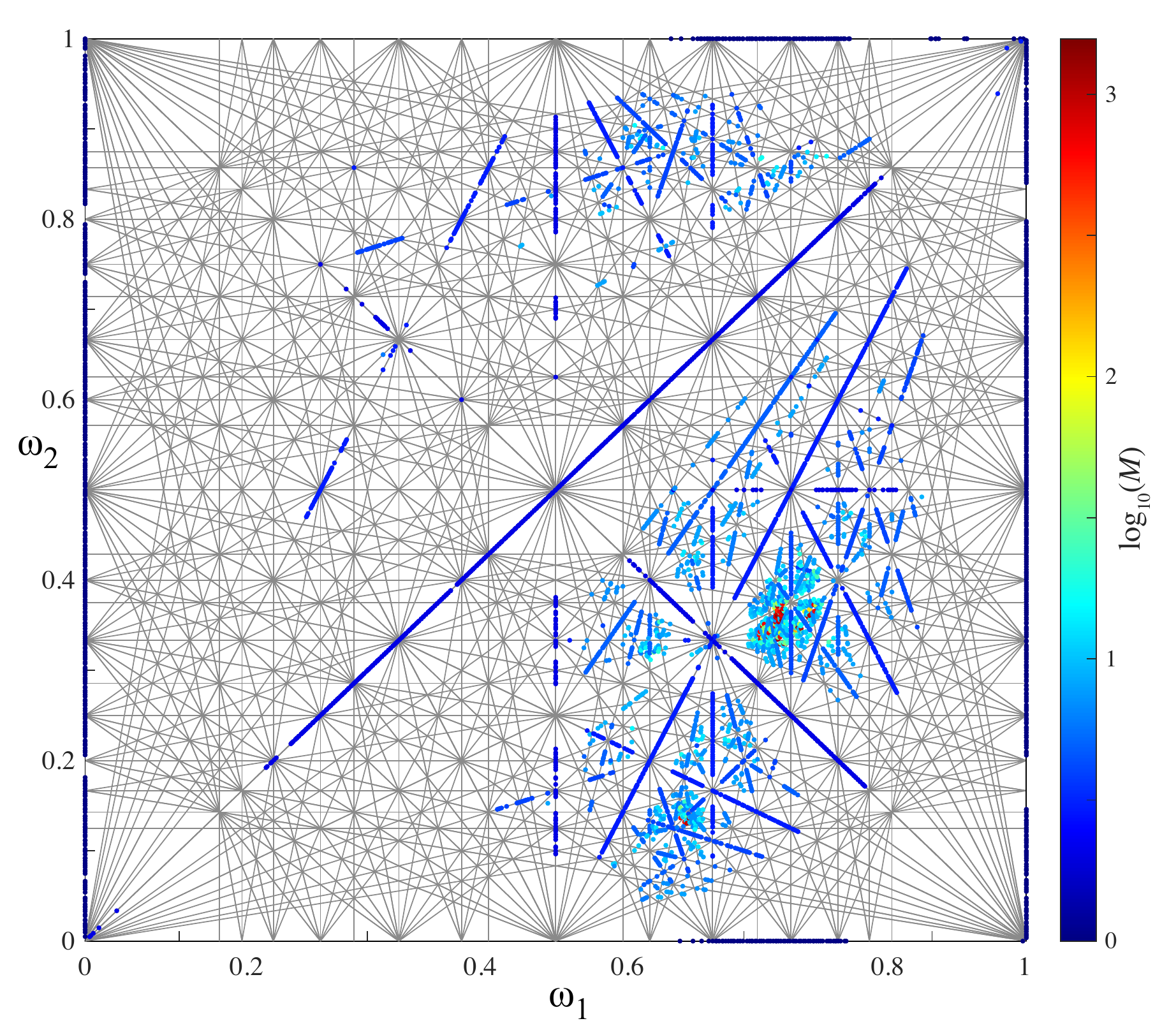}{Frequencies of the $140,338$ nonchaotic orbits on a $1000^2$ grid in \Eq{ydeltaGrid} of initial conditions for $\eps = 0.043$.
The color scale is $\log_{10}(M)$. Only $76$ of these orbits have $\log_{10}(M)>2.4$.
Also shown are the resonant lines, $\cR_{m,n}$, up to order $8$.}{nonChaoticFreqs}{4in}

\InsertFigTwo{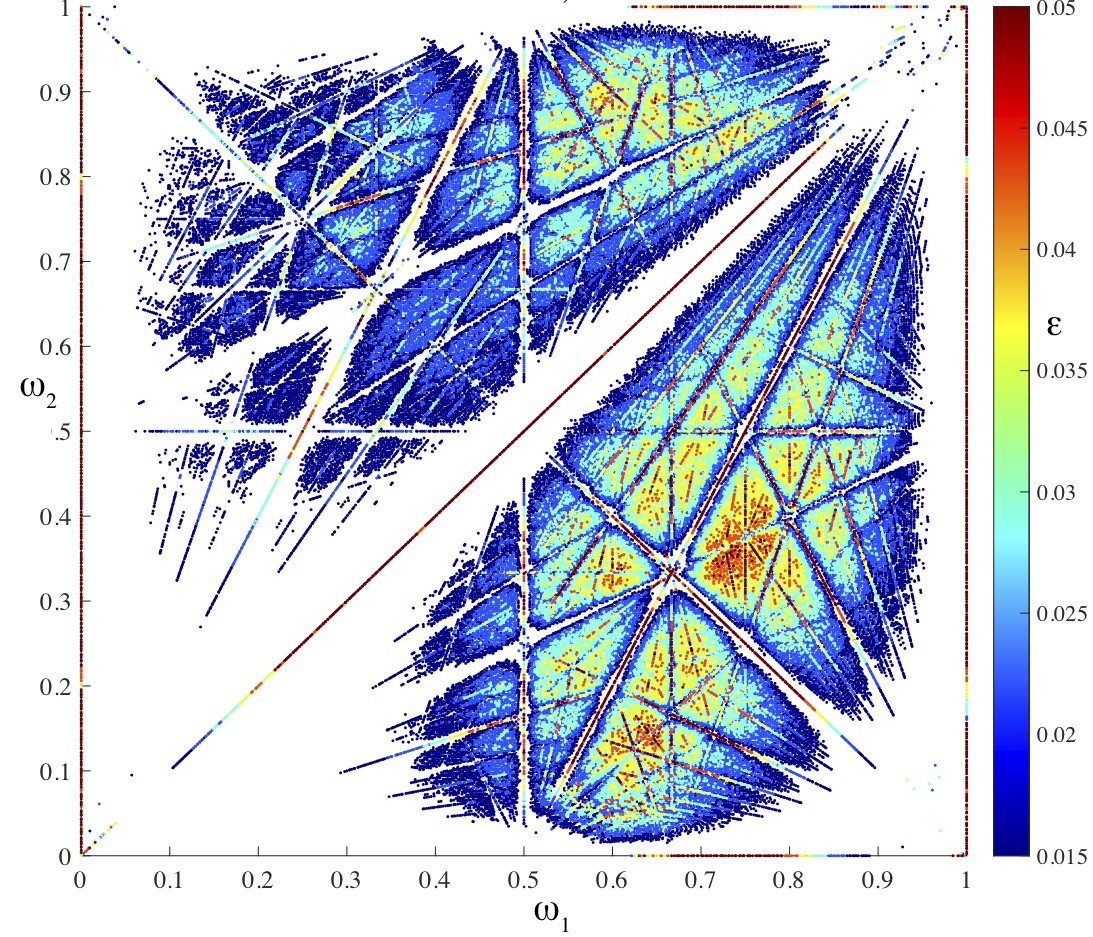}{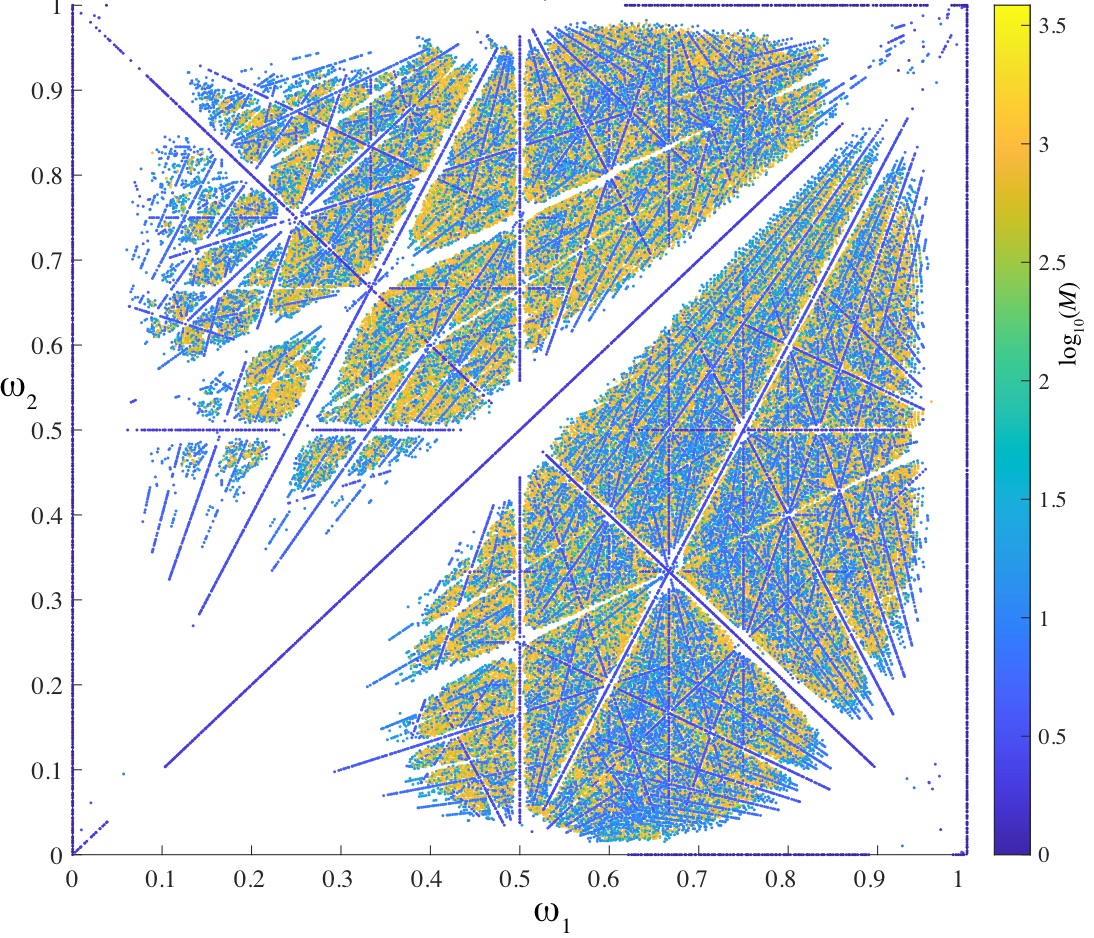}{
(a) The rotation vectors for nonchaotic orbits, using the 
same data of \Fig{crit3d}, now 
viewed in a two-dimensional projection. The color scale gives largest $\eps$ for which
there is a nonchaotic orbit for the given $(y_0,\delta)$.  
(b) The same data, but this time colored using the $\rho$-order, \Eq{Mres}.}
{rot1K}{3in}

Since randomly chosen $\omega$ will almost always be incommensurate, the dynamically obtained frequency vectors that are resonant  should have values of \Eq{Mres} below the bulk of the distribution shown in \Fig{ResonanceRandom}.  This is confirmed in \Fig{histlres}(a), a histogram of resonance orders for the orbits from \Fig{crit3d}. Note how resonant tubes in the dynamics change the histogram from that of the random frequencies in \Fig{ResonanceRandom}. Indeed about $60\%$ of these orbits have $M= 1$ or $2$ corresponding to the largest resonances due to the Fourier terms of the force \Eq{Parabola}, and only $20\%$ have $M > 250$, i.e., the bulk of the domain of \Fig{ResonanceRandom}.

\InsertFigTwo{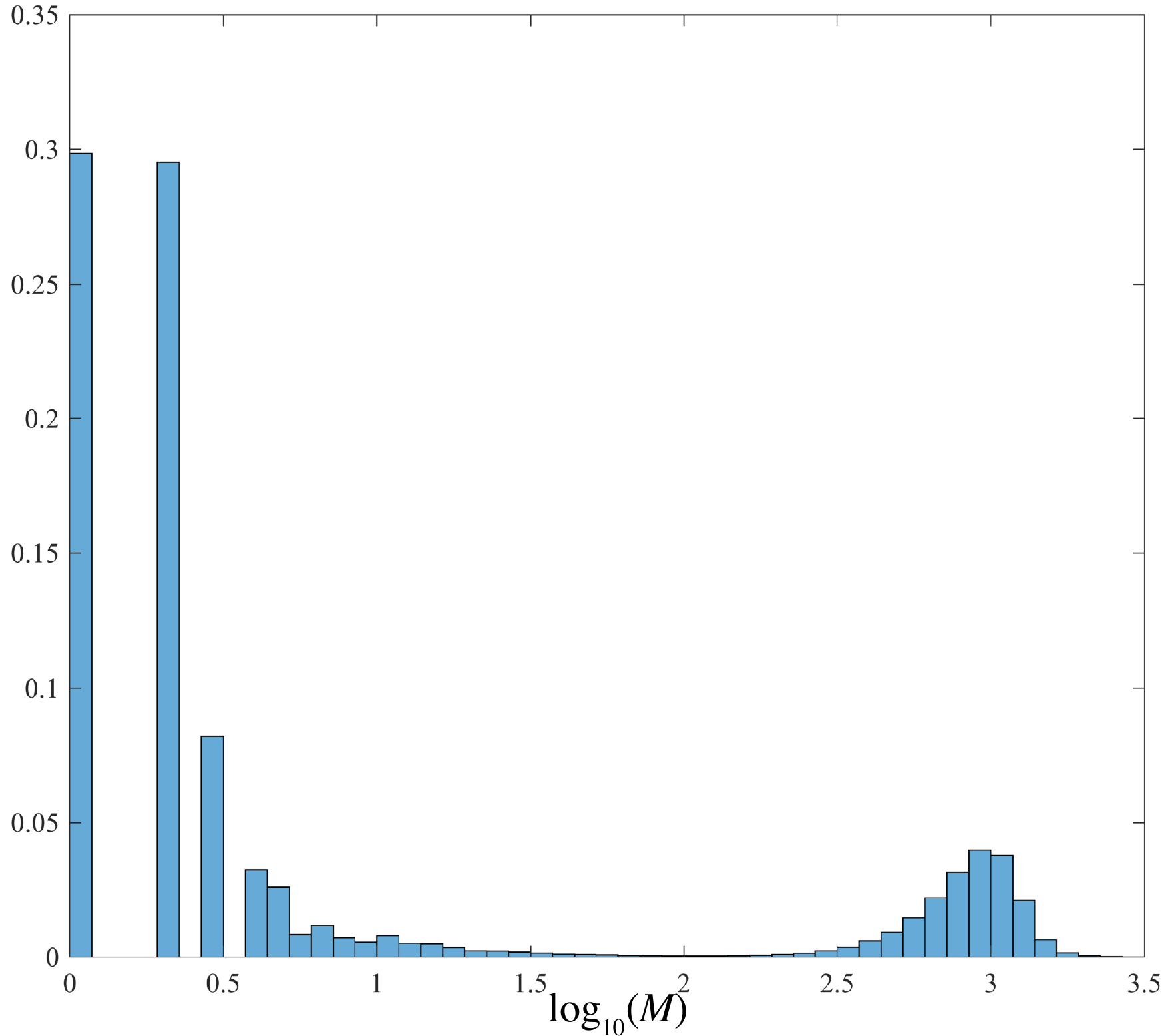}{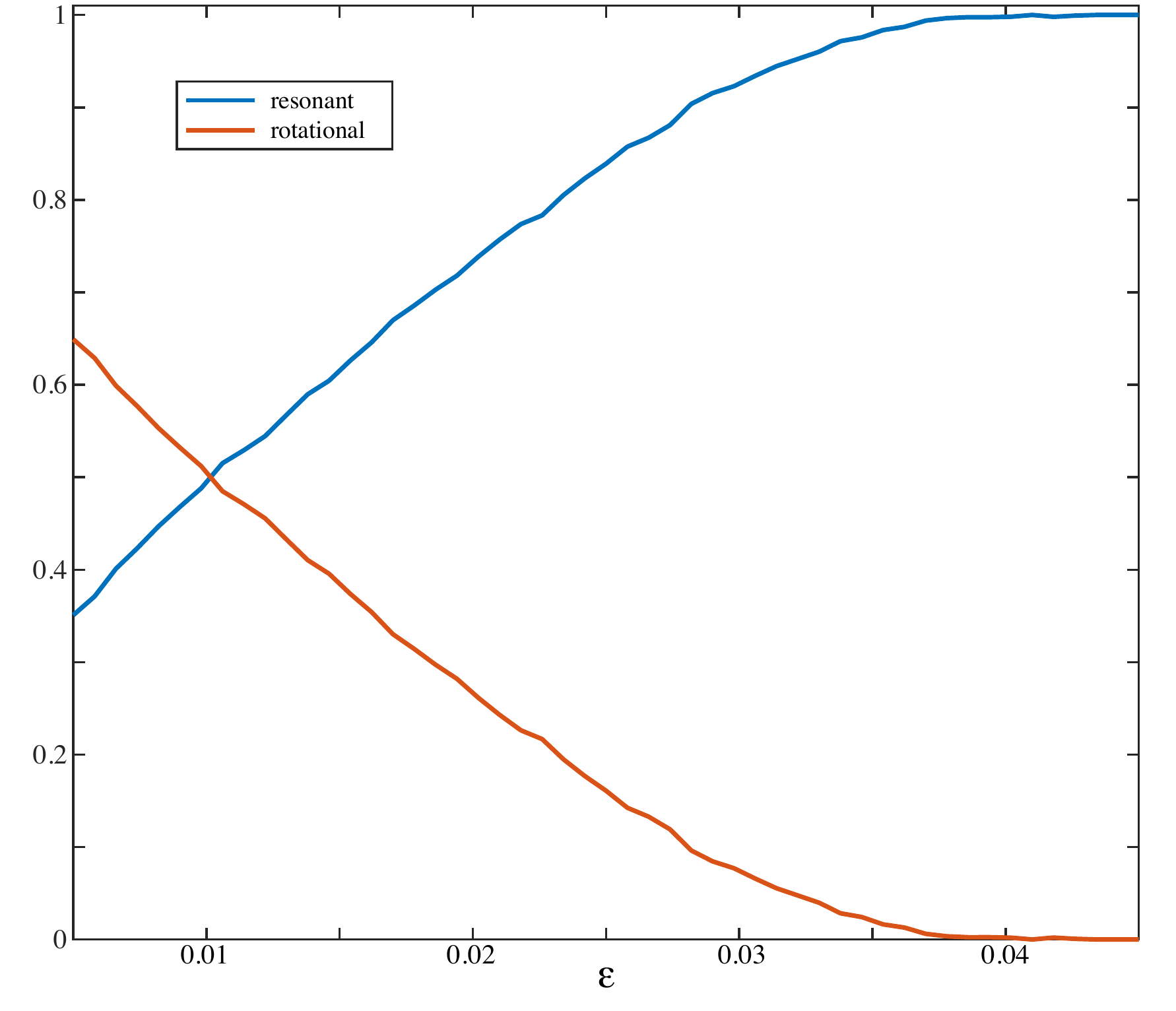}{ (a) Histogram of $\log_{10}(M)$ for the data of \Fig{crit3d}. 
(b) The fractions of nonchaotic orbits with $\omega \in [0,1]^2$ 
that correspond to resonant (blue) and rotational tori (red) as a function of $\eps$.}{histlres}{3in}

For our calculations, we will declare $\omega$ to be resonant if $\log_{10}(M)$ is more than three standard deviations below the mean of the random data of \Fig{ResonanceRandom}. Given that the cutoff \Eq{chaosCutoff} gives at least $11$-digit accuracy in $\omega$, we will use $\rho = 10^{-9}$ so that the computation of $\Delta_{m,n} < \rho$ from \Eq{minDist} has significance. 
In summary we use the cutoff 
\beq{resCutoff}
	M(\omega,10^{-9}) > 10^{2.4} = 251
\eeq
to declare that an orbit ``nonresonant''.
%

Now that we have introduced the full computational method, 
we give some information its computational complexity. 
The total runtime starting with initial conditions and 
distinguishing chaotic from nonchaotic and  resonant from rotational orbits in Matlab 2020b
using a 14-core Intel Xeon W processor at $2.5$ GHz with $64$ GB memory
is approximately $2250$ orbits per minute. 
The computation of $M(\omega,10^{-9})$ for $2250$ frequency vectors takes approximately 
32 seconds, i.e., roughly half the calculation time. 

Applying this criterion to the data in \Fig{rot1K}, separates the $80\%$ of the 
orbits that are resonant, shown in \Fig{crit2D}(a), from the remaining $20\%$ of orbits 
that are not resonant, shown in \Fig{crit2D}(b). We assume that each
of these latter orbits lies on a rotational torus, $\cT_\omega$.

\InsertFigTwo{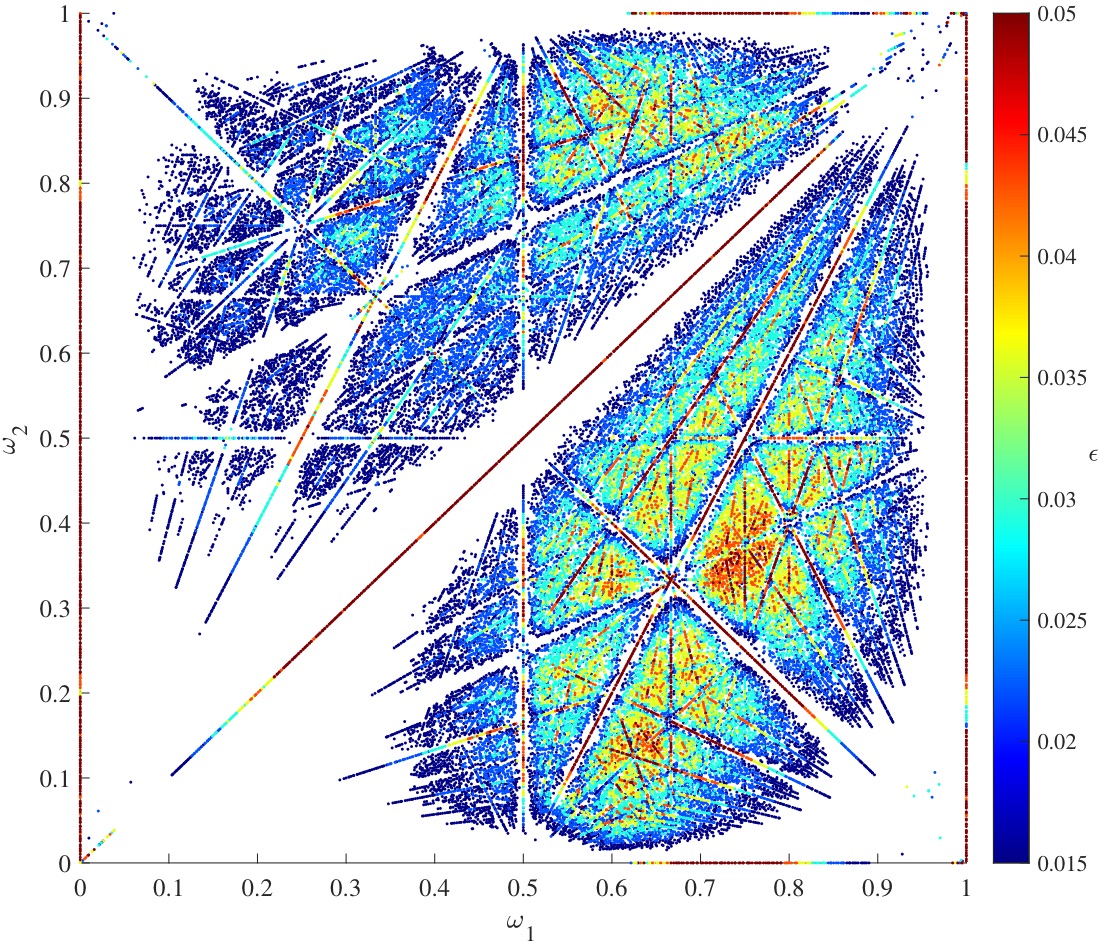}{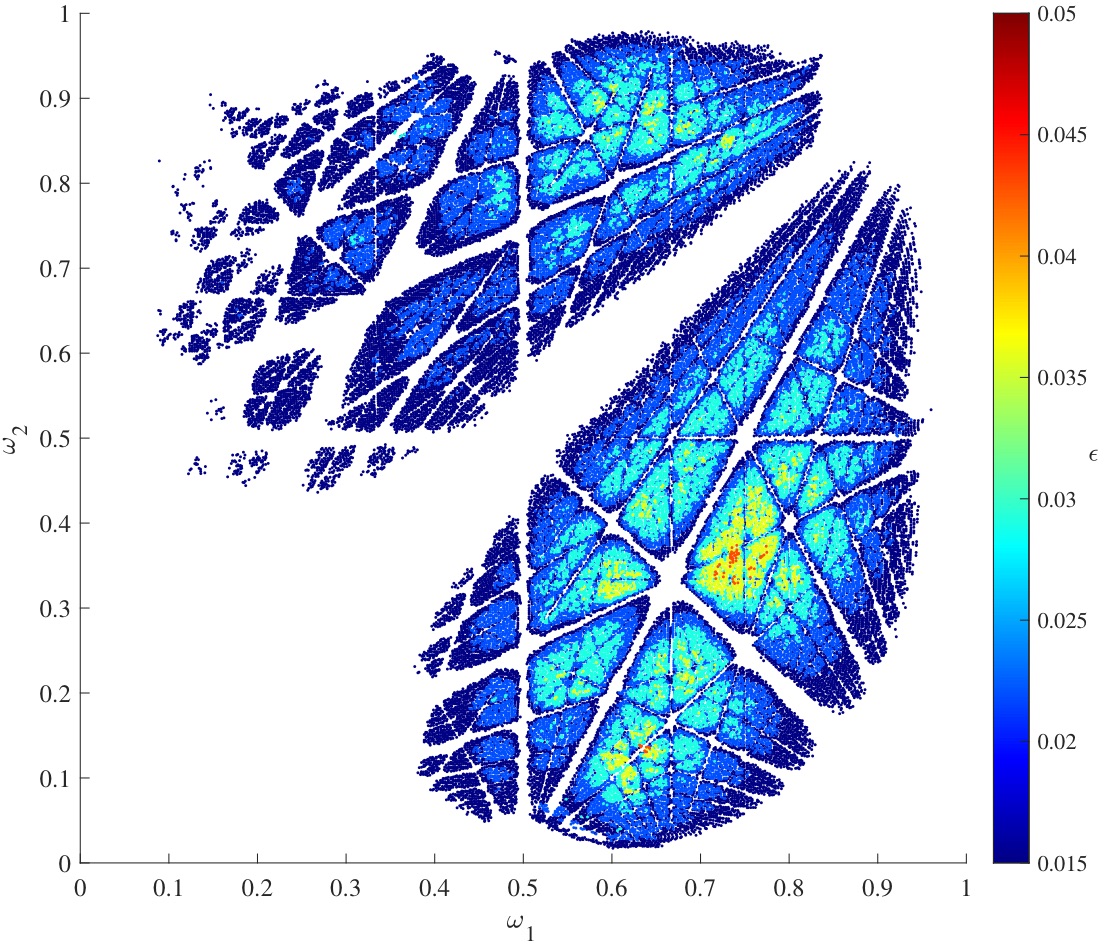}{
Frequency vectors for (a) $1,295,986$ resonant orbits and (b) $317,150$ 
nonresonant orbits, using the data of \Fig{crit3d}. 
These are distinguished by the criterion \Eq{resCutoff}. 
Values are colored  by the largest $\eps$ for which a torus 
persists at the corresponding $(y_0,\delta) \in \cP$. }{crit2D}{3in}

\section{Critical Sets}\label{sec:CriticalSets}

In this section we investigate the robustness of invariant tori as a
function of the perturbation strength $\eps$. In particular we are
interested in finding \textit{critical tori}: those on the threshold of
destruction. For smooth, two-dimensional, twist maps, an invariant
circle with rotation number $\omega \in \cD$ is critical
if has a non-smooth conjugacy to the rigid rotation $\theta \mapsto \theta + \omega \mod 1$. 
This idea was extended to the three-dimensional case in \cite{Fox16}. 
Since we are not computing the conjugacy, we instead define $\eps_c(\omega)$, following \cite{MacKay92c},
to be a value at which a rotational torus $\cT_\omega(\eps_c,y_c,\delta_c)$
\textit{breaks up}; i.e., in any neighborhood $\mathcal{N}$ of
$(y_c,\delta_c)$, there is a $\Delta\eps>0$ such that when $\eps_c <\eps< \eps_c +\Delta\eps$, there is no torus
$\cT_\omega(\eps,y,\delta)$ with the same rotation vector for any
$y,\delta \in \mathcal{N}$.\footnote
{
	One could also look for parameters at which a torus ``re-forms'', so that it does not exist when $\eps < \eps_c(\omega)$.
	Since we start from the integrable case, however, it seems sensible to first look for breakup values.
} 

When $\omega \in \cD$, and $\Omega(y,\delta)$ is a bijection and satisfies a twist condition,
then KAM theory implies that for small enough $\eps>0$ a
torus will exist for some point $(y_0,\delta)$ \cite{Cheng90b, Xia92}. On the other hand, each
resonant torus, $\cT_\omega(0,y_0,\delta)$ for $(y_0,\delta) = \Omega^{-1}(\cR)$, 
generically breaks up at $\eps = 0$. 
Define the critical set 
\beq{CritEps}
\eps_c(\omega) = \left\{ \eps: \cT_\omega(\eps,y_0,\delta) \mbox{ breaks up for some } (y_0,\delta) \in \cP \right\}.
\eeq

For the simplest, two-dimensional case, e.g., the one-parameter Chirikov
standard map, the critical set appears to be a graph over $\omega$, and
each invariant circle---once destroyed---does not reappear \cite{Marmi91}.
However critical set can be much more complicated for maps with several parameters, e.g.,
multiharmonic maps \cite{Baesens94a}, or for nontwist maps \cite{Fuchss06}.
For the standard volume-preserving map, \Eq{Parabola}, we do not actually know whether 
the critical set is simple, with only one breakup for each $\omega$.
Nevertheless, we expect that $\eps_c(\omega ) = 0$ whenever $\omega \in \cR$ and $\eps_c(\omega) > 0$ 
for $\omega \in \cD$. Since both of these sets are dense, the critical surface will 
be nowhere continuous.

\subsection{Locally Robust Tori}\label{sec:localRobust}

As an illustration of the critical set, \Fig{critsurf} shows tori that exist for 
a $100 \times 100$ grid in $\cP$, \Eq{ydeltaGrid},
for $50$ evenly spaced $\eps \in [0.015,0.045]$.
A point $(\omega,\eps)$ is shown in the figure if the parameters $(\eps,y_0,\delta)$ 
give a rotational torus, $\cT_\omega(\eps,y_0,\delta)$, using the criteria described in \Sec{Tori}. 
The upper boundary of the points shown provides a rough approximation of the critical set \Eq{CritEps}.
Of the critical tori, some have locally maximal values of 
$\eps_c(\omega)$, i.e., there is a neighborhood for which all critical tori have smaller $\eps_c$.
We will call such tori \textit{locally robust}.

\InsertFig{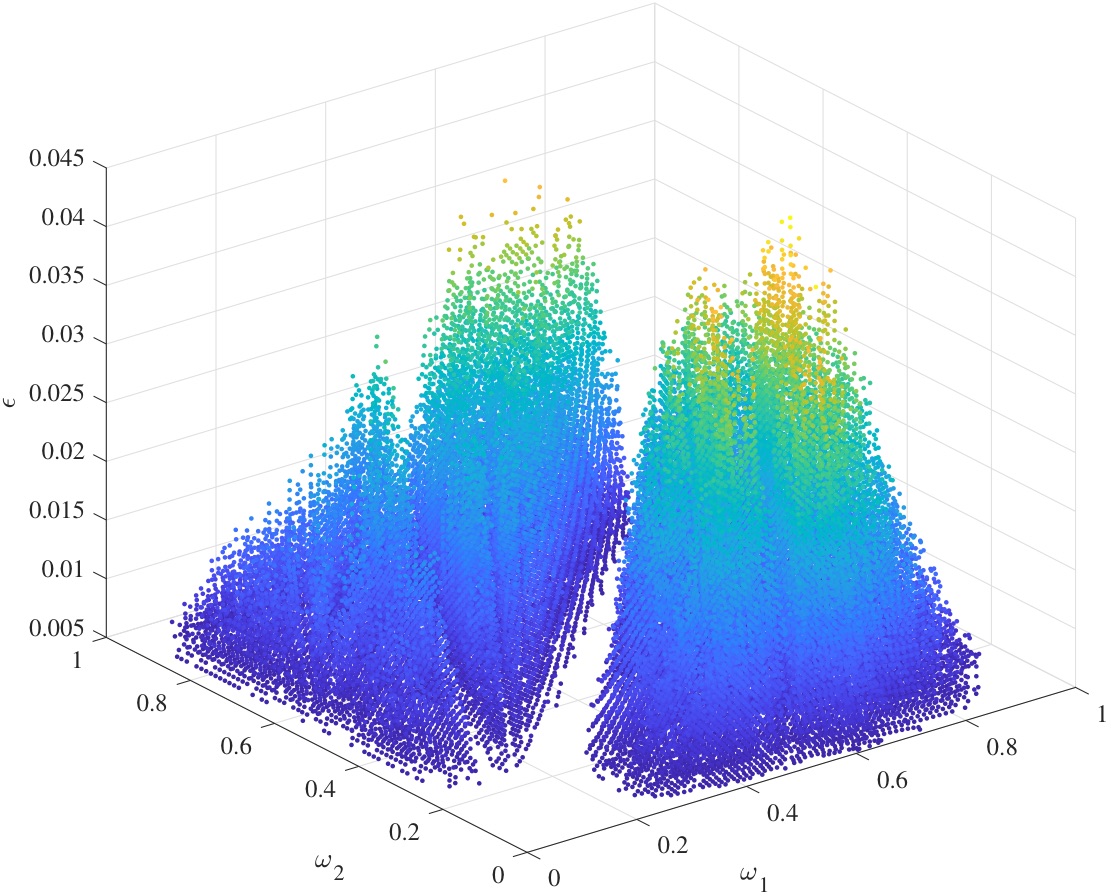}{Rotation vectors corresponding to rotational tori,
using a $100^2$ grid in the domain $\cP$ \Eq{ydeltaGrid}, and 50 $\eps$ values.
}{critsurf}{5in}

To find approximations for the locally robust tori in \Fig{critsurf}, we search for local peaks in sub-regions of the critical set
using a refinement method that does not rely on smoothness. 
In particular, for a fixed subset of $[0,1]^2$, we start at 
$\eps = 0.01$ with a $100 \times 100$ grid of corresponding points in $(y,\delta)$,   
keeping only those $(y,\delta)$ that correspond to a rotational torus in the $\omega$ region. At
each step we refine the grid for the set of parameters that have tori and 
increment $\eps \to \eps + d \eps$. Both the increment $d \eps$ and the number of grid points 
are adapted depending on the tori at the previous $\eps$.
The grid size remains $100 \times 100$ until the grid spacing is below $10^{-12}$ in the $y$ or the $\delta$ direction.
After this point we use a $10 \times 10$ grid. To choose the next $d \eps$, if more than twelve  
tori remain at $\eps$, then $d \eps \to 1.3 \cdot d \eps$. If 
between four and twelve tori remain, $d \eps$ is unchanged, and if
fewer than four tori remain, $d\eps \to d \eps/2$. Finally, if no tori are found, then 
$d \eps \to d \eps/2$, and instead of increasing $\eps$, we decrease it such that $\eps \to \eps - d \eps$.
Our procedure halts once we 
have determined a single value that is isolated on a grid of $10 \times 10$ 
points such that $y_{\max} - y_{\min}$ and $\delta_{\max} - \delta_{\min}$ are both less than $10^{-12}$, and such 
that there is no torus in the same region 
for $\eps + 10^{-14}$. Thus these values should correspond to  local 
peaks up to the corresponding accuracy in $\delta, y, \eps$. 
Dividing $[0,1]^2$ into four regions, the maxima that we 
computed are listed in \Tbl{peakdata}. 

The global maximum, found in quadrant II where $\omega \in [\tfrac12,1]\times[0,\tfrac12]$, 
agrees within $0.6\%$ for $\eps$ and $\omega$ of the results in \cite{Fox13} that were achieved
using Greene's residue method. In that paper, tori were approximated using periodic orbits
chosen on the Kim-Ostlund tree and the most robust torus was represented by a periodic orbit of period $32,316$. 
This torus was estimated to breakup at $\eps = 0.0512 \pm 0.0005$. 
For the initial conditions corresponding to this value, the weighted Birkhoff method identifies the 
orbit as chaotic for $\eps >  0.0512-0.00050005$. Thus the  identification of regular orbits 
using the Greene's residue  method in \cite{Fox13} is slightly less restrictive than the 
identification of regularity that we are using in the current paper. 

\begin{table}
\begin{center}
\begin{tabular}{ c|c c c c c } 
Quadrant		& I  & II & III & IV \\
$\omega$-Region &   $ [0,\tfrac12]^2$		& $[\tfrac12,1]\times [0,\tfrac12]$	& $[0,\tfrac12]\times[\tfrac12,1]$	& $[\tfrac12,1]^2$ \\ [0.5ex] 
 \hline
$\eps$ 		  	&  0.031282089698381	&  0.051261692234977	&   0.032740058025373   &  0.041019021169048 	\\
$\delta$ 	 	& -0.133386500670280 	& -0.334376117328629	&  	-0.743481096516467	& -0.884496372446711 	\\
$y$ 			& -0.119301311749656 	&  0.123097748168231	&	-0.162656510381853	& -0.046691232395606	\\
$\omega_1$	 	&  0.476213927381772 	&  0.734410803700126  	&	0.482238008029131	&  0.641541383714863 	\\
$\omega_2$ 	 	&  0.175290820661118 	&  0.365412254352543	&	0.781945554897404	&  0.890319673258112	\\ 
 \end{tabular}
\caption{Most robust tori in four $\omega$ regions. These were computed by 
successive refinement of a grid in $\cP$.
\label{tbl:peakdata}}
\end{center}
\end{table}

We now give a more in depth computation of the critical set, taking  slices through
\Fig{critsurf} by choosing a curve in $(y_0,\delta)$-space.
Such a slice, fixing $\delta$ and varying $y_0$, is shown in
\Fig{CrossSection}. Points here correspond to tori indicated in
\Fig{zdelta}(d) along the horizontal line at $\delta = -0.4$. Here we
plot the values of $\eps$ for which there is a torus
$\cT_\omega(\eps,y_0,-0.4)$ as $y_0$ varies. The horizontal axis is
taken to be $\omega_1$ since, when $\eps = 0$, $\omega_1 =
\Omega_1(y,-0.4) = y+\gamma$ is a bijection. 
Note that since the rotation vectors in \Fig{CrossSection} are computed on a fixed grid
in $y_0$, they are not true peak values, like those we found by refinement in \Tbl{peakdata}.

In this cross section, as well in similar cross sections for five other $\delta$ values, there appears
to only be one critical torus at any $\omega$.
In particular, the empty holes in the enlargement are sampling artifacts that go away when 
computing on a finer grid.
As the enlargement shows,
for each fixed $y_0$, the curve in $(\eps,\omega_1)$ begins as a
line for small $\eps$, but each bends as $\eps$ grows, especially for those values
approaching a visible resonant region. The curve for fixed $y_0$ does not always slope
in the same direction, as we also will see below in \Fig{sixtori}.
Gaps in these constant $y_0$ curves appear to be due to crossing such resonances.  
Unfortunately, since $\omega$ depends on $(y_0,\delta,\eps)$, the cross section is not
a simple plane $\omega_2 = $ constant. Nevertheless, since $\eps$ is
relatively small, the values of $\omega(y_0,-0.4)$ lie almost on a curve, as shown in
\Fig{CrossSection}(b), that is close to the parabola in given by
\Eq{ydeltaGrid}. This thickened curve has gaps due to resonances and has a maximum
thickness $\Delta\omega_2 \sim 0.001$.
The thickness is largest when  $0.5 < \omega_1 < 0.9$ where, as is seen in \Fig{critsurf}(a), 
tori persist for larger $\eps$. 

Note that the critical set seen in \Fig{CrossSection}(a) is visually similar to that for 
the standard map \cite{MacKay92c}, which is zero on every rational and has local maxima 
on the noble numbers. In our case the zeros occur whenever the cross section intersects a 
resonant line, and the local maxima are perhaps narrower than in the 2D case.

\InsertFigTwo{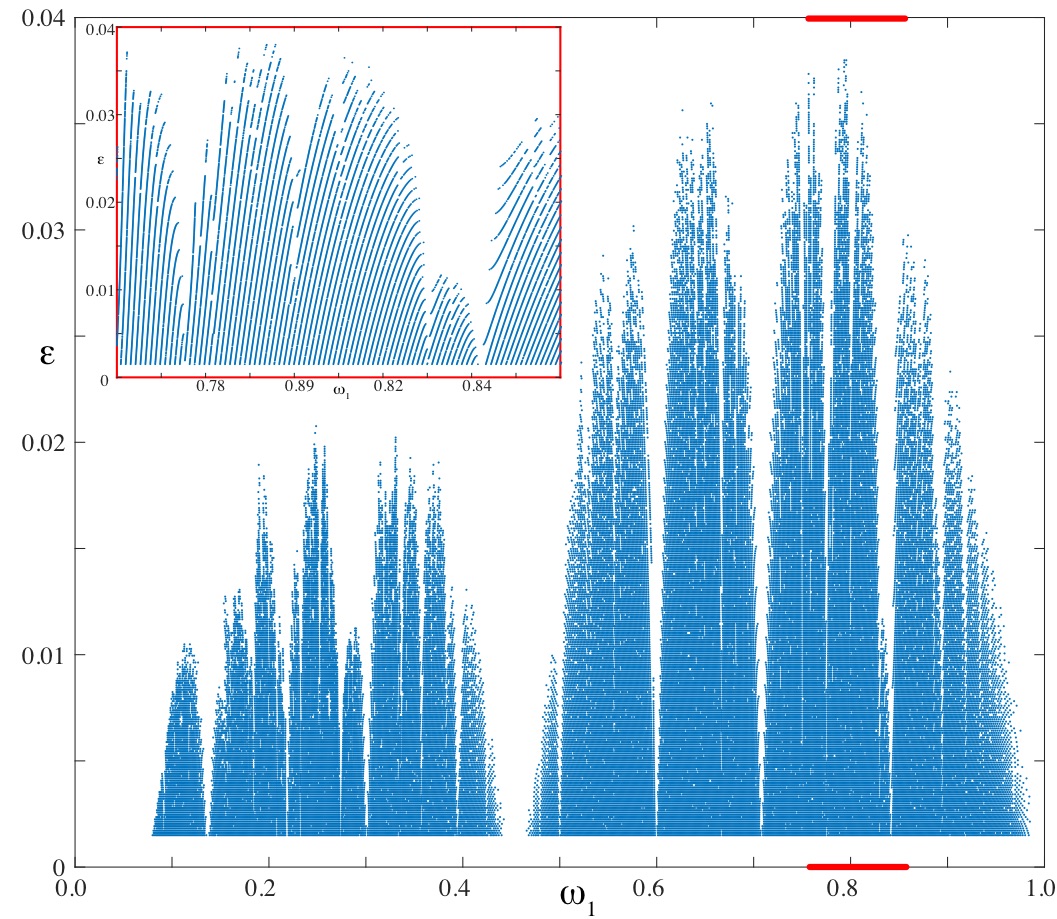}{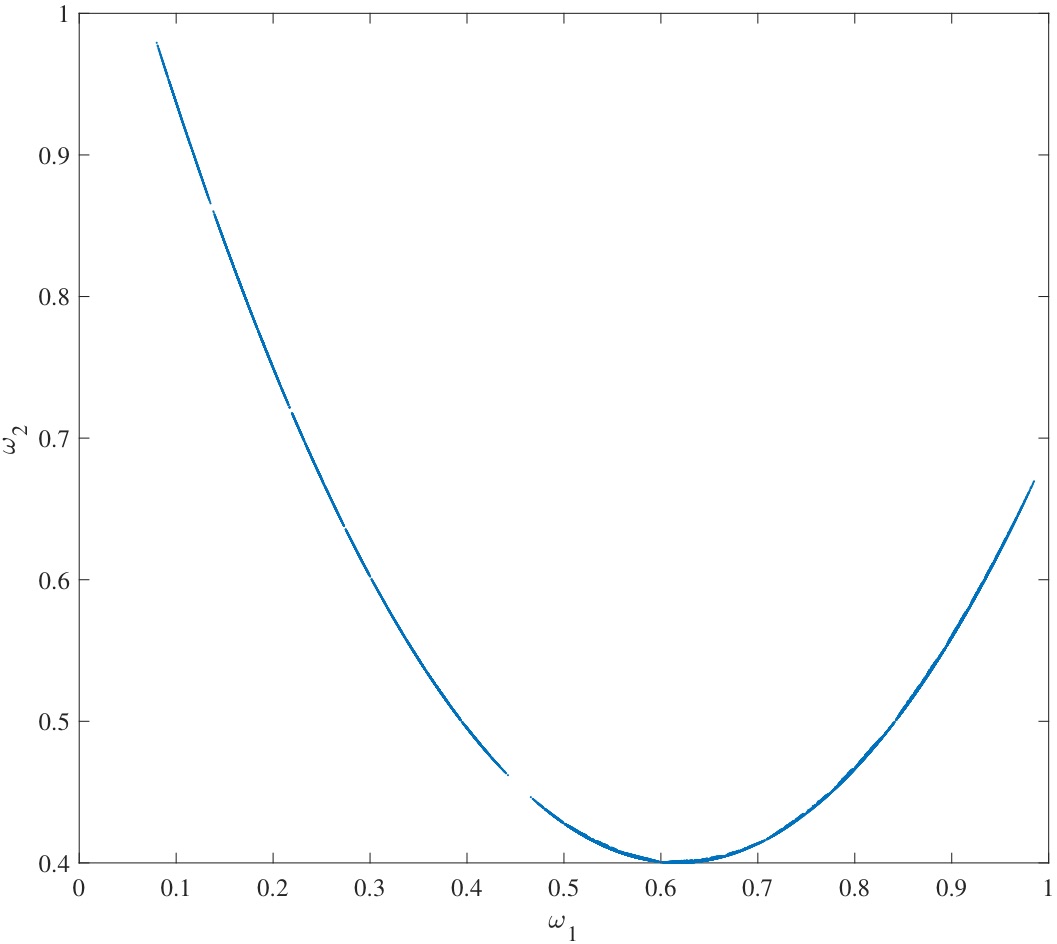}{Rotation vectors for rotational tori 
with fixed $\delta = -0.4$ and for a grid of $1000$  $y \in [-0.05 - \gamma, 1.05-\gamma]$ 
as in \Eq{ydeltaGrid} and 
a grid of $500$  $\eps \in [0.0015,.055]$. 
(a) Data projected onto the $(\omega_1,\eps)$-plane, with an enlargement for a small $\omega_1$ range
shown in the red box. 
(b) A projection of the same data onto the $\omega$-plane.}{CrossSection}{3.3in}


\subsection{Best Approximants}\label{sec:BestApprox}

Lochak \cite{Lochak92} conjectures that the robust two-tori for a symplectic map of the form \Eq{AngleActionMap} with $d=k=2$ will have rotation vectors such that $(\omega,1)$ is an integral basis for the cubic algebraic field of discriminant $49$ generated by
\beq{D49Field}
	 \alpha^3+\alpha^2 -2\alpha-1 \Rightarrow \alpha =  2\cos(2\pi/7) ,
\eeq
(see \App{CubicFields}). The field $\bQ[\alpha]$ has an integral basis $(\omega, 1)$ with
\beq{D49Freq}
	\omega =(\alpha^2-1,\alpha-1) \approx (0.554958132087371, 0.246979603717467).
\eeq
This field has the smallest discriminant amongst all totally real cubic fields. 
An alternative conjecture, probably due to Kim and Ostlund \cite{Kim86}, is that the spiral field, \Eq{SpiralField}, should give the generalization of the noble numbers for two-dimensional maps. The spiral mean is a Pisot (or PV) number: its minimal polynomial has only one root outside the unit circle \cite{Cassels57}. The spiral field is complex with discriminant $-23$, the smallest, in absolute value of all cubic fields; moreover, it is the smallest Pisot number (see \App{CubicFields}). We will consider the vector
\beq{SpiralFreq}
	\omega = (\sigma-1,\sigma^{-1}) \approx ( 0.324717957244746, 0.754877666246693) ,
\eeq
which, together with $1$ gives an integral basis for $\bQ[\sigma]$.
Finally, \cite{Tompaidis96b} considers the field with discriminant $-44$ and the minimal polynomial 
\beq{D44Field}
	\tau^3-\tau^2-\tau -1 \Rightarrow \tau \approx 1.83928675521416, 
\eeq
and chooses a basis vector for $\bQ[\tau]$ that is distinguished by having a repeated sequence in its Jacobi-Perron expansion (see \App{CubicFields}): 
\beq{D44Freq}
	\omega = (\tau-1,\tau^{-1}) = (0.83928675521416, 0.543689012692076) .
\eeq
Here we would like to provide  evidence for/against these conjectures.

As a first attempt, we investigate the Diophantine constants for the
robust frequency vectors that we have found. As discussed in
\App{Diophantine} the simultaneous Diophantine constant can be computed
by finding rational approximations $\omega \approx \frac{p}{q}$, and
computing $\Znorm{q \omega} = \|q \omega -p\|_\infty$. The sequence of
periods, $q_i$, \Eq{Periods}, defined so that $\Znorm{q_i \omega}$
decreases monotonically, corresponding to a sequence of \textit{best}
rational approximations $\omega \approx \frac{p_i}{q_i}$. A frequency
vector is Diophantine if the sequence
\[
	c_s(\omega,q_i) =  q_i \|q_i \omega\|^2_\bZ
\]
is bounded away from zero as $q_i \to \infty$ (see  \Eq{Closeness}-\eqref{eq:DiophantineConst} in \App{Diophantine}).

For example, the sequence of Diophantine constants $c_s(\omega,q_i)$ for \Eq{D49Freq} are shown in \Fig{DiophVsq}(a) 
(see the data in \Tbl{BestCubic} of \App{Diophantine}).
These appear to oscillate quasiperiodically but are bounded below, giving an estimate $c_s(\omega) \approx 0.19$
for the $D=49$ field, at least for $q_i \le 10^8$. 
It is conjectured that there is some integral basis in this field with 
Diophantine constant $\tfrac{2}{7} \approx 0.286$, and that this value is the largest possible for 
$d = 2$ \cite{Cusick74}. The corresponding Diophantine sequence for \Eq{SpiralFreq} in the $D=-23$ field and \Eq{D44Freq} in the $D=-44$ field, are also shown in 
the figure---again they are bounded away from zero as is consistent with the known Diophantine 
property of cubic fields. The dependence of $c_s$ on $q_i$ is less regular for these two vectors than for the first case, and the limit infimum appears smaller, $c_s(\omega) \approx 0.1$.

We show in \Fig{DiophVsq}(b) the sequence of Diophantine constants computed for the four 
robust tori from \Tbl{peakdata}. For these vectors, the values $c_s(\omega,q_i)$ appear to be bounded
away from zero for $q_i \lesssim 10^5$, but the values are smaller than the pure cubic vectors.
Note however, that the peak in quadrant I could have $c_s(\omega) \simeq 0.1$,
though numerical issues cause the value to drop when $q_i > 10^5$.

\InsertFigTwo{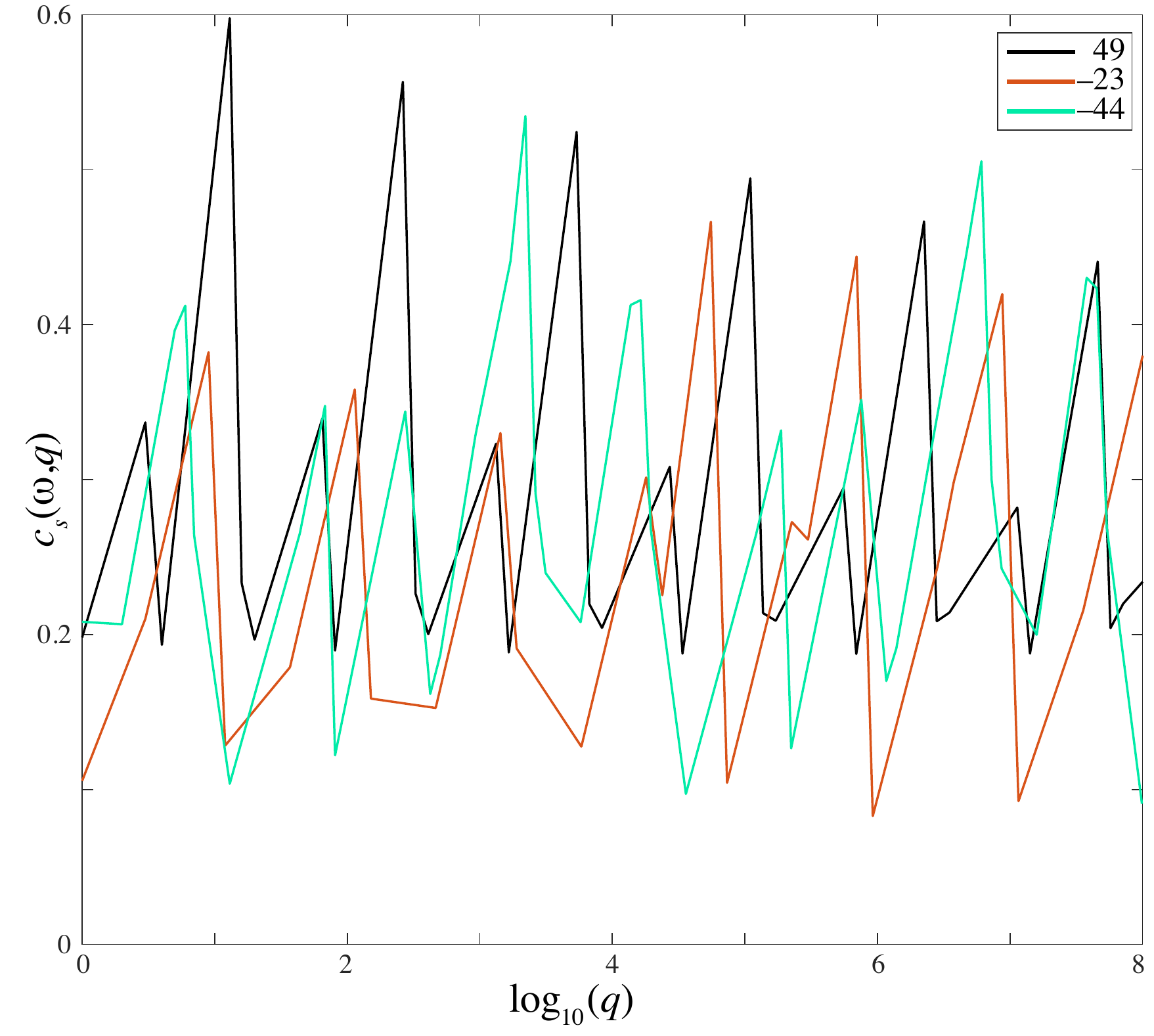}{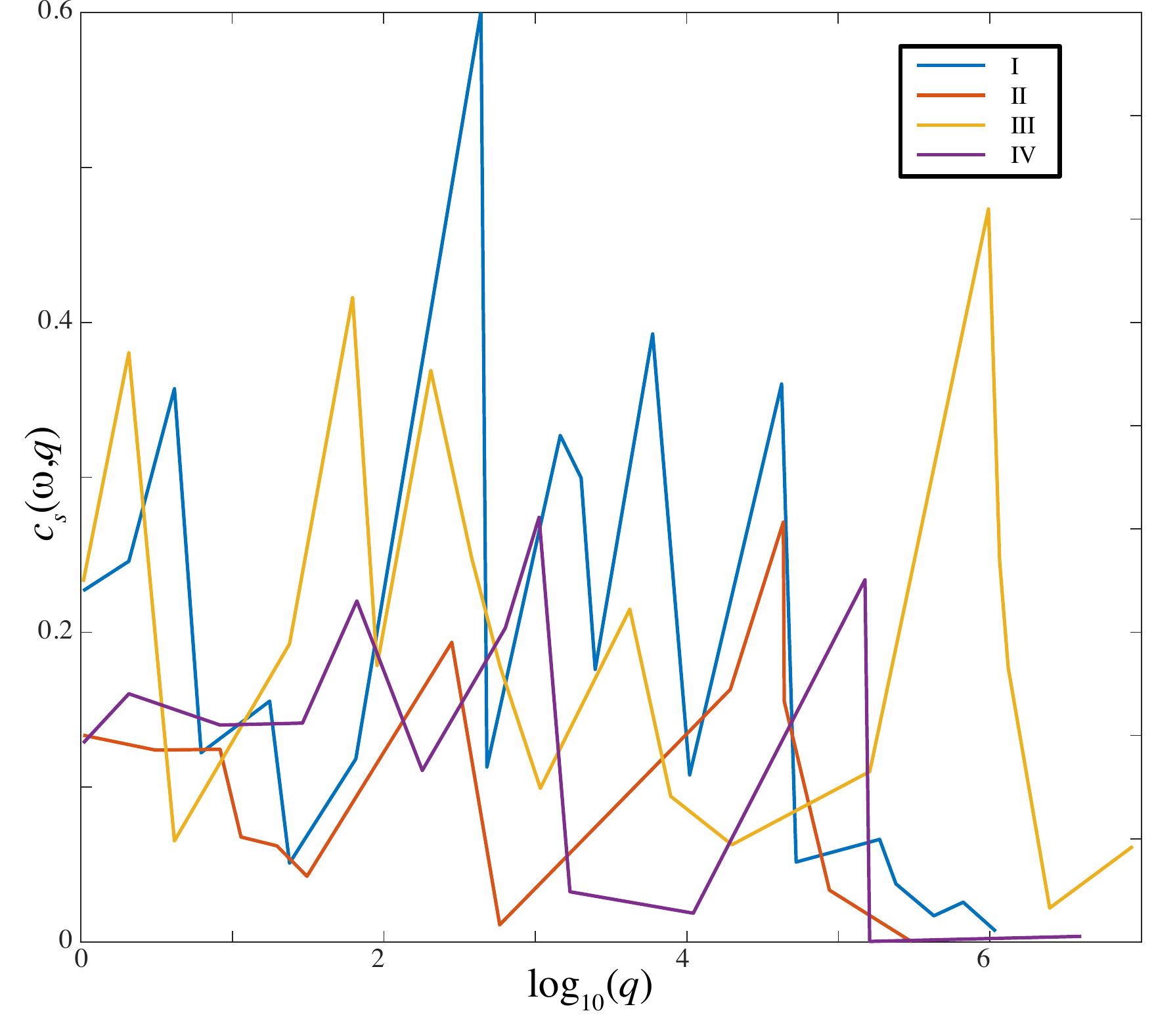}{Diophantine constants for best approximants. 
(a) Frequencies in the cubic field with discriminants $D = 49$ \Eq{D49Freq}, $D = -23$ \Eq{SpiralFreq},
and $D=-44$ \Eq{D44Freq}.
(b) The peak frequencies in the four quadrants of \Tbl{peakdata}.}{DiophVsq}{3in}

Figure \ref{fig:DiophHists} shows histograms of the Diophantine sequence
$c_s(\omega,q_i)$ for three different data sets. The first correspond to randomly
chosen $\omega \in [0,1]^2$. Note that this distribution decreases monotonically,
perhaps consistent with the expectation that there will be near rational vectors
in a random collection.
To construct a second data set, we fix the vector $(\alpha,\alpha^2,1)$, an integral
basis for the $D=49$ cubic field. Multiplication of this vector by any
matrix $A \in SL(3,\bZ)$ results in another integral basis. 
We choose four elementary matrices that generate
this group, and randomly draw a product of $50$ of these matrices to
give a set of integral bases for this field. The resulting Diophantine
constants have the distribution as shown in \Fig{DiophHists}(b). 
Note that this distribution is peaked away from $c_s = 0$.

Finally, in \Fig{DiophHists}(c), we compute $c_s(\omega,q_i)$ for the local peak
data that is obtained as follows. 
A discrete approximation to $\eps_c(\omega)$ is obtained from a
$300 \times 300$ grid  in $\cP$ and $50$ values of $\eps \in [0.005,0.045]$
--a refinement of data set shown in \Fig{critsurf}. 
Then for each bin in $\omega$ of size $0.01$ ($100^2$ bins) we select $\eps_c$ to be the
largest of the $\eps$ for tori with $\omega$ falling in that bin. 
This gives $3066$ bins that have tori with $\eps_c> 0.02$.
Since the true critical surface is not
smooth or continuous, these values almost certainly do not include the
true local maxima for all $\omega$ in each bin; indeed the largest
$\eps_c$ on this grid in $\cP$ is $0.045$, below that of the most robust torus in \Tbl{peakdata}
that we found by refinement.
Nevertheless, this process gives a set of tori that are more robust than their computed
neighbors, so that these tori are, at least, locally robust, if
not true local maxima. 

Note that the maximum of the distribution for the peak tori is shifted to $c_s \sim 0.1$, 
considerably above that of the cubic field, indicating that the peak tori are preferentially
selected to have larger Diophantine constants.
Indeed, even though all three distributions have similar means and standard deviations,
they are statistically different: a Kolmogorov-Smirnov test applied to these data
sets indicates that these distributions are distinct with $p$-values that are extremely small.

\InsertFigThree{Rand_DHist}{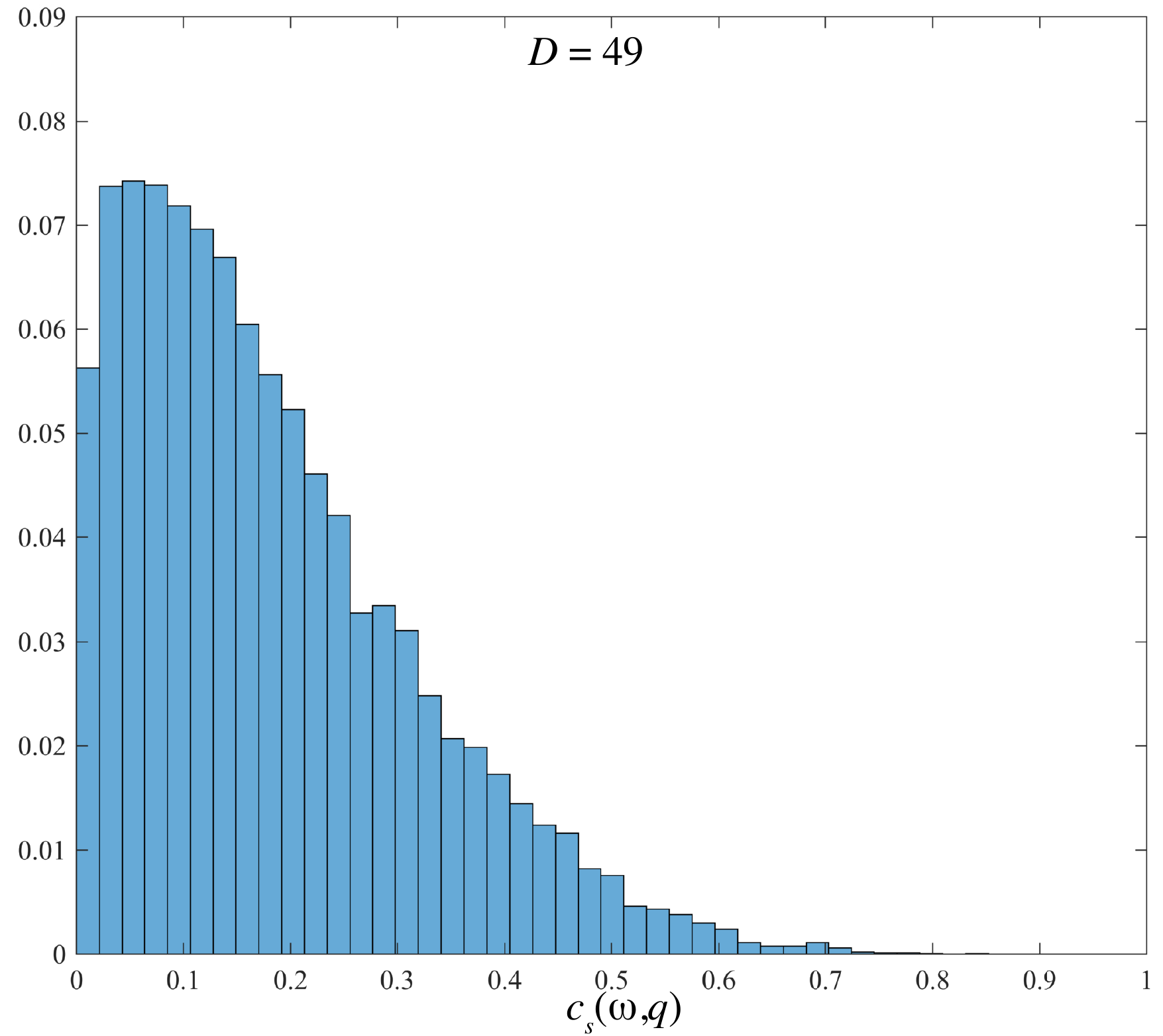}{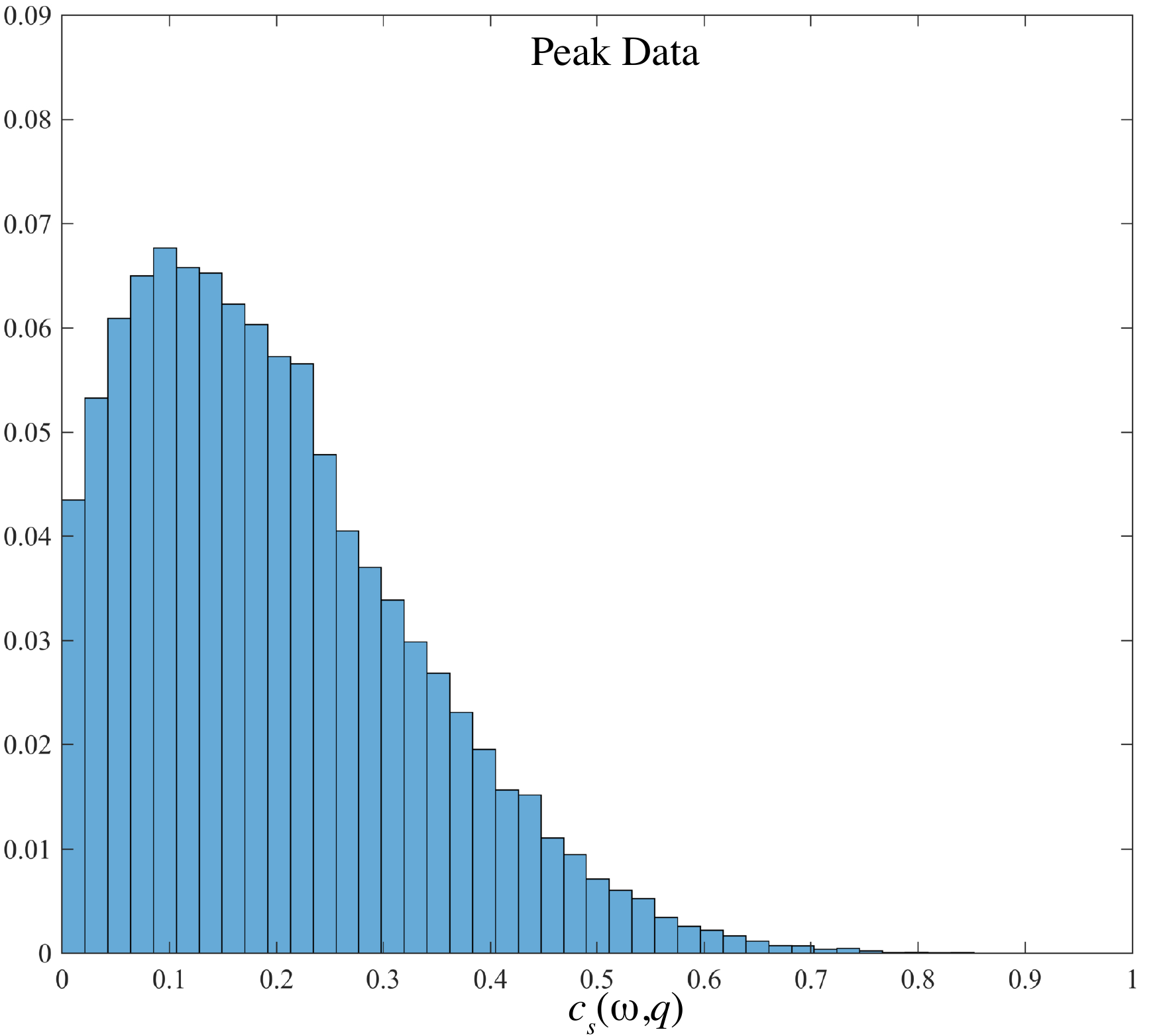}{Histograms of Diophantine constants $c_s(\omega,q_i)$ for the sequence of best periods $\{q_i\}$ for frequency vectors 
$\omega$ computed to $\|q_i\omega\|_\bZ \le \rho= 10^{-11}$. 
(a) $4000$ random frequency vectors chosen from a uniform distribution in $[0,1]^2$ (mean $=0.180$, $\sigma =  0.139$). 
(b) $4000$ vectors that are integral bases in the cubic field with discriminant $D = 49$ (mean $=0.186$, $\sigma =  0.135$) 
(c) $3066$ frequencies of tori that are locally robust on a $100^2$ grid for $\eps > 0.02$ (mean $=0.199$, $\sigma = 0.134$).
}{DiophHists}{2.25in}

The histograms shown \Fig{DiophHists} do not distinguish values of $c_s$
as a function of the period $q_i$. To do this,
\Fig{Peak_Rand} shows two-dimensional histograms with bins for both $c_s$ and
$\log_{10}(q)$. The comparison of the randomly generated data, in the left
pane, with the data from the computed peak rotation vectors in the right
pane, shows again the that the latter has significantly more values of $c_s$ bounded
away from zero.

\InsertFig{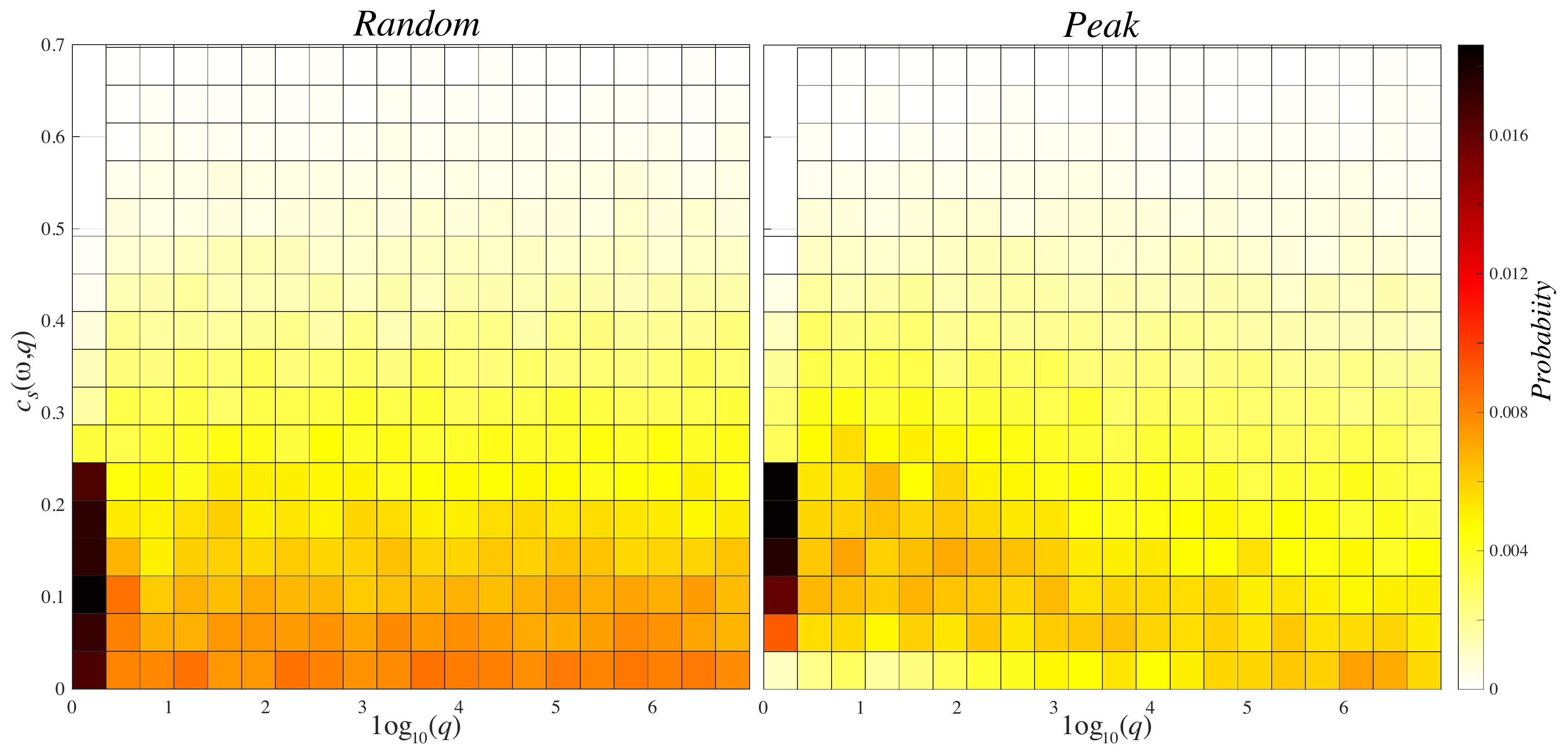}{Diophantine constants vs $\log_{10}(q)$ for 
randomly chosen $\omega$ equi-distributed in $[0,1]^2$ (left panel) and ``peak" tori 
of the volume-preserving map (right panel).}{Peak_Rand}{6in}

In conclusion, our compuations give strong evidence that the more robust tori preferentially have larger values of the simultaneous Diophantine constant, at least up to periods $q_i \sim 10^6$.

\section{Continuation}\label{sec:Continuation}
In order to further test the conjectures on which classes of frequency vectors correspond to  
most robust tori, we study here the breakup of tori $\cT_\omega$ for several fixed rotation vectors.
In particular we will find tori for vectors in the three cubic fields that have the smallest discriminants:
$D = -23$, $49$, and $-44$, recall \Eq{SpiralField}, \Eq{D49Field} and \Eq{D44Field}.

As we noted, each of these fields has properties making it a candidate to generalize of the set of
noble numbers, $\bQ[\phi]$. The field $\bQ[\sigma]$ is generated by the smallest Pisot
number and has bases that are periodic sequences in the Kim-Ostlund generalization of the Farey tree.
The field $\bQ[\alpha]$ is conjectured to have bases with the largest possible simultaneous Diophantine constant.
The field $\bQ[\tau]$ contains a basis with a period-one Jacobi-Perron sequence, one generalization of the
continued fraction, see \App{CubicFields}. 

For each case \Eq{D49Freq}, \Eq{SpiralFreq}, and \Eq{D44Freq}, give vectors for which
$(\omega,1)$ is an integral basis for the respective cubic field. 
Additionally, we will consider the permuted vector $(\omega_2,\omega_1,1)$, 
which gives an additional integral basis for the same field, so that there 
is a set of six vectors, $\omega^*$, that we study.
  
For each frequency vector, we continue the torus with respect to the 
parameter $\eps$, e.g., finding $\cT_{\omega^*}(\eps,y(\eps),\delta(\eps))$, fixing $\omega = \omega^*$.
We find the maximum $\eps$ such that the corresponding torus is nonchaotic--- 
using the $\digT$ criterion \Eq{chaosCutoff}. The torus is found using a predictor-corrector method
starting with the guess $(y,\delta) = \Omega^{-1}(\omega^*)$ at a small value of $\eps$.
Specifically, at each $\eps$, we apply the Matlab root finder {\tt fsolve} to 
find the value of $(\delta(\eps),y(\eps))$ such that the rotation vector is $\omega^*$ when computed
 using $\WB_{10^6}$. We look for $\eps_c$ such that $\digT >11$, but for which 
 $\digT<11$ when $\eps_c < \eps < \eps_c +   10^{-9}$. \Fig{rotcontinue} shows an example of
 the computation for \Eq{D44Freq}. The results for other $\omega$ values appears quite similar. 
 In particular, as seen in \Fig{rotcontinue}(b), $\digT$ drops very quickly near the critical value. 

\InsertFigTwo{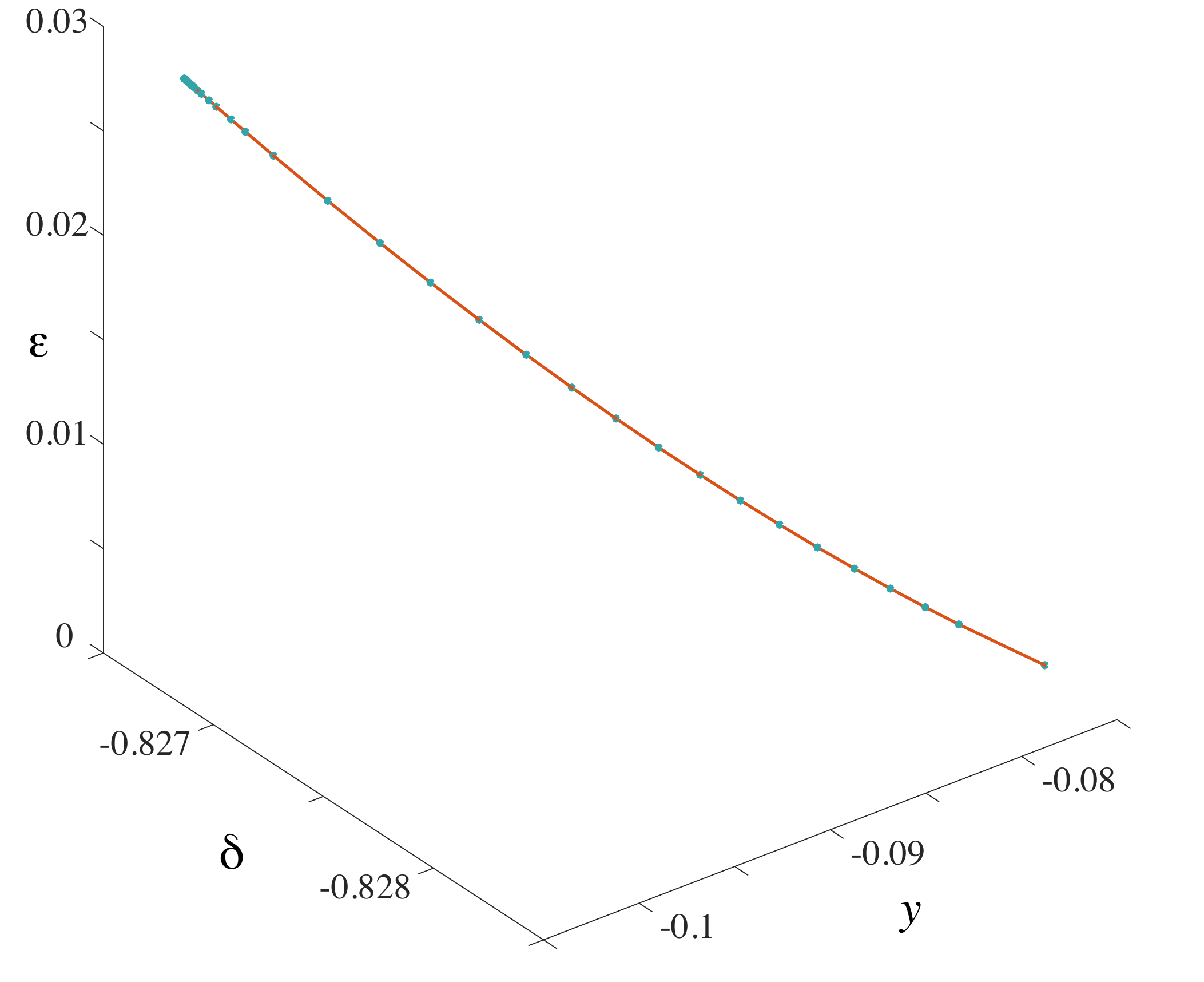}{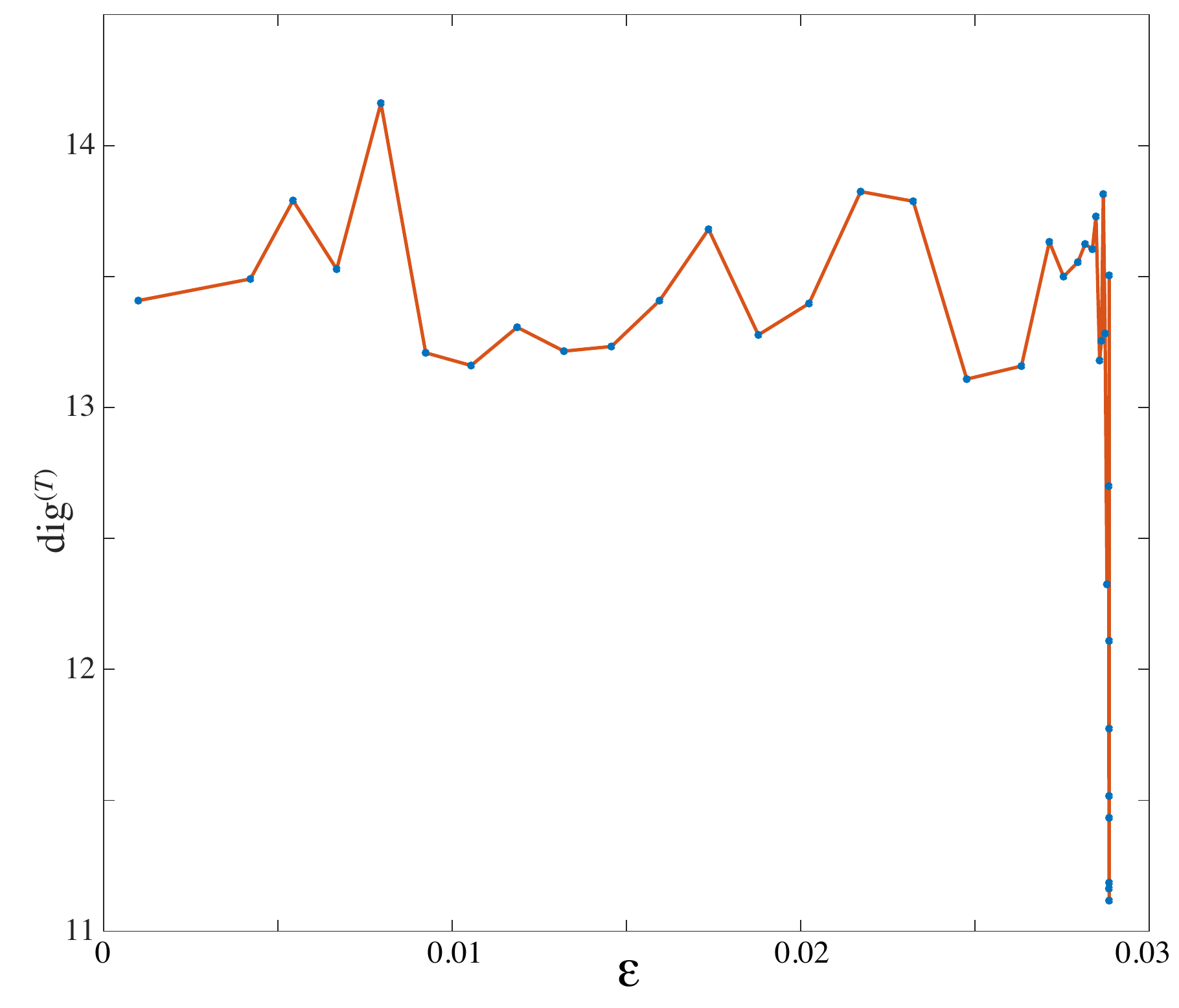}{An example of a rotational torus found using a predictor-corrector continuation method
using the weighted Birkhoff method. (a) The graph $(y(\eps),\delta(\eps))$ for a torus $\cT_\omega(\eps,y,\delta)$ with $\omega = (\tau^{-1},\tau-1)$  in the $D = -44$ field. 
(b) Number of correct digits in $\omega$ as a function of $\eps$.  
This number drops precipitously as $\eps$ approaches the critical value $0.028845453269968$, see \Tbl{cubicdata}. 
In each case, the blue dots are the computed points.}{rotcontinue}{3in}

Results of the continuation method for six frequency vectors are shown in \Fig{sixtori} and \Tbl{cubicdata}. 
None of these frequency vectors correspond to the globally most robust torus, 
nor to the quadrant maxima found in \Tbl{peakdata}. 
The first case shown, the spiral frequency $(\sigma-1,\sigma^{-1})$, was also studied in \cite{Fox13}. 
They found the stability threshold $\eps_c =  0.02590 \pm 5(10)^{-5}$ extrapolated from a sequence of 
periodic orbits up to period $31,572$. Note that our threshold in \Tbl{cubicdata}, 
$\eps_c \simeq 0.02573 = 0.02590 - 1.7(10)^{-4}$, is again slightly more conservative than that 
given by Greene's residue.
 
\InsertFig{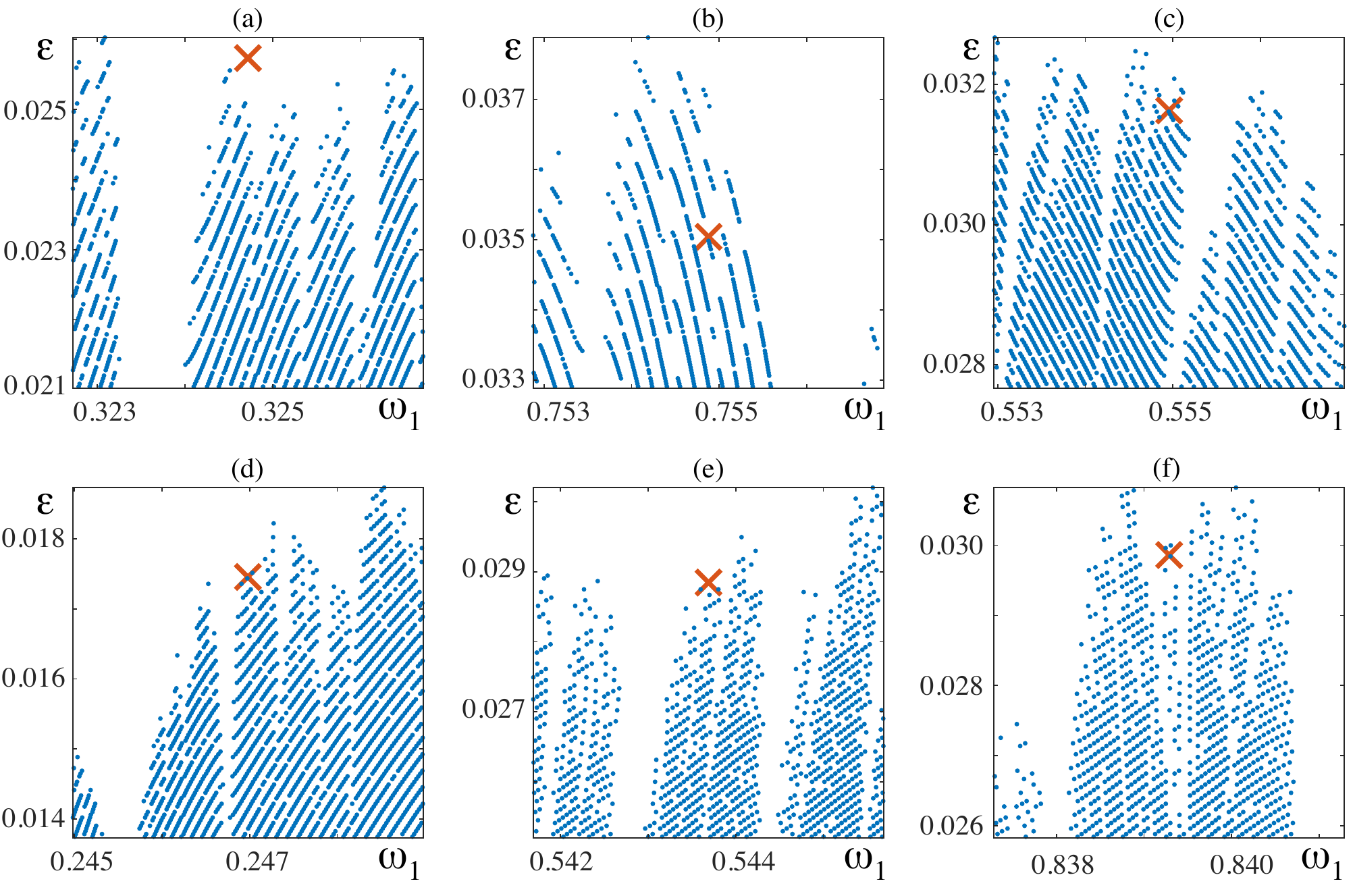}{Tori in the $(\omega_1,\eps)$ plane,
near $\eps_c(\omega^*)$ for the six vectors $\omega^*$ in \Tbl{cubicdata}.
Each point represents a torus at a fixed $\delta = \delta_c$ 
that passes through the corresponding critical point. In each figure, the (red) X marks the location 
of the associated critical torus. The frequencies corresponding to 
nearby tori (blue) are computed on a grid with spacing $10^{-4}$ in $\omega_1$ and 
$4(10)^{-5}$ in $\eps$. }{sixtori}{6in}

While it is not possible to compute whether these tori are locally most robust on 
arbitrarily small neighborhoods, we can quantify the degree to which they fail 
to be local maxima on a fixed small interval in $\omega_1$, fixing a cross-section $\delta = \delta_c$ 
that passes through the corresponding critical point.  
We compute rotational tori for values of $y$ near $y_c$ such that the 
spacing of $\omega_1$ values is approximately $10^{-4}$ (i.e., using the $\eps = 0$ approximation for 
$y$ and $\delta$), 
and the $\eps$ spacing is $4(10)^{-5}$ with  $\eps \le \eps_c + 0.005$. 
Of these, we consider only tori with $|\omega_1-\omega_1^*|<0.002$, 
and which are more robust than $\omega^*$. 
These tori correspond to the points in \Fig{sixtori} that lie above $\eps_c(\omega^*)$, which
is indicated by the (red) X in each panel. 

\begin{table}
\begin{center}
\begin{tabular}{l|c|r r r r r r} 
 &\multicolumn{1}{c}  {$\omega^*$} & \multicolumn{1}{c}{$\eps_c$}& 
		\multicolumn{1}{c}{$y_c$ } &\multicolumn{1}{c}{$\delta_c$} & \multicolumn{1}{c}{ \#M.R.} \\
\hline 
(a) & $(\sigma-1,\sigma^{-1})$    & $0.025731358271922$ & $-0.300341913511639$ & $-0.581991952776833$   & 3\\
(b) &$(\sigma^{-1},\sigma-1)$    & $0.035042379103690$ & $0.137249321586741$  & $-0.285775646047323$   & 180\\
(c) &$(\alpha^2-1,\alpha-1)$	    & $0.031629688390353$ & $-0.046265357816195$ & $-0.237775229395970$   & 61\\
(d) &$(\alpha-1,\alpha^2-1)$     & $0.017453913097431$ & $-0.374459102422933$ & $-0.279201189795316$   & 92 \\
(e) &$(\tau-1, \tau^{-1})$		& $0.029861717573837$ & $-0.444144360895140$ & $0.203236055553548$	  & 47\\
(f) &$(\tau^{-1},\tau-1)$ 		& $0.028845453269968$ & $-0.826759694950616$ & $-0.103754031724428$	  & 44\\
\hline
\end{tabular}
\caption{Critical parameters for tori in cubic fields with discriminants $-23$, $49$ and $-44$, recall
\Eq{SpiralField}, \Eq{D49Freq}, and \Eq{D44Freq}. The labels (a)-(f) indicate
the corresponding panel in \Fig{sixtori}.
For each field, two rotation vectors $\omega^*$ are chosen, related by permutation.
The torus $\cT_{\omega^*}(\eps,y,\delta)$  breaks up at $\eps_c$ and is located at $(y_c,\delta_c)$.
The final column shows the number of tori that are more robust among the tori computed in an interval
$|\omega_1-\omega_1^*|<0.002$ with $\omega_1$ spacing of $10^{-4}$
 and  $\eps<\eps_c+0.005$ with spacing of $4(10)^{-5}$. 
\label{tbl:cubicdata}}
\end{center}
\end{table}

One measure of local robustness is the distance, $\Delta\omega = |\omega_1-\omega_1^*|$,
to the closest, more robust torus.
The first spiral mean vector, (a) in the table and figure, 
has a distance of $\Delta\omega \sim 2(10)^{-3}$; this makes it 
five-times more robust than than (e), the first $D = -44$ torus, and 
more than 25-times more robust than all others.
By this measure the frequency (b) in $D=-23$, 
is least robust since there is a more robust torus within $\Delta\omega \sim 4(10)^{-6}$.
This measure best matches what is observed by eye in \Fig{sixtori}, 
where similarly (a) and (e) appear the most locally robust within a single peak of the
each panel.

Another measure of robustness is the distance $\Delta\eps = \eps_c(\omega)-\eps_c(\omega^*)$ 
for the most robust computed torus in the interval.
By this measure, the first vector, (a), is more robust in the sense that for this torus
$\Delta\eps \sim 3(10)^{-4}$, and this value is at least three-times smaller than that for any of the  other tori.
Torus (b), the second spiral mean case, is the least robust with $\Delta\eps$ more than nine-times larger 
than the value for (a). 

Finally, the last column of \Tbl{cubicdata}, labeled \#M.R., lists the number of more robust tori within 
the computed frequency interval. Again by this measure the first spiral mean, (a), is the most 
robust since the number is least 14-times smaller than any other number. The second spiral 
mean, (b), is the least robust with $180$ nearby, more robust tori. 
 
From these results, though none of the six vectors considered is locally robust over a large range,
both (a) and (e) could be considered locally robust, and 
the first spiral mean case, (a), is the most robust over the range we have considered.
This perhaps provides some indication
that the spiral field remains the best candidate to generalize the noble numbers.
Nevertheless, since there are a countably infinite number of vectors in each field, 
we cannot rule out that another representative would behave more robustly. 



\section{Conclusions}\label{sec:Conclusions}

In this paper we develop criteria to distinguish orbits that lie
on rotational two-tori from those that are chaotic or resonant, as well as to compute 
rotation vectors for the rotational two-tori. Our
model is a three-dimensional, volume-preserving map \Eq{AngleActionMap}
with frequency map and force \Eq{Parabola}, and the primary tools are the
weighted Birkhoff average \Eq{WB} and an algorithm for calculating resonance 
order, \Alg{ResonanceOrder}. 

To distinguish chaotic from regular orbits and 
calculate rotation vectors, we use the weighted Birkhoff average. 
Computation of the rotation vector of the tori to at least 11-digit accuracy
required $2T = 2(10)^6$ iterates of the map---the second half of the
iterates being used to estimate the error \Eq{digits}. However in most
cases as we saw in \Fig{histc8}, many fewer iterates---say
$3(10)^4$---would be sufficient to obtain this accuracy and to
distinguish regular from chaotic dynamics.  Moreover, as we previously
showed in \cite{Sander20}, the weighted Birkhoff average is more
efficient than other techniques, such as fast Lyapunov methods, for this
distinction.

To distinguish tori from resonances, it was important to use ``linear
approximations" rather than the often-used simultaneous approximations to the vector
$\omega$: we look for the closest resonant line \Eq{Resonance} to the
computed rotation number in the sense of Euclidean distance
\Eq{minDist}. Unlike the case of a single frequency (where the
Stern-Brocot tree is the optimal method \cite{Sander20}), there seems to
be no general theory that gives a ``fast'' method for determining
optimal linear approximations. The general theory of optimal resonance
order is not completely understood, though the scaling of resonance
order with tolerance that is seen in \Fig{nonChaoticFreqs} is what would
be expected from the theorems of Dirichlet and Minkowski (see
\App{Diophantine}).

To compute resonance order to a precision $\rho$ we use a brute
force method, recall \Alg{ResonanceOrder}, and unfortunately, this is a
substantial part of the computational cost in our method (about 50\%
of the effort). Nevertheless, is important to find such linear approximations to
eliminate dynamical resonances; such orbits lie on regular tori that are not
``rotational'', instead these enclose isolated invariant circles of the
map (if they exist \cite{Dullin12a}).

After finding rotational tori, in \Sec{CriticalSets} and \Sec{Continuation}
we use a variety of methods, including simultaneous approximations and parameter continuation
to study  the most robust tori, 
and to test conjectures regarding the robustness of tori with rotation vectors in three cubic fields.
Since the number of tori studied here is relatively small 
(compared to the number tested to find the locally robust rotational tori), 
we also use, in \Sec{BestApprox}, a brute force method to compute the sequence of best simultaneous
approximations to $\omega$, e.g., the sequence of ``periods'' \Eq{Periods} of a rotation vector.
Note that, by contrast with linear approximations, there is a fast method for computing the simultaneous
approximation sequence \cite{Clarkson97a}. 

Our results indicate that frequencies of robust tori are discernibly different from
random frequency vectors, see \Fig{DiophHists}, but as in \cite{Fox13}, we 
have been unable to extract a simple number-theoretic property for the globally most robust tori
and the Diophantine constant sequence does not seem to have a simple behavior, \Fig{DiophVsq}.
There is only weak evidence for local robustness of conjectured, low-discriminant
cubic fields, and in particular while some of the elements in the three conjectured cubic fields 
may be locally robust, none of them seem to be associated with the most robust tori 
in the same way that the noble numbers are associated with most robust and locally 
robust invariant circles for 2D maps. 

Our results for computing tori should be compared with previous
techniques that used periodic orbits to approximate the torus and
Greene's residue criterion to estimate their breakup. Greene's residue
is certainly the optimal method for area-preserving maps with a
time-reversal symmetry: it can easily give highly accurate breakup
thresholds for an invariant circle \cite{MacKay83}.  As was shown in
\cite{Fox13}, this idea can be generalized to the volume-preserving case
studied in the current paper because it does have a time-reversal
symmetry, which permits the computation of sequences of symmetric
periodic orbits. In \cite{Fox13}, orbits with periods of order $10^4$ were
used to obtain estimates of breakup thresholds with a relative error in
$\eps_c$ that was estimated to be $0.002$. The results obtained in the current
paper show that these thresholds were slight over-estimates of
$\eps_c(\omega)$, and allow us to refine the breakup threshold to a
relative error of about $4(10)^{-8}$. 
This relies on the ad hoc criterion \Eq{chaosCutoff}, that declares that the
orbit has become chaotic when the accuracy of the weighted Birkhoff
average drops below our threshold of $11$-digits (of course, there is a
similar, ad hoc threshold for the residue criterion). Note, however that the rapid decrease of 
$\digT$ over a narrow parameter range, as seen in \Fig{rotcontinue}, 
is a clear signal of the torus destruction. 

Another advantage of the weighted Birkhoff method used here is that it
can be applied to asymmetric maps (as we did for 2D in \cite{Sander20}),
and does not rely on any sophisticated method to find periodic orbits
that are simultaneous rational approximations to a given incommensurate
frequency vector. Moreover, even though the residue method works very
well in 2D and is successful in 3D, it has not led to accurate methods
that can estimate the breakup thresholds for tori in 4D symplectic maps. One of the problems for the 4D
case is that  there are multiple ``partial traces'' of a symplectic
matrix needed to define a stability threshold. Such considerations are
irrelevant for weighted Birkhoff averages. 

In the future, we plan to continue the study of two-tori for maps in three and possibly four dimensions.  
The fast convergence for the weighted Birkhoff method makes it well suited for extended precision computations
as noted in \cite{Das17}.
We have successfully performed test cases for one- and two-dimensional tori in phase spaces 
of dimensions one, two, and three. In future work, we plan to explore such high precision calculations 
in the hope that extracting more digits will lead to a better
number theoretic understanding of the properties of the rotation vectors for robust tori.

\cleardoublepage
\pagebreak
\appendix

\section{Diophantine Constants}\label{app:Diophantine}
There are two common ways in which a vector $\omega \in\bR^d$ can be approximated by rationals. The first, \textit{simultaneous} approximation, seeks a vector  $(p,q) \in \bZ^{d}\times \bN$ that corresponds to a nearby rational, i.e., $\omega \approx \frac{p}{q}$. The second, \textit{linear} approximation, seeks an integer vector $(m,n) \in \bZ^d \setminus \{0\} \times \bZ$, that corresponds to a nearby resonance, i.e., $m \cdot \omega -n \approx 0$.

Define the pseudo-norm
\beq{Znorm}
	\Znorm{\omega} = \inf_{p \in \bZ^d} \|\omega - p\|_\infty
\eeq
that computes the distance to the nearest point on an integer grid.
Lochak \cite{Lochak92} defines the \textit{periods}, $q_i$ of $\omega$ as the sequence of positive integers so that, $q_0 = 1$, and $\forall q < q_{i+1}$,
\beq{Periods}
	\Znorm{q_{i}\omega} \le \Znorm{ q\omega}.
\eeq
Thus if $p_i$ is the integer vector such that $\Znorm{q_{i}\omega} = \|q_i\omega -p_i\|_\infty$, then $\frac{p_i}{q_i}$ are (strong) \textit{best approximants} of $\omega$. 

Similarly we can define a set of nearest resonances to $\omega$ as a sequence of nonzero integer vectors $m_i$ so that $m_0 =(1,...,1)$ and whenever $0 < \|m\| \le \|m_{i+1}\|$ then
\[
	\Znorm{m_i \cdot \omega} \le \Znorm{m\cdot \omega}, \quad m,m_i \in \bZ^d \setminus \{0\}
\]
The sequence of resonance orders, $\|m_i\|_1 = M_i$ thus obtained is unique, even though the resonant sequence itself may not be.

Theorems of Dirichlet and Minkowski, give a bound on the goodness of these approximations:

\begin{thm}[\cite{Baker84,Cassels57,Schmidt91}]\label{thm:Minkowski} For any $\omega \in \bR^d$ and for any $K > 0$, there exist  $(m, n) \in \bZ^d\setminus \{0\} \times \bZ$ with $\|m\|_\infty \le K$ such that 
\beq{MinkowskiRes}
	|m \cdot \omega - n| < \frac{1}{K^d}.
\eeq
Similarly whenever at least one $\omega_i$ is irrational, then for any $Q>0$,  there exist $0 < |q| \le Q$ and $p \in \bZ^d$ such that 
\beq{MinkowskiRat}
	\|q\omega - p\|^d_\infty < \frac{1}{Q}.
\eeq
\end{thm}

\noindent
In \Eq{MinkowskiRes} we say that $(m,n)$ is a ``near'' resonance for $\omega$ and, for \Eq{MinkowskiRat}, that the vector $p/q$ is a ``good'' rational approximation to $\omega$. Note the complementary placement of the $d^{th}$ power in these two expressions.

Based on these complementary notions of approximation, there are also two senses in which 
$\omega$ can be strongly ``irrational", or Diophantine. Defining the linear and simultaneous ``closeness" parameters
\bsplit{Closeness}
    c_l(\omega,m) &= \|m\|^d_\infty \Znorm{m \cdot \omega} , \\
	c_s(\omega, q) &=  q \|q \omega\|^d_\bZ ,
\esplit
then the associated Diophantine constants are
\bsplit{DiophantineConst}
	c_l(\omega) &= \liminf_{\|m\|_\infty \to \infty} c_l(\omega,m) , \\
	c_s(\omega) & = \liminf_{q \to \infty}  c_s(\omega, q) .
\esplit
A vector $\omega$ is (linear, simultaneous) Diophantine if $c_{l,s}(\omega) > 0$. A theorem of Dirichlet implies that if $\theta$ is an algebraic irrational of degree $d+1$, then for $\omega = (\theta,\theta^2,\ldots,\theta^d)$, $c_s(\omega) > 0$; thus the vector $\omega$ is (simultaneous) Diophantine. For example, every quadratic irrational is Diophantine for $d=1$.

Note that, by \Th{Minkowski}, $c_l, c_s \le 1$.  When $d = 1$ these constants are trivially equal 
and it is known that for any $\omega$, $c_s(\omega) \le c_s(\phi) = \frac{1}{\sqrt{5}}$,
where $\phi$ is the golden mean \cite[Thm. 194]{HardyWright79}. As noted in \cite{Lochak92}, for $d = 2$ it has been proven by Davenport that the upper bounds on $c_l$ and $c_s$ over $\omega \in \bR^2$ are the same, and that these upper bounds are at least $\tfrac27$.  Furthermore, Adams showed that there are integral bases for cubic fields for which $c_s(\omega) = \tfrac27$. It been conjectured that for any $\omega$, $c_s \le \tfrac27$, and that the only numbers for which $c_s$ is near $\tfrac27$ are integral bases of a real cubic field \cite{Cusick74}. Indeed Cusick conjectures that there is an integral basis for the discriminant $49$ field \Eq{D49Field} that achieves this value.

Computation of the sequence of periods \Eq{MinkowskiRat} can be done efficiently using algorithms developed by Clarkson \cite{Clarkson97a}, and this would be especially important if high precision computations are required to attempt to estimate the asymptotic Diophantine constant. 
We simply use the brute force method of computing $\|q\omega\|_Z$ for each natural number up to some $Q$.
The resulting sequence of periods and approximations to the Diophantine constant for the vector \Eq{D49Freq} 
are shown in \Tbl{BestCubic}.  A comparison between several cubic irrationals is shown in \Fig{DiophVsq}(a).

Note that if we know $\omega$ to precision $\rho$,  then $c_s(\omega,q)$ can 
be computed with precision $q \rho$. Using the cutoff $\digT = 11$, this implies that periods must be limited 
so that $q \ll 10^{11}$, e.g., $q \lesssim 10^7$ so that $c_s$ can be computed with 4 digit accuracy.

\begin{table}[htp]
\begin{center}
\begin{tabular}{r r r | c c}
$p_1$ & $p_2$ & $q$ & $\|q\omega\|_Z$ & $c_s(\omega,q)$\\
1 &	0 &	1& 	0.445041867912628&	0.198062264195161\\ 
2 &	1 &	3& 	0.335125603737885&	0.336927510842046\\ 
2 &	1 &	4& 	0.219832528349486&	0.193305362082111\\ 
7 &	3 &	13& 	0.214455717135830&	0.597886309959160\\ 
9 &	4 &	16& 	0.120669886602055&	0.232979544520846\\ 
11 &	5 &	20& 	0.099162641747430&	0.196664590366583\\ 
36 &	16 &	65& 	0.072278585679150&	0.339572606605588\\ 
45 &	20 &	81& 	0.048391300922901&	0.189679158405874\\ 
146 &	65 &	263& 	0.046011261021278&	0.556780505022018\\ 
182 &	81 &	328& 	0.026267324657880&	0.226310929055850\\ 
227 &	101 &	409& 	0.022123976265050&	0.200193363242583\\ 
737 &	328 &	1328& 	0.015600587970539&	0.323206442195226\\ 
919 &	409 &	1656& 	0.010666736687313&	0.188418473697497\\ 
2984 &	1328 &	5377& 	0.009876233796604&	0.524472547765830\\ 
3721 &	1656 &	6705& 	0.005724354173708&	0.219710986884052\\ 
4640 &	2065 &	8361& 	0.004942382513946&	0.204235358627254\\ 
15066 &	6705 &	 27148& 	0.003369907963133&	0.308300280752348\\ 
18787 &	8361 &	 33853& 	0.002354446209210&	0.187661294078270\\ 
61001 &	 27148 &	109920& 	0.002120956116414&	0.494470156865191\\ 
76067 &	 33853 &	137068& 	0.001248951841262&	0.213809728033248\\ 
\end{tabular}
\end{center}
\caption{Best approximants for \Eq{D49Freq}, a vector in the $D = 49$ cubic field.}
\label{tbl:BestCubic}
\end{table}

\section{Cubic fields and Jacobi-Perron}\label{app:CubicFields}
An algebraic field is an extension of $\bQ$ to include some family of algebraic numbers.
Consider the monic polynomial
\beq{Polynomial}
	p(x) = x^3 - kx^2 - lx -m
\eeq
for $k,l,m\in \bZ$ and $m \neq 0$. A root, $\tau$ of such a polynomial
is an algebraic integer and generates a cubic field $\bQ[\tau] =
\{a+b\tau+c\tau^2: a,b,c \in \bQ\}$. The integers in such a field
correspond to the restriction of $a,b,c \in \bZ$. An integral basis of
such a field is a vector $(1,\tau,\sigma)$ that generates the integers. 
The fields can be characterized by the discriminant, $D= k^2l^2 - 4k^3m
+ 4l^3 - 18klm - 27m^2$ of the polynomial \Eq{Polynomial}. When $D <
0$ there are two complex roots and when $D>0$ all the roots are real.

The cubic field with the smallest $|D|$ is that of the spiral mean, which has
$D = -23$, the minimal polynomial \Eq{Polynomial} with $(k,l,m)=
(0,1,1)$ with real root $\sigma$, and integral basis
$(1,\sigma,\sigma^2)$. Siegel showed this root is the smallest Pisot number (an algebraic number
that  is the unique root of its minimal polynomial outside the unit circle) \cite{Siegel44}. We have used
the fact that, for the vector $\omega$ in \Eq{SpiralField},  $(1,\omega)$ forms an integral basis for
$\bQ[\sigma]$. An alternative polynomial for this field has $(k,l,m) =
(-1,0,1)$. The real root in this case is $\sigma^{-1}$.

When $m = 1$, and $\tau$ is a Pisot number, Tompaidis noted that there is
an integral basis of the cubic field for which the Jacobi-Perron
algorithm (JPA), a generalization of the continued fraction, is periodic.
Given a vector $\omega \in [0,1]^2$ the
JPA generates a sequence of integer vectors $r_i =
(p_i,q_i) \in \bN^3$ so that $p_i/q_i$ is a rational approximation to
$\omega$ and these vectors obey a recursion
\[
	r_{n+1}= k_{n+1} r_n + l_{n+1} r_{n-1} + r_{n-2} .
\]
The coefficients $(k_n,l_n)$  of this recursion are determined by iterating the map
\bsplit{JPA}
	\omega &\mapsto \left( \frac{1}{\omega_2}, \frac{\omega_1}{\omega_2}\right) - (k,l) ,\\
	(k,l) &= \left( \left\lfloor \frac{1}{\omega_2}\right\rfloor, 
	 			\left\lfloor \frac{\omega_1}{\omega_2}\right\rfloor \right) .
\esplit
The algorithm is initialized with $r_{-2} = (0,1,0)$, $r_{-1} = (1,0,0)$ and $r_0 = (0,0,1)$.
We give the resulting expansions for several frequency vectors in the first four cubic fields in \Tbl{CubicFields}.
Here the sequences $(k_n,l_n)$ are always eventually periodic, with the repeated portion enclosed in $[\,]$. For example the first spiral mean case gives a period-two sequence.

\begin{table}[htp]
\begin{center}
\begin{tabular}{r| l l l l}
 $D$     & $p(x)$     &  Root                 &  $\omega$                & Jacobi-Perron  \\
 \hline
 $-23$   & $x^3-x -1$ & $1.324717957244745$& $(\sigma-1,\sigma^{-1})$ & 
 								$\left[ \tv{1}{0}, \tv{2}{0} \right]$\\
		 &            & & $(\sigma^{-1},\sigma-1)$ & $\tv{3}{2},\left[\tv{2}{0},\tv{4}{0}\right]$\\
 \hline
 $-31$   & $x^3-x^2 -1$ &$1.465571231876768$ & $(\kappa-1,\kappa^{-1})$ & $\left[\tv{1}{0}\right]$\\
         &           & & $(\kappa^{-1},\kappa-1)$  & $\tv{2}{1} , \left[ \tv{2}{0},\tv{3}{0} \right]$\\
\hline
 $-44$  & $x^3 - x^2 - x -1$ & $1.839286755214161$ & $(\tau-1,\tau^{-1})$  &$\left[\tv{1}{1} \right]$ \\
         &           & & $(\tau^{-1},\tau-1)$  & ${\tv{1}{0}}^2,\tv{3}{1}, \left[\tv{1}{0} ,
        									      {\tv{2}{0}}^3 , \tv{1}{0}, \tv{4}{0}\right]$\\
\hline
 $49$  &$x^3+x^2-2x-1$ & $1.246979603717467$ & $(\alpha-1,\alpha^{-1})$ & $\left[\tv{1}{0},\tv{3}{0}\right]$\\
         &           & & $(\alpha^{-1},\alpha-1)$ & $\tv{1}{0},\tv{2}{1},\left[\tv{1}{0},\tv{3}{0}\right]$\\
         &           & & $(\alpha-1,\alpha^2-1)$  & $\tv{4}{2} , \left[\tv{4}{0}, \tv{5}{0} \right]$\\
\end{tabular}
\end{center}
\caption{Jacobi-Perron sequences for vectors in the four cubic fields with $|D| < 50$. The last column shows the sequence generated by \Eq{JPA}. All are eventually periodic, and the periodic portions are enclosed in brackets.}
\label{tbl:CubicFields}
\end{table}

There are two cases in the table for which the sequences are period-one, analogous to
the golden mean for the ordinary continued fraction. 
The resulting sequences are then the values $(k,l)$ in the polynomial \Eq{Polynomial}
with $m=1$. In order for this to occur, the root must be a Pisot number, with $k < \tau < k+1$ and one must
choose $\omega = (\tau-k,\tau^{-1})$. 
Following Tompaidis, we might think of these vectors as generalized versions of the golden mean.

More generally it is known there are only finitely many periodic JPA sequences associated with each unit in an algebraic field \cite{Adam11}. When $d=1$, every eventually periodic continued fraction is a quadratic irrational; however, it is apparently not known how to characterize the vectors that have eventually periodic JPA sequences.

\bibliographystyle{alpha}
\bibliography{RotVector}

\newcommand{\etalchar}[1]{$^{#1}$}
\begin{thebibliography}{FWAM06}

\bibitem[AC15]{Abud15}
C.~V. Abud and I.~L. Caldas.
\newblock On {S}later's criterion for the breakup of invariant curves.
\newblock {\em Physica D}, 308:34--39, 2015.
\newblock \url{https://doi.org/10.1016/j.physd.2015.06.005}.

\bibitem[ACP06]{Altmann06}
E.G. Altmann, G.~Cristadoro, and D.~Paz.
\newblock Nontwist non-{H}amiltonian systems.
\newblock {\em Phys. Rev. E}, 73(5):056201, 2006.
\newblock \url{http://link.aps.org/abstract/PRE/v73/e056201}.

\bibitem[ACS92]{Artuso92}
R.~Artuso, G.~Casati, and D.L. Shepelyansky.
\newblock Break-up of the spiral mean torus in a volume-preserving map.
\newblock {\em Chaos, Solitons \& Fractals}, 2(2):181--190, 1992.
\newblock \url{https://doi.org/10.1016/0960-0779(92)90007-A}.

\bibitem[AR11]{Adam11}
B.~Adam and G.~Rhin.
\newblock Periodic {J}acobi-{P}erron expansions associated with a unit.
\newblock {\em J. Théorie des Nombres de Bordeaux}, 23(3):527--539, 2011.
\newblock \url{http://www.jstor.org/stable/44011250}.

\bibitem[Bak84]{Baker84}
A.~Baker.
\newblock {\em A Concise Introduction to the Theory of Numbers}.
\newblock Cambridge Univ. Press, Cambridge, 1984.

\bibitem[BM93]{Bollt93}
E.M. Bollt and J.D. Meiss.
\newblock Breakup of invariant tori for the four-dimensional semi-standard map.
\newblock {\em Physica D}, 66(3\&4):282--297, 1993.
\newblock \url{https://doi.org/10.1016/0167-2789(93)90070-H}.

\bibitem[BM94]{Baesens94a}
C.~Baesens and R.S. MacKay.
\newblock The one to two-hole transition in cantori.
\newblock {\em Physica D}, 71:372--389, 1994.
\newblock \url{https://doi.org/10.1016/0167-2789(94)90005-1}.

\bibitem[Cas57]{Cassels57}
J.W.S Cassels.
\newblock {\em An Introduction to Diophantine Approximation}.
\newblock Cambridge University Press, Cambridge, 2nd printing edition, 1957.

\bibitem[Chi79]{Chirikov79a}
B.V. Chirikov.
\newblock A universal instability of many-dimensional oscillator systems.
\newblock {\em Phys. Rep.}, 52(5):263--379, 1979.
\newblock \url{https://doi.org/10.1016/0370-1573(79)90023-1}.

\bibitem[Cla97]{Clarkson97a}
L.V. Clarkson.
\newblock {\em Approximation of Linear Forms by Lattice Points with
  Applications to Signal Processing}.
\newblock Ph{D} thesis, Australian National University, 1997.

\bibitem[CS90]{Cheng90b}
C.-Q. Cheng and Y.-S. Sun.
\newblock Existence of periodically invariant curves in 3-dimensional
  measure-preserving mappings.
\newblock {\em Celestial Mech. and Dyn. Astron.}, 47:293--303, 1990.
\newblock \url{https://doi.org/10.1007/BF00053457}.

\bibitem[Cus72]{Cusick72}
T.W. Cusick.
\newblock Formulas for some {D}iophantine approximation constants.
\newblock {\em Math. Ann.}, 197:182--188, 1972.
\newblock \url{https://doi.org/10.1007/BF01428224}.

\bibitem[Cus74]{Cusick74}
T.W. Cusick.
\newblock The two-dimensional {D}iophantine approximation constant.
\newblock {\em Monatshefte für Mathematik}, 78:297--304, 1974.
\newblock \url{https://doi.org/10.1007/BF01294641}.

\bibitem[DDS{\etalchar{+}}16]{Das16b}
S.~Das, C.B. Dock, Y.~Saiki, M.~Salgado-Flores, E.~Sander, J.~Wu, and J.A.
  Yorke.
\newblock Measuring quasiperiodicity.
\newblock {\em Euro. Phys. Lett.}, 114:40005, 2016.
\newblock \url{https://doi.org/10.1209/0295-5075/114/40005}.

\bibitem[DM12]{Dullin12a}
H.R. Dullin and J.D. Meiss.
\newblock Resonances and twist in volume-preserving mappings.
\newblock {\em SIAM J. Appl. Dyn. Sys.}, 11:319--349, 2012.
\newblock \url{https://doi.org/10.1137/110846865}.

\bibitem[DSSY16]{Das16a}
S.~Das, Y.~Saiki, E.~Sander, and J.A. Yorke.
\newblock Quasiperiodicity: Rotation numbers.
\newblock In C.~Skiadas, editor, {\em The Foundations of Chaos Revisited: From
  Poincar\'e to Recent Advancement}, Understanding Complex Systems. Springer,
  Cham, 2016.
\newblock \url{https://doi.org/10.1007/978-3-319-29701-9_7}.

\bibitem[DSSY17]{Das17}
S.~Das, Y.~Saiki, E.~Sander, and J.A. Yorke.
\newblock Quantitative quasiperiodicity.
\newblock {\em Nonlinearity}, 30(11):4111, 2017.
\newblock \url{https://doi.org/10.1088/1361-6544/aa84c2}.

\bibitem[DY18]{Das18}
S.~Das and J.A. Yorke.
\newblock Super convergence of ergodic averages for quasiperiodic orbits.
\newblock {\em Nonlinearity}, 31(2):491--501, 2018.
\newblock \url{https://doi.org/10.1088/1361-6544/aa99a0}.

\bibitem[EV01]{Efstathiou01}
K.~Efstathiou and N.~Voglis.
\newblock A method for accurate computation of the rotation and twist numbers
  for invariant tori.
\newblock {\em Physica D}, 158:151--163, 2001.
\newblock \url{https://doi.org/10.1016/S0167-2789(01)00299-8}.

\bibitem[FM13]{Fox13}
A.M. Fox and J.D. Meiss.
\newblock Greene's residue criterion for the breakup of invariant tori of
  volume-preserving maps.
\newblock {\em Physica D}, 243(1):45--63, 2013.
\newblock \url{https://doi.org/10.1016/j.physd.2012.09.005}.

\bibitem[FM16]{Fox16}
A.M. Fox and J.D. Meiss.
\newblock Computing the conjugacy of invariant tori for volume-preserving maps.
\newblock {\em SIAM J. Appl. Dyn. Sys.}, 15(1):557--579, 2016.
\newblock \url{http://epubs.siam.org/doi/abs/10.1137/15M1022859}.

\bibitem[FS73]{Froeschle73}
C.~Froeschl\'e and J.P. Scheidecker.
\newblock Numerical study of a four dimensional mapping {II}.
\newblock {\em Astron. and Astrophys.}, 22:431--436, 1973.
\newblock \url{http://adsabs.harvard.edu/abs/1973A%26A....22..431F}.

\bibitem[FWAM06]{Fuchss06}
K.~Fuchss, A.~Wurm, A.~Apte, and P.J. Morrison.
\newblock Breakup of shearless meanders and ``outer'' tori in the standard
  nontwist map.
\newblock {\em Chaos}, 16:033120, 2006.
\newblock \url{https://doi.org/10.1063/1.2338026}.

\bibitem[GMS10]{Gomez10a}
G.~G\'omez, J.M. Mondelo, and C~Sim\'o.
\newblock A collocation method for the numerical {F}ourier analysis of
  quasi-periodic functions. {I}: Numerical tests and examples.
\newblock {\em Discrete Contin. Dyn. Syst. Ser. B}, 14:41--74, 2010.
\newblock \url{https://doi.org/10.3934/dcdsb.2010.14.41}.

\bibitem[Gre79]{Greene79}
J.M. Greene.
\newblock A method for determining a stochastic transition.
\newblock {\em J. Math. Phys.}, 20:1183--1201, 1979.
\newblock \url{https://doi.org/10.1063/1.524170}.

\bibitem[HCF{\etalchar{+}}16]{Haro16}
A.~Haro, M.~Canadell, J-L Figueras, A.-L. Josep, and M.~Mondelo.
\newblock {\em The Parameterization Method for Invariant Manifolds from
  Rigorous Results to Effective Computations}.
\newblock Springer International, 2016.
\newblock \url{https://doi.org/10.1007/978-3-319-29662-3}.

\bibitem[HM88]{Hu88}
B.~Hu and J.M. Mao.
\newblock Transitions to chaos in higher dimensions.
\newblock In B.L. Hao, editor, {\em Directions in Chaos}, volume~1, pages
  206--271. World Scientific, Singapore, 1988.

\bibitem[HW79]{HardyWright79}
G.H. Hardy and E.M. Wright.
\newblock {\em An Introduction to the Theory of Numbers}.
\newblock Oxford Univ. Press, Oxford, 1979.

\bibitem[KM89]{Kook89b}
H.-t. Kook and J.D. Meiss.
\newblock Periodic-orbits for reversible, symplectic mappings.
\newblock {\em Physica D}, 35(1-2):65--86, 1989.
\newblock \url{https://doi.org/10.1016/0167-2789(89)90096-1}.

\bibitem[KO86]{Kim86}
S.~Kim and S.~Ostlund.
\newblock Simultaneous rational approximations in the study of dynamical
  systems.
\newblock {\em Phys. Rev. A}, 34:3426--3434, 1986.
\newblock \url{https://doi.org/10.1103/PhysRevA.34.3426}.

\bibitem[LFC92]{Laskar92}
J.~Laskar, C.~Froeschl\'e, and A.~Celletti.
\newblock The measure of chaos by the numerical analysis of the fundamental
  frequencies. {A}pplication to the standard mapping.
\newblock {\em Physica D}, 56:253--269, 1992.
\newblock \url{https://doi.org/10.1016/0167-2789(92)90028-L}.

\bibitem[LM10]{Levnajic10}
Z.~Levnajić and I.~Mezić.
\newblock Ergodic theory and visualization. {I}. {Mesochronic} plots for
  visualization of ergodic partition and invariant sets.
\newblock {\em Chaos}, 20(3):033114, 2010.
\newblock \url{https://doi.org/10.1063/1.3458896}.

\bibitem[Loc92]{Lochak92}
P.~Lochak.
\newblock Canonical perturbation theory via simultaneous approximation.
\newblock {\em Russ. Math. Surveys}, 47(6):59--140, 1992.
\newblock \url{https://doi.org/10.1070/RM1992v047n06ABEH000965}.

\bibitem[LV09]{Luque09}
A.~Luque and J.~Villanueva.
\newblock Numerical computation of rotation numbers of quasi-periodic planar
  curves.
\newblock {\em Physica D}, 238(20):2025--2044, 2009.
\newblock \url{https://doi.org/10.1016/j.physd.2009.07.014}.

\bibitem[LV14]{Luque14}
A.~Luque and J.~Villanueva.
\newblock Quasi-periodic frequency analysis using averaging-extrapolation
  methods.
\newblock {\em SIAM J. Dyn.Sys.}, 13(1):1--46, 2014.
\newblock \url{https://doi.org/10.1137/130920113}.

\bibitem[Mac83]{MacKay83}
R.S. MacKay.
\newblock A renormalisation approach to invariant circles in area-preserving
  maps.
\newblock {\em Physica D}, 7:283--300, 1983.
\newblock \url{https://doi.org/10.1016/0167-2789(83)90131-8}.

\bibitem[May88]{Mayer88}
D.H. Mayer.
\newblock On the distribution of recurrence times in nonlinear systems.
\newblock {\em Lett. Math. Phys.}, 16(2):139--143, 1988.
\newblock \url{https://doi.org/10.1007/BF00402021}.

\bibitem[Mei12]{Meiss12a}
J.D. Meiss.
\newblock The destruction of tori in volume-preserving maps.
\newblock {\em Comm. Nonl. Sci. Numer. Simul.}, 17:2108--2121, 2012.
\newblock \url{https://doi.org/10.1016/j.cnsns.2011.04.014}.

\bibitem[MS91]{Marmi91}
S.~Marmi and J.~Stark.
\newblock On the standard map critical function.
\newblock {\em Nonlinearity}, 5(3):743--761, 1991.
\newblock \url{https://doi.org/10.1088/0951-7715/5/3/007}.

\bibitem[MS92]{MacKay92c}
R.S. MacKay and J.~Stark.
\newblock Locally most robust circles and boundary circles for area-preserving
  maps.
\newblock {\em Nonlinearity}, 5:867--888, 1992.
\newblock \url{http://iopscience.iop.org/0951-7715/5/4/002}.

\bibitem[Sch91]{Schmidt91}
W.M. Schmidt.
\newblock {\em Diophantine Approximations and Diophantine Equations}, volume
  1467 of {\em Lect. Notes in Math.}
\newblock Springer-Verlag, New York, 1991.
\newblock \url{https://doi.org/10.1007/BFb0098246}.

\bibitem[Sie44]{Siegel44}
C.L. Siegel.
\newblock Algebraic integers whose conjugates lie in the unit circle.
\newblock {\em Duke Math. J.}, 11:597--602, 1944.
\newblock \url{https://doi.org/10.1215/S0012-7094-44-01152-X}.

\bibitem[SM20]{Sander20}
E.~Sander and J.~D. Meiss.
\newblock Birkhoff averages and rotational invariant circles for
  area-preserving maps.
\newblock {\em Physica D}, 411:132569, 2020.
\newblock \url{https://doi.org/10.1016/j.physd.2020.132569}.

\bibitem[SMS{\etalchar{+}}18]{Santos18}
M.S. Santos, M.~Mugnaine, J.D. Szezech, A.M. Batista, I.L. Caldas, M.S.
  Baptista, and R.L. Viana.
\newblock Recurrence-based analysis of barrier breakup in the standard nontwist
  map.
\newblock {\em Chaos}, 28(8):085717, 2018.
\newblock \url{https://doi.org/10.1063/1.5021544}.

\bibitem[SMS{\etalchar{+}}19]{Santos19}
M.S. Santos, M.~Mugnaine, J.D. Szezech, A.M. Batista, I.L. Caldas, and R.L.
  Viana.
\newblock Using rotation number to detect sticky orbits in {H}amiltonian
  systems.
\newblock {\em Chaos}, 29(4):043125, 2019.
\newblock \url{https://doi.org/10.1063/1.5078533}.

\bibitem[SSC{\etalchar{+}}13]{Szezech13}
J.D. Szezech, A.B. Schelin, I.L. Caldas, S~R. Lopes, P.J. Morrison, and R.L.
  Viana.
\newblock Finite-time rotation number: {A} fast indicator for chaotic dynamical
  structures.
\newblock {\em Phys. Lett. A}, 377:452--456, 2013.
\newblock \url{https://doi.org/10.1016/j.physleta.2012.12.013}.

\bibitem[Ste96]{Stewart96}
I.~Stewart.
\newblock Tales of a neglected number.
\newblock {\em Sci. Am.}, 274(6):102--103, 1996.
\newblock \url{https://doi.org/10.1038/scientificamerican0696-102}.

\bibitem[SV06]{Seara06}
T.M. Seara and J.~Villanueva.
\newblock On the numerical computation of {D}iophantine rotation numbers of
  analytic circle maps.
\newblock {\em Physica D}, 217(2):107--120, 2006.
\newblock \url{https:/doi.org/10.1016/J.Physd.2006.03.013}.

\bibitem[Tom96]{Tompaidis96b}
S.~Tompaidis.
\newblock Numerical study of invariant sets of a quasiperiodic perturbation of
  a symplectic map.
\newblock {\em Experiment. Math.}, 5:211--230, 1996.
\newblock \url{https://doi.org/10.1080/10586458.1996.10504589}.

\bibitem[VBK96]{Vrahatis96}
M.N. Vrahatis, T.C. Bountis, and M.~Kollmann.
\newblock Periodic orbits and invariant surfaces of 4d nonlinear mappings.
\newblock {\em Int. J. Bif. and Chaos}, 6(8):1425--1437, 1996.
\newblock \url{https://doi.org/10.1142/s0218127496000849}.

\bibitem[Wal12]{Waldschmidt12}
M.~Waldschmidt.
\newblock Recent advances in {D}iophantine approximation.
\newblock In D.~Goldfeld, J.~Jorgenson, J.~Jones, D.~Ramakrishnan, K.~Ribet,
  and J.~Tate, editors, {\em Number Theory, Analysis and Geometry}, pages
  659--704. Springer, Boston, MA, 2012.
\newblock \url{https://doi.org/10.1007/978-1-4614-1260-1}.

\bibitem[Xia92]{Xia92}
Z.~Xia.
\newblock Existence of invariant tori in volume-preserving diffeomorphisms.
\newblock {\em Erg. Th. Dyn. Sys.}, 12(3):621--631, 1992.
\newblock \url{https://doi.org/10.1017/S0143385700006969}.

\bibitem[ZHS01]{Zhou01b}
J.L. Zhou, B.~Hu, and Yi-Sui Sun.
\newblock Universal behaviour on the break-up of the spiral mean torus.
\newblock {\em Chin. Phys. Lett}, 18(12):1550--1553, 2001.
\newblock \url{https://doi.org/10.1088/0256-307X/18/12/303}.

\bibitem[ZTRK07]{Zou07}
Y.~Zou, M.~Thiel, M.C. Romano, and J.~Kurths.
\newblock Characterization of stickiness by means of recurrence.
\newblock {\em Chaos}, 17:043101, 2007.
\newblock \url{http://link.aip.org/link/?CHAOEH/17/043101/1}.

\end{thebibliography}

\end{document}